\newif\ifarxiv

\arxivtrue

\newif\ifthesis
\thesisfalse

\ifarxiv
	\documentclass{article}
\else
	\documentclass[a4paper,fleqn]{cas-sc}
\fi
\ifarxiv
    \usepackage[a4paper, portrait, margin=3cm]{geometry}
\fi
\usepackage{amssymb}
\usepackage{amsmath}
\usepackage[english]{babel}
\usepackage{eurosym}
\usepackage{tabularx,booktabs}
\usepackage{graphicx}
\usepackage{enumitem}
\usepackage{biblatex}
\usepackage{csquotes}
\usepackage{amsthm}
\usepackage{xspace}
\usepackage{hyperref}
\usepackage{algorithm}
\usepackage{algpseudocodex}
\usepackage{multirow}   
\usepackage{placeins}   
\usepackage{subcaption} 
\ifarxiv
    \usepackage{authblk}    
    \usepackage{orcidlink}  
\fi

\newcommand{\R}{\mathbb{R}}

\renewcommand{\vec}[1]{\mathbf{#1}}
\newcommand{\mat}[1]{\mathbf{#1}}
\newcommand{\transpose}{^\mathsf{T}}

\newcommand{\abs}[1]{\left|#1\right|}
\newcommand{\norm}[1]{\left\|#1\right\|}

\renewcommand{\O}[1]{\mathcal{O}\left(#1\right)}

\newcommand{\etal}{et~al.~}

\newtheorem  {remark}       {Remark}
\numberwithin{theorem}      {section}
\numberwithin{lemma}        {section}
\numberwithin{corollary}    {section}
\numberwithin{definition}   {section}
\numberwithin{remark}       {section}

\ifarxiv
\else
	\shorttitle{Architectural principles for two-level preconditioning}
    \shortauthors{H.A.~Melchers et~al.}
\fi

\bibliography{references/deep-learning}
\bibliography{references/additional-references}
\bibliography{references/deep-preconditioning}
\bibliography{references/greens-functions}
\bibliography{references/misc}
\bibliography{references/neural-operators}
\bibliography{references/preconditioning}

\begin{document}
\graphicspath{{figures/}}

\title{When can a neural operator replace a coarse solve? Architectural principles for two-level preconditioning}
\ifarxiv
	\author[1]{Hugo Melchers\,\orcidlink{0009-0003-4921-4253}\footnote{Corresponding author. Email: \url{mailto:h.a.melchers@tue.nl}}}
	\author[1]{Victorita Dolean\,\orcidlink{0000-0002-5885-1903}}
	\author[1]{Michael Abdelmalik\,\orcidlink{0000-0001-7373-5110}}
	\affil[1]{Eindhoven University of Technology, Eindhoven, Netherlands}
\fi

\ifarxiv
	\date{November 2026}
\fi

\maketitle

\begin{abstract}
    Neural operators are increasingly used as accelerators inside classical numerical methods, but it is rarely clear which architectural ingredients matter for which application. We answer this question for one important use case: the coarse-space correction inside a two-level preconditioner for discretised linear partial differential equations. We systematically vary four DeepONet-like architectures along two design axes: input discretisation (sampling versus integration against a basis) and source-term linearity. In doing this, we show that the favourable corner of this 2$\times$2 design is occupied by a single architecture, the Neural Green's Operator (NGO), and that moving away from it produces predictable failure modes: highly non-normal preconditioned systems, breakdown of preconditioned conjugate gradients on self-adjoint problems, and stagnation on non-self-adjoint ones. Used as a coarse-space correction, the NGO matches the iteration count of an exact coarse solve on 2D diffusion and advection-diffusion problems. The principle generalises: integrating inputs against the basis used for the output is what allows a neural operator to serve as a Galerkin-type coarse-space correction.
\end{abstract}

\FloatBarrier\section{Introduction}
\label{sec:introduction}

Solving large, sparse linear systems arising from the discretisation of partial differential equations (PDEs) is a workhorse computation across the physical sciences. Iterative Krylov solvers, such as conjugate gradients (CG) and the generalised minimal residual method (GMRES) and their variants, dominate this regime, but their convergence depends critically on the choice of preconditioner. Constructing preconditioners that are simultaneously \emph{robust} (effective across a parameter range), \emph{cheap to apply} (ideally with cost scaling linearly in the number of unknowns), and \emph{compatible} with the Krylov method of choice (e.g.\ symmetric and positive definite for CG) remains a central challenge of numerical linear algebra. Two-level preconditioners, which combine a local \emph{smoother} with a global \emph{coarse-space correction}, are an established and often highly effective approach \cite{dolean2015domain, tang2009comparison, trottenberg2001multigrid}, with multigrid the canonical instance \cite{trottenberg2001multigrid}.

Over the last decade, machine learning has emerged as a complementary route. Neural operators \cite{kovachki2023neural, li2020fourier, lu2021learning} approximate maps between function spaces and so, in principle, learn the solution operator of a parametric PDE directly. By themselves, however, neural operators may not deliver the accuracy required by downstream applications, and they offer no guarantee that the predicted field satisfies the discrete equations. A natural way to harness their speed without sacrificing reliability is to use them as \emph{preconditioners} inside a classical Krylov method: an accurate neural operator then accelerates an iterative solver of which the final accuracy is governed only by the discretisation and a user-set tolerance, not by the network \cite{rudikov2024fcgno, zhang2024blending}. This division of labour is the line of work this paper contributes to.

\paragraph*{Families of learned preconditioners}
Existing approaches to learning preconditioners can be grouped into three broad families.
\begin{itemize}
	\item \emph{Algebraic preconditioners} take the matrix \(\mat{A}\) alone and produce, via a neural network, an approximation to \(\mat{A}^{-1}\) or to a sparse factorisation thereof. Examples include direct inverse approximation \cite{ruelmann2018prospects}, learned incomplete LU/Cholesky factorisations \cite{häusner2025learning, häusner2024neural2, li2024deep, li2023learning}, and graph-neural-network preconditioners \cite{chen2025graph, trifonov2026learning}. These methods are PDE-agnostic and therefore do not exploit the structure of the underlying equation.
	\item \emph{Hybrid iterative methods} alternate a neural operator with a classical smoother inside a fixed-point iteration. The Hybrid Iterative Numerical Transferable Solver (HINTS) \cite{zhang2024blending} is the most prominent representative, with closely related schemes proposed in \cite{hu2025hybrid, kopaničáková2025deeponet}. Because the neural component is non-linear, these methods are generally combined with flexible variants of the Krylov method, e.g.\ flexible GMRES (F-GMRES) \cite{saad1993flexible} or flexible CG \cite{rudikov2024fcgno}. Recent analyses, however, show that HINTS-style methods can stagnate or diverge unless the training paradigm is chosen carefully \cite{wu2026deep}.
	\item \emph{Multigrid- and deflation-based learned preconditioners} use a learned component as a coarse-grid solver \cite{azulay2023multigrid, lerer2024multigrid, li2025neuralpreconditioning} or as a deflation/coarse subspace \cite{kopaničáková2025leveraging}. This is the family the present work belongs to. The methodology of Kopaničáková \etal\cite{kopaničáková2025leveraging} is particularly close to ours: there, a DeepONet trunk net is trained explicitly so that its basis can be repurposed as a deflation subspace inside preconditioned CG (PCG).
\end{itemize}

\noindent The methods above are increasingly competitive with classical preconditioners, but they typically introduce a learned component of which the internal structure is opaque and whose effectiveness is documented empirically rather than tied to specific architectural choices.

\paragraph*{Our approach}
A two-level preconditioner has two components: the smoother and the coarse-space correction. The smoother is well understood and largely PDE-specific; the coarse-space correction, by contrast, is exactly what a neural operator naturally approximates: a global map from a right-hand side to an approximate solution, expressed as a linear combination of basis functions. This raises the following question: \emph{which architectural properties of a neural operator govern its effectiveness as a coarse-space correction?}

To answer this question, we keep the smoother and the coarse basis \(\varphi_1, \dots, \varphi_m\) fixed and vary only the neural operator that maps right-hand sides to coarse-space coefficients. We compare four related architectures: the Deep Operator Network (DeepONet, \cite{lu2021learning}), the Resolution-Invariant Neural Operator (RINO, \cite{bahmani2025resolution}), the Variationally Mimetic Operator Network (VarMiON, \cite{patel2024variationally}), and the Neural Green's Operator (NGO, \cite{melchers2026neural}). These architectures are chosen because they represent all four options corresponding to two architectural features: first, the RINO and NGO integrate their inputs against a basis, while the DeepONet and VarMiON sample their inputs at a fixed set of sensor nodes. Second, the VarMiON and NGO preserve the linearity of the solution in the source terms, while in the DeepONet and RINO this dependence is non-linear. This design lets us attribute every difference in preconditioning behaviour to one of two architectural axes: how the input is discretised, and whether the model preserves the linear dependence of the PDE solution on the source term.

\paragraph*{Contributions}
We make four contributions.
\begin{enumerate}
	\item \textbf{An architectural principle for learned coarse corrections.} We show that the way a neural operator discretises its inputs, i.e.~by sampling versus by integration against a basis, is the dominant architectural factor controlling its effectiveness as a coarse-space correction inside a two-level preconditioner. Sampling-based architectures (DeepONet, VarMiON) produce learned operators of which the row and column spaces differ structurally, leading to spectra with non-negligible imaginary parts even on self-adjoint problems and to PCG breakdown on diffusion. Integration-based architectures (RINO, NGO) remove this asymmetry.
	\item \textbf{The role of linearity.} We show that whether a neural operator preserves the linear dependence of the PDE solution on the source term is a secondary but consequential design axis. Linearity is especially important when the spectrum of the discrete operator is structurally complex (advection-diffusion, Helmholtz), separating the NGO from the RINO.
	\item \textbf{NGO-based coarse corrections are competitive with exact coarse solves.} On the diffusion and advection-diffusion equations, NGO-based coarse corrections match exact coarse solves in iteration count, while requiring no linear solve at the coarse level. For the Helmholtz equation they remain the best of the four learned architectures but perform worse than exact coarse corrections. Additionally, the fixed-size coarse space is a fundamental limitation for high wave numbers.
	\item \textbf{NGO-based coarse corrections can outperform algebraic multigrid.} When solving diffusion equations with PCG, the two-level preconditioner with NGO-based coarse correction is computationally cheap to construct and to apply, resulting in an improvement in wall-time needed to obtain the solution when compared to algebraic multigrid (AMG) with PyAMG. However, the computational complexity of the neural operator is found to be a crucial point of improvement for NGOs with larger coarse spaces.
\end{enumerate}

\noindent We also identify two methodological conditions on which all of the above hinges: (i) some neural operators require that the training data set populates the entire coarse space, not just a subset of it (Section~\ref{sec:experiments}); and (ii) basis functions supported on Dirichlet boundaries must be excluded from the coarse correction. Both lessons are easily missed and have, to our knowledge, not been articulated in this form in prior work.

\paragraph*{Scope}
This work primarily focuses on \emph{iteration counts}, not wall-clock time. The neural operators we use are mostly evaluated on relatively small 2D problems (up to \(\sim\!10^5\) unknowns), with a single 3D test problem with up to \(5\cdot10^5\) unknowns; a fair runtime comparison against optimised classical preconditioners would be most useful for problems that are more computationally intensive, and is left to follow-up work. However, we perform a preliminary performance comparison which shows that neural operator-based coarse solvers can be more computationally efficient than their exact counterparts and make two-level preconditioners competitive with algebraic multigrid for diffusion problems. We also do not claim that an NGO is the ``best'' learned preconditioner for any of these PDEs in absolute terms; our claim is that, among neural-operator-based coarse-space corrections built from a fixed basis, the integrate-and-linear architecture systematically outperforms its sample-based and non-linear counterparts, and matches an exact coarse solve on a few test problems.

\paragraph*{Outline}
Section~\ref{sec:preliminaries} reviews two-level preconditioning and the four neural-operator architectures we test, emphasising their structural differences. Section~\ref{sec:nos-for-preconditioning} describes how each architecture is adapted to act as a coarse-space correction on a finite-difference or finite-volume discretisation. Section~\ref{sec:experiments} presents the numerical experiments on diffusion, advection-diffusion, and Helmholtz problems, including a spectral analysis that explains the observed differences. Section~\ref{sec:experiments} also includes a timing comparison (Section~\ref{subsec:vs-amg}) and tests the approach on a 3D problem (Section~\ref{subsec:3d-diffusion}). Section~\ref{sec:conclusions} draws conclusions and lists open problems.

\FloatBarrier\section{Preliminaries}
\label{sec:preliminaries}

This section recalls the two ingredients combined in the rest of the paper. Section~\ref{subsec:two-level-preconditioning} fixes notation and reviews two-level preconditioning, and Section~\ref{subsec:neural-operators} reviews the four neural-operator architectures of which the use as coarse-space corrections will be compared.

\subsection{Two-level preconditioning}
\label{subsec:two-level-preconditioning}

A two-level preconditioner combines a local correction, the \emph{smoother}, with a global correction defined on a low-dimensional subspace, the \emph{coarse-space correction}. The construction goes back to deflation \cite{nicolaides1987deflation}, balancing domain-decomposition methods \cite{mandel1993balancing}, and the multigrid coarse-grid correction \cite{trottenberg2001multigrid}. A general two-level preconditioner has the following components, following the unified notation of Tang \etal\cite{tang2009comparison}:

\begin{itemize}
	\item A large \(n \times n\) linear system \(\mat{A}\,\vec{u} = \vec{b}\), arising from the discretisation of a linear PDE.
	\item A matrix \(\mat{Z} \in \R^{n \times m}\) with \(m < n\), called the \emph{coarse basis}; its column space is the \emph{coarse space}.
	\item The \emph{coarse matrix} \(\mat{E} = \mat{Z}\transpose\,\mat{A}\,\mat{Z} \in \R^{m \times m}\).
	\item The \emph{coarse-space correction} \(\mat{Q} = \mat{Z}\,\mat{E}^{-1}\,\mat{Z}\transpose\), which approximates the action of \(\mat{A}^{-1}\) on the coarse space.
	\item A \emph{smoother} \(\mat{M}^{-1}\), typically derived from an incomplete factorisation (e.g.\ ILU or incomplete Cholesky) or from a domain-decomposition method such as additive Schwarz; see \cite[Ch.\ 2]{toselli2005domain} for a textbook treatment.
\end{itemize}

The two components \(\mat{Q}\) and \(\mat{M}^{-1}\) must then be combined to form a single preconditioner. The simplest combination is the \emph{additive} one,
\begin{equation}
	\label{eq:prec-additive}
	\mat{P}_{\text{AD}} = \mat{M}^{-1} + \mat{Q},
\end{equation}
in which the smoother and the coarse-space correction are simply summed. An alternative is to apply one component to a residual from which the contribution of the other component has already been projected out. For the two-component case relevant here, this gives the two adapted-deflation forms
\begin{subequations}
	\label{eq:prec-deflated}
	\begin{align}
		\mat{P}_{\text{A-DEF1}} & = \mat{Q} + \mat{M}^{-1}\,(\mat{I} - \mat{A}\,\mat{Q}), \label{eq:prec-adef1}      \\
		\mat{P}_{\text{A-DEF2}} & = \mat{M}^{-1} + \mat{Q}\,(\mat{I} - \mat{A}\,\mat{M}^{-1}), \label{eq:prec-adef2}
	\end{align}
\end{subequations}
which differ only in the order in which the two corrections are applied. A symmetric variant of these forms is the balancing-Neumann-Neumann preconditioner \cite{mandel1993balancing},
\begin{equation}
	\label{eq:prec-bnn}
	\mat{P}_{\text{BNN}} = \mat{Q} + (\mat{I} - \mat{Q}\,\mat{A})\,\mat{M}^{-1}\,(\mat{I} - \mat{A}\,\mat{Q}),
\end{equation}
which is symmetric whenever \(\mat{A}\), \(\mat{Q}\), and \(\mat{M}^{-1}\) are. Equivalently, a matrix-vector product of \eqref{eq:prec-adef1}--\eqref{eq:prec-bnn} with some vector \(\vec{v}\) can be implemented as two or three steps of fixed-point iteration on the un-preconditioned residual: starting from \(\vec{w}_0 = \vec{0}\), each component is applied in turn to the current residual \(\vec{v} - \mat{A}\,\vec{w}_j\) and the result is added to \(\vec{w}_j\). We use this fixed-point form in our implementation.

The four labels AD, A-DEF1, A-DEF2, and BNN are standard \cite{tang2009comparison}, and the four preconditioners \eqref{eq:prec-additive}--\eqref{eq:prec-bnn} are exactly those compared in our experiments. Note that in domain-decomposition methods the smoother \(\mat{M}^{-1}\) may itself be a combination of subdomain solvers, opening up further hybrid combinations \cite{toselli2005domain}. We take \(\mat{M}^{-1}\) as fixed throughout and do not consider such nested hybrids here.

While two-level preconditioning has been shown to be very effective, one difficulty is the matrix \(\mat{E}^{-1}\) appearing in \(\mat{Q}\). Each application of \(\mat{Q}\) requires solving an \(m \times m\) linear system \(\mat{E}\,\vec{w} = \vec{v}\). This solve does not have to be exact as long as the resulting preconditioner is still effective; in multigrid methods, for example, \(\mat{E}^{-1}\) is approximated recursively. The starting point of this paper is that an accurate neural operator can play the same role: a single pass can replace a linear solve at the coarse level.

The above definitions hold for a single linear system \(\mat{A}\vec{u} = \vec{b}\), arising from a single discretisation of a single linear PDE. In practice, we can also consider \emph{families of linear PDEs} in which the differential equation depends on some coefficients, which we abstractly denote as \(\theta\). For example, in a diffusion equation, \(\theta\) may be the diffusivity coefficient, which may be constant or may vary throughout the PDE domain. Typically, the matrix \(\mat{A}\) of the discretised equations will then depend on \(\theta\), and therefore the coarse matrix \(\mat{E}\) will also be dependent on \(\theta\). With a slight abuse of notation, we use \(\theta\) as a general symbol to refer to such PDE coefficients, as well as to the diffusion coefficient specifically. In the advection-diffusion equation, the PDE coefficients should then be understood to include both the diffusion coefficient and the velocity field \(\vec{c}\).

\subsection{Neural operators}
\label{subsec:neural-operators}

Whereas standard machine-learning models learn maps between vector spaces, \emph{neural operators} learn maps between function spaces \cite{kovachki2023neural}. They are a natural choice for approximating PDE solution operators, where the inputs and outputs, i.e.~the coefficients, source terms, boundary data, and solution itself, are all functions on the PDE domain or its boundary.

A wide variety of neural operator architectures has been developed since \cite{li2020fourier, lu2021learning}, including DeepONets and their variants \cite{bahmani2025resolution, jin2022mionet, kopaničáková2025deeponet, melchers2026neural, patel2024variationally}, Fourier Neural Operators (FNO) \cite{li2020fourier} and their extensions \cite{bonev2023spherical, li2020multipole}, transformer-based operators \cite{cao2021choose, hao2023gnot}, Convolutional Neural Operators (CNOs, \cite{raonic2024convolutional}), and resolution-invariant operators for complex geometries \cite{huang2025resolution}. These architectures differ markedly in how they handle the input fields and how they represent the output function.

This work focuses on the \emph{DeepONet-like} family: architectures of which the output is a linear combination of a finite set of basis functions. We restrict to this family for two reasons. First, this restriction enables a more direct comparison: every architecture we consider produces a function in the same finite-dimensional space, so differences in performance can be attributed to differences in the model rather than to differences in the output representation. Second, models with this property can be interpreted directly as coarse-space corrections inside a two-level preconditioner: the basis functions \(\varphi_1, \dots, \varphi_m\), evaluated at the discretisation nodes, define the coarse basis \(\mat{Z}\), and the coefficients produced by the network play the role of the coarse-space coefficients (see Section~\ref{sec:nos-for-preconditioning}).

Within this family, the four architectures we compare differ along two axes:
\begin{enumerate}
	\item how the input fields are \emph{discretised}: by sampling at sensor nodes, or by integration against a function basis;
	\item whether the architecture \emph{preserves the linear dependence} of the PDE solution on the source-term inputs.
\end{enumerate}
The four resulting combinations are realised by DeepONet, RINO, VarMiON, and NGO respectively, as summarised in Table~\ref{tab:no-design-axes}.

\begin{table}[h]
	\centering
	\caption{The four neural-operator architectures compared in this work, classified along the two design axes of which we will investigate the effect on preconditioning behaviour.}
	\label{tab:no-design-axes}
	\begin{tabular}{l|cc}
		                                       & \emph{Linear in source}               & \emph{Non-linear in source}       \\ \hline
		\emph{Samples inputs at sensor nodes}  & VarMiON \cite{patel2024variationally} & DeepONet \cite{lu2021learning}    \\
		\emph{Integrates inputs against basis} & NGO \cite{melchers2026neural}         & RINO \cite{bahmani2025resolution} \\
	\end{tabular}
\end{table}

\paragraph*{DeepONet} Introduced by Lu \etal\cite{lu2021learning} and based on the universal-approximation result of Chen and Chen \cite{chen1995universal}, the Deep Operator Network (DeepONet) is a widely used neural operator architecture. Each input field is sampled at a fixed set of \emph{sensor nodes}, and the resulting fixed-size vectors are concatenated and passed to a deep network (the \emph{branch net}) that outputs coefficients in a basis represented by a second deep network (the \emph{trunk net}). In Table~\ref{tab:no-design-axes}, the DeepONet sits in the top-right cell: it samples its inputs (top row) and depends non-linearly on all input fields, including the source term (right column).

\paragraph*{RINO} Bahmani \etal\cite{bahmani2025resolution} introduced two related architectures, the Resolution-Invariant DeepONet (RI-DeepONet) and the Resolution-Invariant Neural Operator (RINO). Both replace the sampling step of the DeepONet by projection of the input fields onto a function basis: the inputs to the branch net are basis coefficients rather than pointwise values. This makes the architecture ``resolution-invariant'', in the sense that the input fields can be specified on any set of points equipped with quadrature weights, rather than on a fixed sensor grid. The RINO additionally fixes the basis (rather than learning a trunk) and trains the model to predict basis coefficients of the solution; the RI-DeepONet keeps a learnable trunk and trains for solution accuracy. To enable a like-for-like comparison with the other architectures, we use a hybrid architecture that we will simply refer to as the RINO throughout this paper: the basis is fixed (as in the original RINO) but the model is trained to minimise the error in the solution function (as in the RI-DeepONet). In Table~\ref{tab:no-design-axes}, this RINO sits in the bottom-right cell: it integrates its inputs (bottom row) but, like the DeepONet, depends non-linearly on the source term (right column).

\paragraph*{VarMiON} Patel \etal\cite{patel2024variationally} introduce the Variationally Mimetic Operator Network (VarMiON), which is aimed specifically at \emph{linear} PDEs. The key architectural choice is to handle the PDE coefficients separately from the source terms: the coefficients enter a non-linear branch network that produces a parameter-dependent matrix, and the source-term inputs (sampled at sensor nodes) enter linearly through a matrix-vector multiplication with this matrix. The resulting model is non-linear in the PDE coefficients but \emph{linear in the source terms}, mimicking the structure of the underlying linear PDE solution operator. The VarMiON is one of several DeepONet variants that separate the dependence on different input fields; another is the Multiple-Input Operator Network (MIONet) \cite{jin2022mionet}, with closely related modifications appearing in \cite{kopaničáková2025deeponet, kopaničáková2025leveraging}. We focus on the VarMiON because it most closely matches the linear structure of the PDE problems considered here. For other (non-linear) PDEs other MIONet-like architectures may be more suitable. In Table~\ref{tab:no-design-axes}, the VarMiON sits in the top-left cell: it samples its inputs (top row) and is linear in the source terms (left column).

\paragraph*{Neural Green's Operator}
The Neural Green's Operator (NGO) is described in Melchers \etal\cite{melchers2026neural}. For our purposes, the NGO can be thought of as combining the architectural features of the RINO and the VarMiON. First, it does not sample its input fields but discretises these by integration against the basis functions. Second, it preserves the linear dependence of the output on the source terms while learning the non-linear dependence on the PDE coefficients. In Table~\ref{tab:no-design-axes}, the NGO therefore sits in the bottom-left cell. This combination is not arbitrary. The exact solution operator of a linear PDE is, by definition, the integral against a Green's function, and is linear in the source terms. The NGO is the only one of the four architectures we consider that mimics both of these structural properties. The central empirical question of this paper is whether mimicking these structural properties also makes the NGO a better coarse-space correction, and, if so, which of the two properties is responsible.

\FloatBarrier\section{Neural operators as coarse-space corrections}
\label{sec:nos-for-preconditioning}

A two-level coarse-space correction is determined by two ingredients that play very different roles. The first is the \emph{coarse basis} \(\mat{Z} \in \R^{n\times m}\), which fixes the subspace \(\mathrm{span}(\mat{Z})\) into which residuals are projected and within which the correction is sought; the choice of \(\mat{Z}\) is a structural decision about \emph{which} part of \(\mat{A}^{-1}\) the correction is meant to capture. The second is the \emph{coarse operator}
\begin{equation}
	\label{eq:coarse-operator}
	\mat{E} = \mat{Z}\transpose\,\mat{A}\,\mat{Z} \in \R^{m \times m},
\end{equation}
which is a Galerkin restriction of \(\mat{A}\) to the coarse space and which determines \emph{how} a residual is mapped to its coarse-space best approximant. Together they yield the canonical coarse-space correction \(\mat{Q} = \mat{Z}\,\mat{E}^{-1}\,\mat{Z}\transpose\), which can be read as a three-stage pipeline:
\begin{equation}
	\label{eq:coarse-pipeline}
	\underbrace{\R^n \xrightarrow{\;\mat{Z}\transpose\;} \R^m}_{\text{coarse projection}}
	\;\xrightarrow{\;\mat{E}^{-1}\;}\;
	\underbrace{\R^m \xrightarrow{\;\mat{Z}\;} \R^n}_{\text{coarse interpolation}}.
\end{equation}
The first and third stages are determined entirely by \(\mat{Z}\); only the middle stage involves a non-trivial linear solve, which is the expensive part of applying \(\mat{Q}\). It is therefore natural to ask whether the first part of $\mat{Q}$, i.e. \(\mat{E}^{-1} \mat{Z}\transpose\), can be approximated by a learned surrogate that maps coarse-space residuals to coarse-space coefficients, while leaving \(\mat{Z}\) untouched.

This separation of the coarse operator from the coarse basis is the starting point of the present paper, and it is what distinguishes the architectural question we ask from the basis-design question asked in other recent work: in Kopaničáková \etal\cite{kopaničáková2025leveraging}, for instance, the basis itself is learned (a DeepONet trunk net is trained so that its output basis functions form an effective deflation subspace), so \(\mat{Z}\) is the object being adapted and the action of \(\mat{E}^{-1} \) on the coarse space is performed exactly. We do the opposite: \(\mat{Z}\) is fixed at the start of the experiments using a problem-agnostic basis (a uniform B-spline grid), and what varies between methods is only the surrogate for \(\mat{E}^{-1} \mat{Z}\transpose\). This is what allows us to attribute every difference between the four methods we test to architectural choices in the surrogate, rather than to differences in the basis.

Each of the four neural operators considered in this work expresses its output as a linear combination of fixed basis functions \(\varphi_1, \dots, \varphi_m\), evaluated on the discretisation nodes. Identifying the matrix \(\mat{Z} = \begin{bmatrix} \boldsymbol{\varphi}_1 & \cdots & \boldsymbol{\varphi}_m \end{bmatrix}\) of nodal values with the coarse basis, a neural operator becomes a candidate substitute for the composite operator \(\mat{E}^{-1}\,\mat{Z}\transpose\) appearing in the first two stages of \eqref{eq:coarse-pipeline}: given a right-hand side, it returns a vector of coarse-space coefficients that is then mapped back to the fine grid by multiplication with \(\mat{Z}\). What the four architectures of Table~\ref{tab:no-design-axes} differ in is precisely how they implement this composite operator: whether the projection onto the coarse space is realised by sampling or by integration against the basis, and whether the dependence on the right-hand side is linear or non-linear. We make these implementations explicit in Section~\ref{subsec:four-preconditioners}.

The substitution of \(\mat{E}^{-1}\mat{Z}\transpose\) by a neural operator is attractive for two reasons. First, as is standard for preconditioners, the accuracy of the final solution depends only on the discretisation and the solver tolerance; the accuracy of the neural operator only affects the convergence rate. Second, evaluating a neural operator does not require a linear solve at the coarse level, but only a single forward pass.

\begin{remark}[Scope of the architectural axes]
	\label{rem:scope}
	The two design axes of Table~\ref{tab:no-design-axes} apply to the surrogate for \(\mat{E}^{-1}\mat{Z}\transpose\), not to the basis \(\mat{Z}\) itself. In particular, if the basis \(\mat{Z}\) is itself the wrong object, no choice of architecture for the surrogate can recover what is missing from the basis. This happens if its dimension is too small to capture the relevant near-kernel of \(\mat{A}\), as may happen for high-wave number Helmholtz (Section~\ref{subsec:helmholtz}). The architectural question we study is therefore well-posed only on top of an adequate basis.
\end{remark}

\begin{remark}[Training neural operators]
	Of particular note is that all models are trained \emph{exclusively as surrogate models}. By this, we mean that they are trained to map the PDE coefficients and right-hand-side terms to an approximation to the true PDE solution \(u\). The loss function is a simple relative squared \(L^2\) error of the approximate solution (approximated using an appropriate quadrature rule), meaning that no direct effort is made to (i) learn the full coarse operator \(\mat{E}^{-1}\), or (ii) ensure favourable spectral properties of the resulting preconditioner. This means that the neural operators are completely discretisation-agnostic: in our case, the training data is generated by solving finite-element discretisations (see Appendix~\ref{app:pde-details}), but the trained neural operators are used to precondition finite-difference or finite-volume discretisations.
\end{remark}

\subsection{Adapting a neural operator to the discrete linear system}
\label{subsec:adapting}

A neural operator is naturally a map between functions defined on the PDE domain \(\Omega\) (and its boundary \(\partial\Omega\)). A preconditioner, by contrast, is a map between vectors in \(\R^n\) whose meaning depends on the discretisation. To use a neural operator as a coarse-space correction we must therefore choose a consistent interpretation of the input vector \(\vec{v} \in \R^n\) as a function, and of the output function as a vector \(\vec{w} \in \R^n\) approximating \(\mat{A}^{-1}\,\vec{v}\).

We construct discretised systems \(\mat{A}\,\vec{u} = \vec{b}\) using second-order accurate finite-difference or finite-volume schemes. In all cases the components of \(\vec{u}\) are interpreted as nodal values of the solution \(u(\vec{x})\), and the components of \(\vec{b}\) as the corresponding nodal values of the source term \(f(\vec{x})\).\footnote{In a finite-volume method the components of \(\vec{b}\) and \(\vec{u}\) technically correspond to cell-wise averages of \(f\) and \(u\); however, these are approximated by the values at the cell barycentres to second-order accuracy, so the interpretation does not reduce the order of accuracy of the discretisation.} The full discretisations of the three PDEs considered are detailed in Appendix~\ref{app:pde-details}, and the network architectures themselves are summarised in Appendix~\ref{app:model-architecture-training}.

For each of the four architectures the output is a function expressed in the basis \(\varphi_1, \dots, \varphi_m\). Converting it to an approximation of \(\mat{A}^{-1}\,\vec{v}\) amounts to evaluating the output function at the discretisation nodes; equivalently, left-multiplying the predicted coefficients by \(\mat{Z}\). The four architectures therefore differ \emph{only} in how the input vector \(\vec{v}\) is mapped into a function input, and in how the coefficients are then computed.

The four architectures fall into two categories by the distinction in Table~\ref{tab:no-design-axes}:
\begin{itemize}
	\item \emph{Integration-based architectures (RINO, NGO).} These models take input functions only via integrals of the form \(\int_\Omega \varphi_j\,f\,\mathrm{d}x\). They have the property sometimes called ``resolution invariance'': the input function may be supplied at any set of points equipped with quadrature weights producing a quadrature rule for these integrals. In particular, no separate sampling step is needed: the components of \(\vec{v}\) are interpreted directly as the function values of \(f\) at the discretisation nodes, with the quadrature weights given by the cell measures. We collect these weights into a diagonal matrix \(\mat{W} \in \R^{n \times n}\), so that the discretised inner products of \(f\) with the basis are \(\mat{Z}\transpose\,\mat{W}\,\vec{v} \in \R^m\).
	\item \emph{Sampling-based architectures (DeepONet, VarMiON).} These models sample the input function on a fixed set of sensor nodes \(\{\vec{x}_1^{\mathrm{s}}, \dots, \vec{x}_m^{\mathrm{s}}\}\) that is, in general, distinct from the discretisation nodes. In our experiments, the two grids do not coincide, as the sensor grid is fixed at training time but the linear system can come from any discretisation. Therefore, the input function values must be interpolated from the discretisation nodes to the sensor nodes. We use nearest-neighbour interpolation throughout, and write the resulting interpolation/sampling map as a (sparse) matrix \(\mat{S} \in \R^{m \times n}\), so that the discretised input to a sampling-based model is \(\mat{S}\,\vec{v} \in \R^m\).
\end{itemize}

The PDEs we consider all have non-zero source terms as well as non-trivial boundary data, and the neural operators take all of these as separate inputs. In some cases, however, source and boundary contributions are combined into a single right-hand-side vector \(\vec{b}\). When using a neural operator as a preconditioner we therefore set the boundary inputs to zero and treat \(\vec{v}\) entirely as a discretisation of \(f\).

\subsection{The four resulting preconditioners}
\label{subsec:four-preconditioners}

With the input and output discretisations defined, each of the four neural operators induces a preconditioner of the same overall structure
\begin{equation}
	\label{eq:no-precond-general}
	\vec{w} = \mat{Z}\,\mathcal{N}(\boldsymbol{\theta};\,\vec{v}),
\end{equation}
where the matrix \(\mat{Z}\) realises the coarse-interpolation stage of the pipeline \eqref{eq:coarse-pipeline} and \(\mathcal{N}(\boldsymbol{\theta};\,\vec{v}) \in \R^m\) is a learned surrogate for the composite map \(\mat{E}^{-1}\,\mat{Z}\transpose\). Note that \(\theta\) refers to the PDE coefficients (e.g.~diffusion coefficients, velocity fields, etc.) that affect the matrices \(\mat{A}\) and \(\mat{E}\) of the full and coarse linear systems. In particular, \(\theta\) does \emph{not} refer to the trainable parameters of the neural network, which are left implicit. The vector \(\boldsymbol\theta\) is used to refer to a vector of pointwise values of the coefficients \(\theta\). The four architectures of Table~\ref{tab:no-design-axes} differ only in how \(\mathcal{N}\) is implemented: specifically, in how the coarse-projection stage \(\mat{Z}\transpose\) is realised (sampling vs.\ integration) and in whether the dependence on \(\vec{v}\) is kept linear:
\begin{itemize}
	\item \textbf{DeepONet.} The branch network takes the concatenation \([\mat{S}\,\vec{v};\,\mat{S}\,\boldsymbol{\theta}]\) of sampled source-term and coefficient values, and outputs basis coefficients directly. The map \(\vec{v} \mapsto \mathcal{N}(\boldsymbol{\theta};\,\vec{v})\) is non-linear in \(\vec{v}\).
	\item \textbf{RINO.} The branch network takes the projections of source terms and coefficients onto the basis, written here as \([\mat{G}^{-1}\,\mat{Z}\transpose\,\mat{W}\,\vec{v};\,\mat{G}^{-1}\,\mat{Z}\transpose\,\mat{W}\,\boldsymbol{\theta}]\), and outputs basis coefficients directly. Here, \(\mat{G}\) is the Gram matrix, given by \(\mat{G}_{ij} = \int_\Omega \varphi_i(\vec{x})\,\varphi_j(\vec{x})\,d\vec{x}\). Like the DeepONet, the map \(\vec{v} \mapsto \mathcal{N}(\boldsymbol{\theta};\,\vec{v})\) is non-linear in \(\vec{v}\), but the input is integrated against the basis rather than sampled.
	\item \textbf{VarMiON.} A non-linear network \(\mat{C}(\boldsymbol{\theta}) \in \R^{m \times m}\) is produced from the (sampled) PDE coefficients, and the source term enters \emph{linearly} through a matrix-vector product. The induced preconditioner therefore takes the closed form
	      \begin{equation}
		      \label{eq:varmion-preconditioner}
		      \vec{w} = \mat{Z}\,\mat{C}(\boldsymbol{\theta})\,\mat{S}\,\vec{v}.
	      \end{equation}
	\item \textbf{NGO.} As in the VarMiON, a non-linear network produces a matrix \(\mat{C}(\boldsymbol{\theta}) \in \R^{m \times m}\); the source term enters linearly, but through integration against the basis rather than through sampling. The induced preconditioner is
	      \begin{equation}
		      \label{eq:ngo-preconditioner}
		      \vec{w} = \mat{Z}\,\mat{C}(\boldsymbol{\theta})\,\mat{Z}\transpose\,\mat{W}\,\vec{v}.
	      \end{equation}
\end{itemize}

The closed forms \eqref{eq:varmion-preconditioner}--\eqref{eq:ngo-preconditioner} make the structural difference between sampling-based and integration-based architectures explicit. The NGO preconditioner has the form \(\mat{Z}\,\mat{C}\,\mat{Z}\transpose\mat{W}\): its row space is the span of the basis functions, its column space is the same span, and the only quantity learned from data is the small matrix \(\mat{C}(\boldsymbol{\theta})\). The VarMiON preconditioner has the form \(\mat{Z}\,\mat{C}\,\mat{S}\): the column space is still the span of the basis functions, but the row space is the span of (interpolated) Dirac deltas at the sensor nodes. The DeepONet and RINO have the same row spaces as the VarMiON and NGO respectively, but with the matrix \(\mat{C}\) replaced by a non-linear function of \(\vec{v}\).

This structural observation is the basis of the spectral analysis in Section~\ref{sec:experiments}: when the row and column spaces of a preconditioner do not match those of \(\mat{A}\transpose\), the preconditioned matrix \(\mat{A}\,\mat{P}\) cannot in general be made symmetric, and the spectrum acquires non-negligible imaginary components even on self-adjoint problems.

\begin{remark}
	Another important aspect of using neural operators as preconditioners is the distribution of the right-hand-side vectors \(\vec{v}\): inside right-preconditioned GMRES, for example, the preconditioner \(\mat{P}\) is applied to a set of vectors \(\vec{v}_1, \dots, \vec{v}_k\) that form an orthonormal basis for a \(k\)-dimensional Krylov space. Of these \(k\) vectors, only \(\vec{v}_1\) corresponds to the PDE problem being solved (and only in the first GMRES cycle). Inside a preconditioned solve, the neural operator may therefore be applied to right-hand sides that lie outside of the training data distribution. This is investigated further in Section~\ref{subsec:diffusion}.
\end{remark}

\paragraph*{A symmetric NGO preconditioner for conjugate gradients}
The closed form \eqref{eq:ngo-preconditioner} for the NGO preconditioner makes it possible to enforce symmetry. This is of particular interest for self-adjoint problems such as diffusion, where the canonical Krylov method is preconditioned conjugate gradients (PCG), which requires (or benefits substantially from) a symmetric positive-definite preconditioner.

When the discretisation grid is uniform, the diagonal of quadrature weights is approximately a multiple of the identity, \(\mat{W} \approx w\,\mat{I}\), with corrections only at boundary nodes. The NGO-induced preconditioner \eqref{eq:ngo-preconditioner} then approximately equals \(w\,\mat{Z}\,\mat{C}(\boldsymbol{\theta})\,\mat{Z}\transpose\); its symmetric part is
\begin{equation}
	\label{eq:ngo-sym-preconditioner}
	\mat{P}_{\mathrm{NGO,sym}} = \frac{w}{2}\,\mat{Z}\,\bigl(\mat{C}(\boldsymbol{\theta}) + \mat{C}(\boldsymbol{\theta})\transpose\bigr)\,\mat{Z}\transpose,
\end{equation}
which is symmetric by construction and positive definite whenever the symmetric part of \(\mat{C}(\boldsymbol{\theta})\) is positive definite. Section~\ref{sec:experiments} reports PCG iteration counts for both \eqref{eq:ngo-preconditioner} and \eqref{eq:ngo-sym-preconditioner}, and shows that the symmetrisation usually has only a minor effect on iteration count, but can in some cases prevent a slow-down of the convergence rate of PCG due to indefiniteness of the preconditioner. None of the other three architectures admits an analogous symmetrisation: the VarMiON has the wrong row/column structure (\eqref{eq:varmion-preconditioner} is a map from \(\R^n\) to \(\mathrm{span}(\mat{Z})\) but its transpose is not), and the DeepONet and RINO are non-linear in \(\vec{v}\) so the very notion of a symmetric part requires a linearisation that need not be unique.

\paragraph*{Reuse of trained models}
Since training a neural operator comes at a significant up-front computational cost, it is desirable to use a single trained model to precondition many linear solves. The neural operators tested here allow a high level of reuse: since the PDE coefficients \(\theta\) are passed to the neural operator, a single neural operator can be used to precondition not just a single PDE, but a family of linear PDEs parametrised by \(\theta\). Furthermore, the neural operators are trained independently of the discretisation used to obtain the linear systems \(\mat{A}\vec{u} = \vec{b}\), which means that a single neural operator can also be used to precondition different \emph{discretisations} of the same PDE. As such, a single neural operator can, once trained, be used as a preconditioner that works across (i) a range of PDE coefficients \(\theta\), (ii) a range of right-hand-side terms (leading to different vectors \(\vec{b}\)), and (iii) a range of discretisations. However, the neural operators are trained with a fixed coarse basis \(\mat{Z}\). Therefore, changing the coarse basis \(\mat{Z}\) (for example, changing the dimension \(m\)) requires re-training the neural operator.

\FloatBarrier\section{Numerical experiments}
\label{sec:experiments}

This section tests the predictions of Section~\ref{sec:nos-for-preconditioning} on three prototypical PDEs of increasing difficulty: diffusion (self-adjoint and positive definite), advection-diffusion (non-self-adjoint), and Helmholtz (indefinite). Section~\ref{sec:nos-for-preconditioning} predicts that the integration-vs-sampling axis of Table~\ref{tab:no-design-axes} may be decisive for self-adjoint problems, with the linearity axis more relevant on non-self-adjoint problems. The three test problems isolate these predictions in turn.

The general protocol is the same in all three cases. We generate a training data set by sampling and solving a large number of realisations of the PDE; train each neural operator on \(80\%\) of the data set, using another \(10\%\) as validation data to prevent overfitting (see Appendix~\ref{app:model-architecture-training}); and finally discretise the remaining \(10\%\) of the realisations and solve the resulting linear systems with various Krylov methods, comparing different choices of coarse-space correction. The full PDE definitions are in Appendix~\ref{app:pde-details}, the smoothers in Appendix~\ref{app:smoothers}, and the network architectures and training procedure in Appendix~\ref{app:model-architecture-training}. In particular, all neural operators are chosen to have approximately the same number of trainable parameters.

\paragraph*{Baselines and benchmarks}
The natural classical benchmark for these problems is multigrid \cite{trottenberg2001multigrid}, which is the gold-standard preconditioner for diffusion and remains effective for advection-diffusion in many regimes (the Helmholtz equation classically requires a different approach as will be discussed in Section~\ref{subsec:helmholtz}). Because the present work focuses on architectural principles for the learned coarse correction rather than on optimal classical methods, we use as our reference benchmark the \emph{exact} two-level preconditioner (i.e., the same combination structure with a direct solve of the coarse system \(\mat{E}\,\vec{w} = \vec{v}\)) alongside two further baselines: the smoother alone (no coarse correction), and no preconditioning at all. The target is for neural operator-based coarse corrections to match the performance of the exact coarse solves. A learned coarse correction can already be considered ``helpful'' if it results in lower iteration counts than the smoother-only preconditioning; if it leads to higher iteration counts than the smoother-only baseline, then adding the neural coarse solve is actually counterproductive.

For each PDE, the emphasis is placed on the mean iteration counts across 1000 realisations of the relevant equation, indexed by mesh size. In cases where not all linear solves converged, the effect of these divergent cases on the means is discussed. For the advection-diffusion and Helmholtz equations, we also show the dependence of the iteration counts on a relevant parameter, namely the Péclet number and wave number, respectively. Finally, we show convergence histories (relative residual norm versus iteration) for one representative test problem per PDE, which shows the qualitative convergence behaviour that the mean iteration counts smooth over.

For most test problems, we deliberately do not report wall-clock times. Our implementation is unoptimised \texttt{Python} on relatively small 2D problems, so timings would not be representative of what an optimised implementation can achieve on either the classical or the learned side of the comparison. Iteration counts are essentially implementation-independent, and therefore cleanly separate algorithmic effects from implementation effects. However, we compare the wall-clock time for the diffusion problem in Section~\ref{subsec:vs-amg}, which shows that machine learning-based coarse corrections can significantly speed up two-level preconditioners and result in performance that is competitive with classical methods.

\paragraph*{Protocol}
Unless stated otherwise, we use right preconditioning with restarted F-GMRES with restart parameter 20. The relative tolerance is set to \(\varepsilon = 10^{-6}\,(320 h)^2\), so that the finest discretisation \(h^{-1} = 320\) corresponds to a tolerance of \(10^{-6}\) and coarser meshes use a proportionally looser tolerance, matching the second-order discretisation error. Alongside these in-distribution tests, we also run a stricter set of out-of-distribution tests with tolerance \(\varepsilon = 10^{-10}\) and right-hand sides drawn from an i.i.d.\ Gaussian distribution rather than from the training data distribution. The latter tests verify that the Krylov methods do not stagnate when the right-hand side leaves the training distribution, a phenomenon documented for HINTS by Wu \etal\cite{wu2026deep}.

\begin{remark}
	We use PDE-specific smoothers throughout, which generally outperform PDE-agnostic choices such as Jacobi or Gauss-Seidel. This choice has no bearing on the conclusions of the paper: the same smoother is used with every coarse-space correction we compare, so all differences reported can be attributed to the coarse component alone. To investigate the effect of the smoother, a brief comparison to some PDE-agnostic smoothers is performed in Section~\ref{subsec:diffusion} (Table~\ref{tab:diffusion-effect-of-smoother}).
\end{remark}

\paragraph*{Coarse space}
As explained in Section~\ref{sec:nos-for-preconditioning}, the neural operators approximate PDE solutions on a fixed basis \(\varphi_i, i = 1, \dots, m\), which also determines the coarse basis \(\mat{Z}\). In two-level preconditioning, the size \(m\) of the coarse space must be chosen together with the smoother to ensure uniformly fast convergence. This is because the smoother and coarse correction essentially precondition the ``small'' and ``large'' length scales of \(\mat{A}\), respectively. To ensure that the combined preconditioner is effective, the combination of the smoother and coarse correction must together precondition all length scales. Therefore, a smaller coarse space, which preconditions only the very largest length scales, requires that the smoother grows in complexity to precondition the length scales not handled by the coarse correction.

In the experiments in this work, we choose the basis functions \(\varphi_i\) as a \(12 \times 12\) cubic B-spline basis, resulting in \(m=144\) basis functions. In Section~\ref{subsec:diffusion}, we show that the methodology we use is equally effective for larger or smaller coarse spaces. Section~\ref{subsec:vs-amg} gives insight into the trade-off in complexity of the smoother and coarse-space correction.

\subsection{Diffusion equation}
\label{subsec:diffusion}

The diffusion equation is the simplest of our three test problems: self-adjoint, positive definite, and amenable to multigrid \cite{trottenberg2001multigrid}. We use it for a number of purposes. First, as a baseline to check that learned coarse-space corrections can in principle accelerate Krylov methods. Second, to study the effect of the training data distribution on the resulting preconditioner. Third, to expose an interaction between the coarse space and the discretisation of Dirichlet boundary conditions that turns out to matter substantially. Fourth, since diffusion is normally solved with PCG, to test whether the predictions of Section~\ref{subsec:four-preconditioners} on the symmetry of the learned preconditioner (and hence its compatibility with PCG) hold up empirically. Fifth, the diffusion problem is used to investigate the generalisability of our approach to out-of-distribution PDE coefficients, different smoothers, and larger or smaller coarse spaces. Sixth, we use this test problem to study the variability of the neural operator-based preconditioners due to randomised weight initialisation and data shuffling during training.

\begin{figure}
	\centering
	\begin{subfigure}[t]{0.7\textwidth}
		\includegraphics[width=\textwidth]{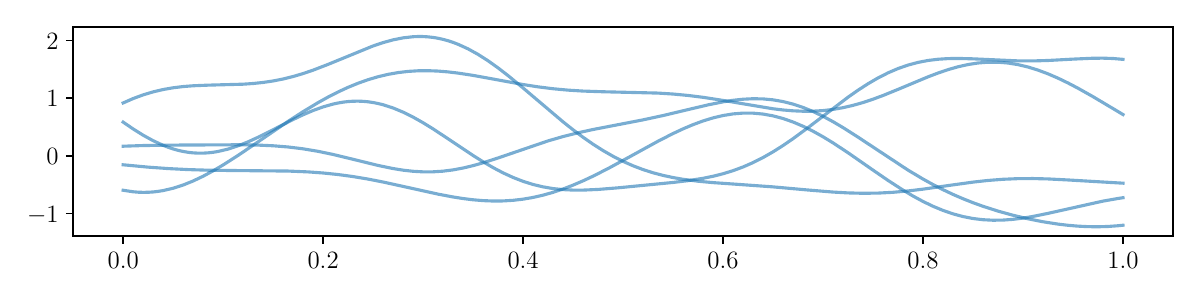}
		\caption{Dataset 1}
	\end{subfigure}
	\begin{subfigure}[t]{0.7\textwidth}
		\includegraphics[width=\textwidth]{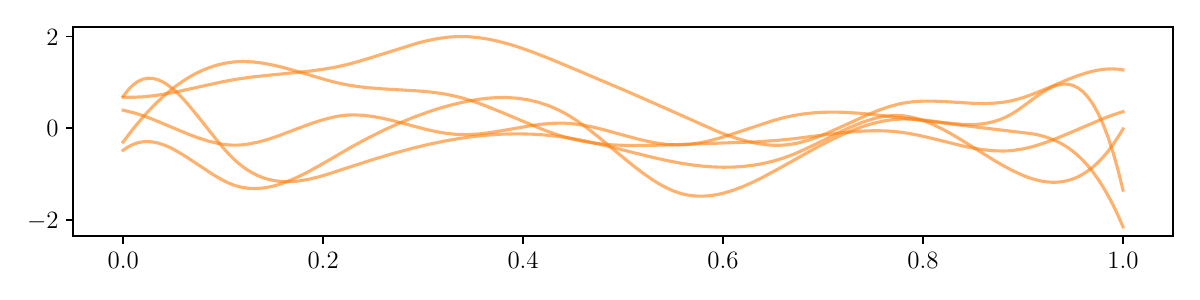}
		\caption{Dataset 2}
	\end{subfigure}
	\begin{subfigure}[t]{0.7\textwidth}
		\includegraphics[width=\textwidth]{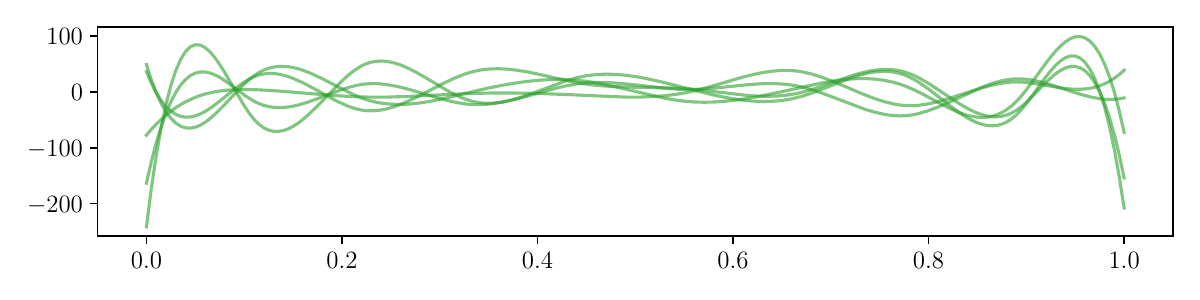}
		\caption{Dataset 3}
	\end{subfigure}
	\caption{1D source-term samples from the three diffusion training distributions. Datasets 2 and 3 both cover the full coarse space; Dataset 3 has the broadest spectral content.}
	\label{fig:diffusion-rhs-distributions}
\end{figure}

\paragraph*{Effect of the training data distribution}
A right-preconditioned Krylov method \(\mat{A}\,\mat{P}\,\vec{y} = \vec{b}\) applies the preconditioner not only to \(\vec{b}\) but to every Krylov-basis vector and every restart residual it generates. The preconditioner must therefore be a good approximation to \(\mat{A}^{-1}\) on the entire space \(\R^n\), not just for the vector \(\vec{b}\). For a coarse-space correction this means in particular that the neural operator must be accurate on \emph{every} right-hand side that lies in the coarse space, and not just on those in the training distribution. The training data must therefore populate the entire coarse space. To test this, we train all four architectures on three different distributions of source terms (see Figure~\ref{fig:diffusion-rhs-distributions}):
\begin{itemize}
	\item \emph{Dataset 1:} Gaussian Random Fields (GRFs) with squared-exponential kernel and characteristic length scale \(\lambda = 0.2\). Functions sampled from this distribution are very smooth and low-rank in the coarse basis.
	\item \emph{Dataset 2:} expansions in the same \(12 \times 12\) B-spline basis used by the neural operators, with i.i.d.\ Gaussian coefficients. By construction, this populates the full \(m\)-dimensional coarse space.
	\item \emph{Dataset 3:} as Dataset 2, but with coefficients sampled as i.i.d.\ Gaussians multiplied by the inverse Gram matrix of the basis. This distribution also covers the entire coarse space, but with much more high-frequency content.
\end{itemize}
The diffusivity \(\theta\) is drawn from the same distribution throughout (see Appendix~\ref{app:pde-diffusion} for details). The large absolute magnitude of Dataset 3 samples is harmless: the VarMiON and NGO are linear in the source by construction, and the DeepONet and RINO have positive-homogeneity rescaling enforced in the model (Appendix~\ref{app:model-architecture-training}).

\begin{figure}
	\centering
	\begin{subfigure}[t]{0.45\textwidth}
		\includegraphics[width=\textwidth]{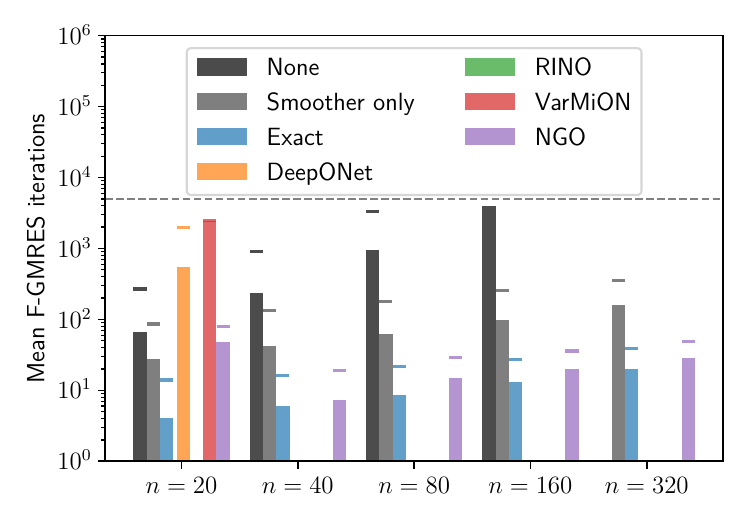}
		\caption{Dataset 1}
	\end{subfigure}
	\begin{subfigure}[t]{0.45\textwidth}
		\includegraphics[width=\textwidth]{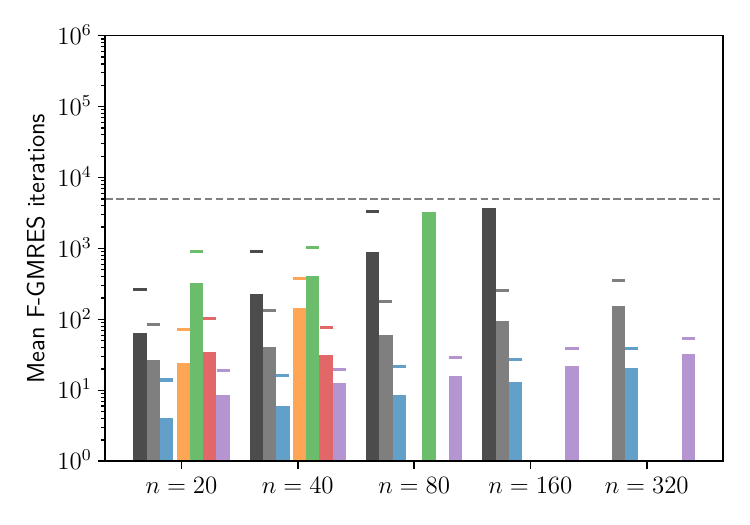}
		\caption{Dataset 2}
	\end{subfigure}
	\begin{subfigure}[t]{0.45\textwidth}
		\includegraphics[width=\textwidth]{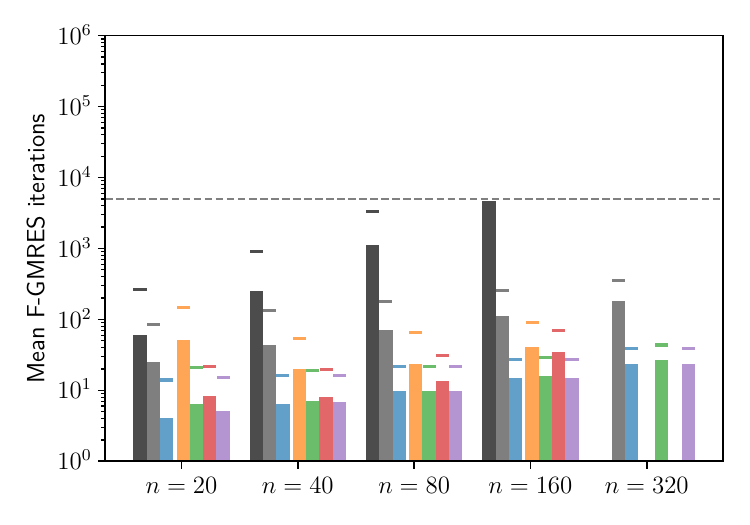}
		\caption{Dataset 3}
	\end{subfigure}
	\begin{subfigure}[t]{0.45\textwidth}
		\includegraphics[width=\textwidth]{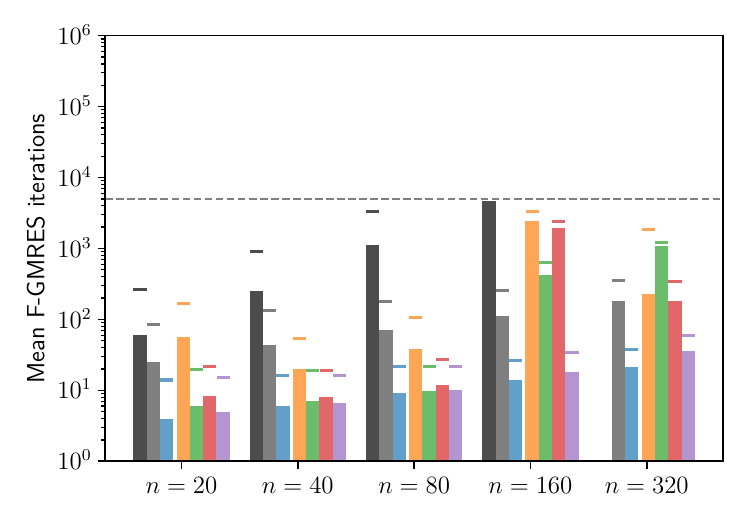}
		\caption{Dataset 3, with Dirichlet-supported basis}
	\end{subfigure}
	\caption{Mean F-GMRES iterations for diffusion problems, averaged over 1000 realisations. Bars: discretisation-dependent tolerance with the in-distribution right-hand side; lines: tolerance \(10^{-10}\) with random right-hand side. Missing entries indicate consistent convergence failure; the dashed line is the iteration limit (5000).}
	\label{fig:diff-prec-bars}
\end{figure}

\begin{figure}
	\centering
	\begin{subfigure}{0.9\textwidth}
		\includegraphics[width=\textwidth]{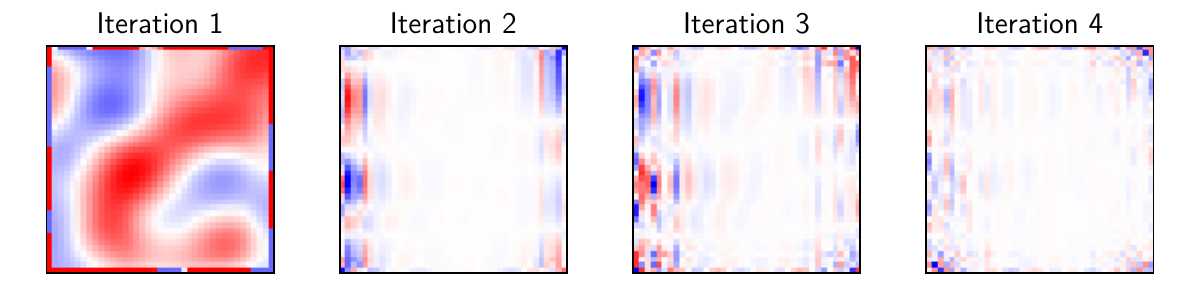}
		\caption{Dataset 1}
	\end{subfigure}
	\begin{subfigure}{0.9\textwidth}
		\includegraphics[width=\textwidth]{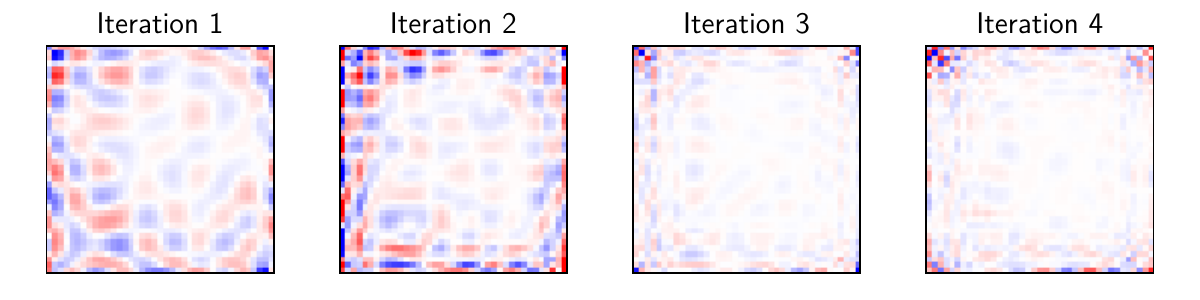}
		\caption{Dataset 3}
	\end{subfigure}
	\caption{The first four vectors passed to the exact A-DEF1 preconditioner during F-GMRES on diffusion problems from Datasets~1 and~3, reshaped to the discretisation grid. Even when the original source term is smooth (Dataset~1, top row), later GMRES iterations produce vectors of broader spectral content that the preconditioner must handle.}
	\label{fig:diff-gmres-rhss}
\end{figure}

The preconditioning results across the three datasets are shown in Figure~\ref{fig:diff-prec-bars}. The dataset has a dominant effect on preconditioner quality: most architectures fail to converge on Dataset~1, perform poorly on Dataset~2, and become reasonably effective on Dataset~3. Note that all tests are \emph{in-distribution} (each architecture is tested on samples from its training distribution); the failures on Datasets~1 and~2 are not generalisation failures but \emph{coarse-space coverage} failures. The smooth GRF samples of Dataset~1 only excite a low-rank slice of the coarse space, so the coarse-space coefficients seen during training span only a low-dimensional subspace; the resulting learned preconditioner is then a good approximation to \(\mat{A}^{-1}\) only on that subspace, not on the full coarse space the Krylov method actually probes. Figure~\ref{fig:diff-gmres-rhss} illustrates this: the first vector passed to the preconditioner during F-GMRES on a Dataset-1 problem is smooth, but later iterations supply progressively higher-frequency vectors that the model should therefore also see during training.

On Dataset~3, the integration-based architectures (RINO, NGO) outperform the sampling-based ones (DeepONet, VarMiON) consistently, and the NGO matches the exact two-level preconditioner. This is an indication that an effective coarse-space correction needs the row space of the learned operator to coincide with the column space, which only the integration-based architectures provide.

\begin{figure}
	\centering
	\includegraphics[width=0.5\textwidth]{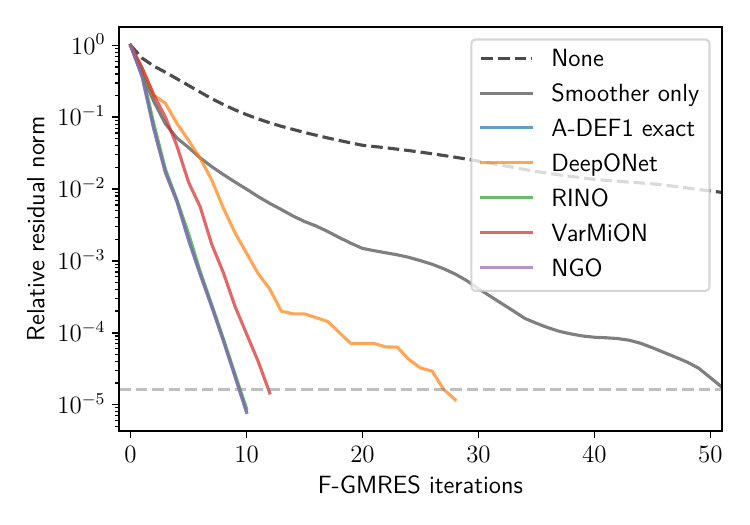}
	\caption{Relative residual norm versus F-GMRES iteration for one representative diffusion problem from Dataset~3 on the \(80 \times 80\) discretisation. The integration-based architectures (RINO, NGO) match the exact coarse-solve curve; the VarMiON converges slightly more slowly, and the DeepONet begins to converge more slowly after roughly 15 iterations.}
	\label{fig:diff-conv-history}
\end{figure}

\begin{figure}
	\centering
	\includegraphics[width=0.5\textwidth]{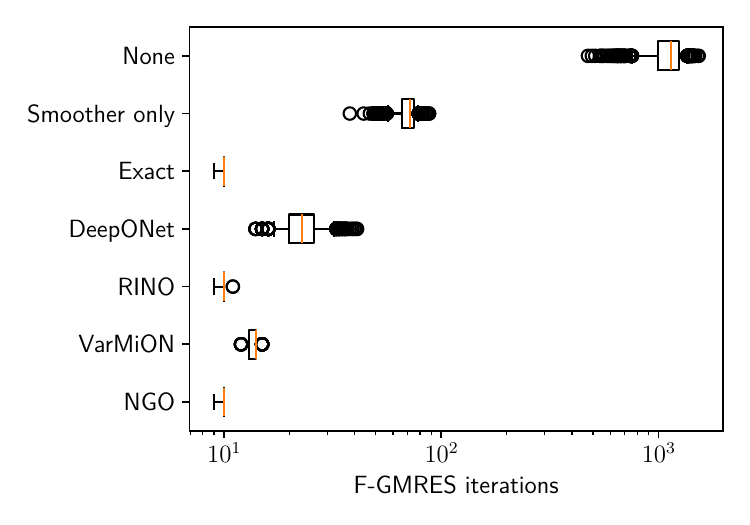}
	\caption{Distribution of F-GMRES iteration counts across the 1000 test problems on the \(80 \times 80\) Dataset-3 mesh, by preconditioner. Box plots show median, interquartile range, and 5th--95th percentiles; outliers above the 95th percentile are plotted individually.}
	\label{fig:diff-iter-distribution}
\end{figure}

The bar charts of Figure~\ref{fig:diff-prec-bars} compress two distinct effects into a single number per method: how fast a typical solve converges, and how often a solve fails. For Dataset~3, all methods converge within 5000 iterations for all 1000 problems, with the only exceptions being the \(160 \times 160\) discretisation without preconditioner and the methods omitted for the \(320 \times 320\) discretisation: un-preconditioned F-GMRES converges too slowly to reach the specified tolerance within 5000 iterations, while the DeepONet and VarMiON preconditioners cause F-GMRES to stagnate, failing to converge entirely. We refer to Appendix~\ref{app:results-table} for the precise iteration counts. Figure~\ref{fig:diff-conv-history} shows the convergence history (the relative residual norm against iteration count) for one representative problem from Dataset~3 on the \(80 \times 80\) discretisation. The exact, RINO, and NGO preconditioners produce smoothly decreasing residuals at comparable rates; the VarMiON exhibits a slower but still monotone decrease; the DeepONet begins to slow down after an initial drop. Figure~\ref{fig:diff-iter-distribution} shows the distribution of iteration counts across all 1000 test problems on the same mesh, summarised as box plots. The exact, VarMiON, RINO, and NGO distributions are tightly concentrated, while the DeepONet distribution has a heavy upper tail, reflecting a substantial fraction of test problems on which this method performs substantially worse than average.

\paragraph*{Handling Dirichlet boundaries}
A second methodological condition is summarised in the bottom-right panel of Figure~\ref{fig:diff-prec-bars}. In all the previous tests, basis functions \(\varphi_j\) supported on the Dirichlet boundary were excluded from the coarse-space correction (their coefficients forced to zero in the network output). Including these basis functions changes the picture qualitatively: at the two finest discretisations (\(160 \times 160\) and \(320 \times 320\)), the DeepONet, RINO, and VarMiON all become counterproductive: adding the coarse correction makes F-GMRES converge \emph{slower} than the smoother alone. The reason is structural. Our training data is generated using a stabilised weak form with weak imposition of Dirichlet boundary conditions (Appendix~\ref{app:pde-diffusion}), in which a Nitsche-type penalty term grows with the inverse mesh size. The corresponding coarse matrix \(\mat{E} = \mat{Z}\transpose\,\mat{A}\,\mat{Z}\) inherits eigenvalues that scale with \(h^{-1}\) along the directions of basis functions supported on \(\Gamma_\mathrm{D}\). The neural operator, however, is mesh-independent: it is trained at one discretisation but applied at any mesh size. The mismatch between the mesh-dependent coarse-matrix spectrum and the mesh-independent learned operator destroys the preconditioner's effectiveness as the mesh is refined. Removing the boundary-supported basis functions removes the unstable directions from the coarse space, restoring effectiveness. We use this ``trimmed'' coarse space throughout the remainder of the paper.

This boundary-handling condition is, to our knowledge, not articulated in prior work, and it has the convenient feature of being trivially enforceable: it does not require re-training the model, only zeroing out a subset of its output coefficients.

\begin{figure}
	\centering
	\includegraphics[width=0.9\textwidth]{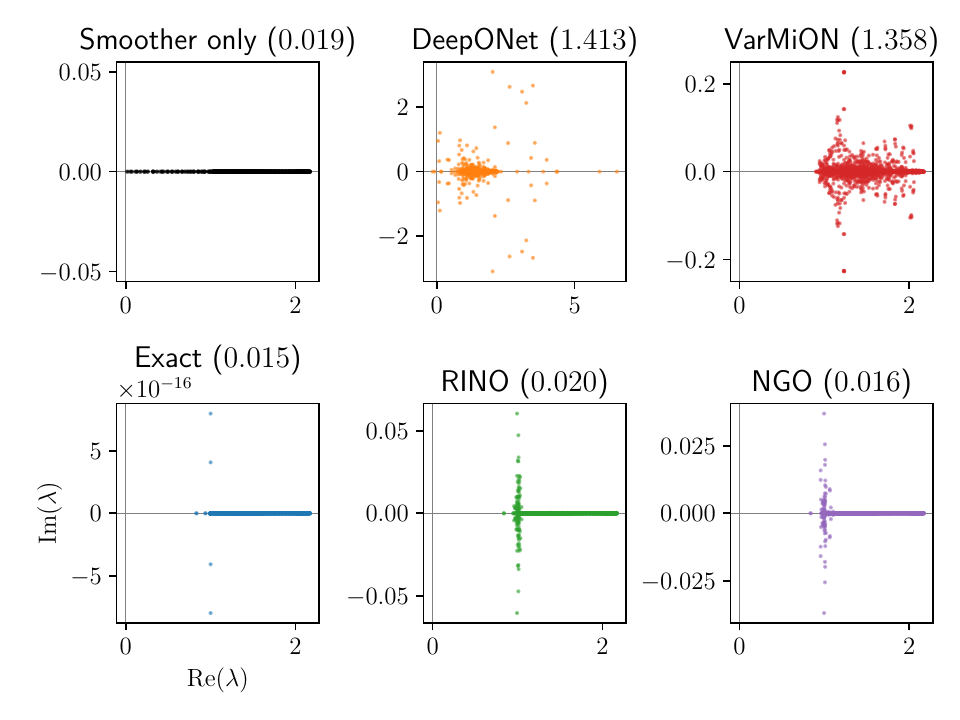}
	\caption{Eigenvalues of preconditioned systems for the six preconditioners (smoother only, exact, DeepONet, RINO, VarMiON, NGO), applied to a \(40 \times 40\) discretisation of diffusion with \(\theta \equiv 1.0\). The numbers in parentheses are the commutator errors as defined in \eqref{eq:commutator-error}. Sampling-based architectures (DeepONet, VarMiON) produce preconditioned systems that are highly non-symmetric and non-normal; integration-based ones (RINO, NGO) approximate the (real) spectrum of the exact preconditioner.}
	\label{fig:diff-prec-spectra}
\end{figure}

\paragraph*{Spectral analysis}
Figure~\ref{fig:diff-prec-spectra} shows the spectra of the four neural-operator-induced preconditioned matrices, plus the exact and smoother-only baselines, on a \(40 \times 40\) discretisation with constant diffusivity \(\theta \equiv 1\). For each preconditioner, we additionally report the relative commutator error of the preconditioned systems, defined as
\begin{align}
	c(\mat{A}) = \frac{\Vert \mat{A}\mat{A}\transpose - \mat{A}\transpose\mat{A} \Vert_F}{\Vert \mat{A}\transpose\mat{A} \Vert_F},
	\label{eq:commutator-error}
\end{align}
which is a measure of the non-normality of the matrix \(\mat{A}\). Generally, GMRES is expected to converge more quickly for linear systems in which the preconditioned matrix \(\mat{A}\) has a spectrum that is well-clustered away from the origin, and is close to normal, meaning that
\(c(\mat{A})\) is small~\cite{embree2022how}. Notice that preconditioning with only the smoother does not result in a well-conditioned system: there are still some eigenvalues that are close to zero. The addition of a coarse space is required to effectively precondition the small eigenvalues. The spectrum of exact two-level preconditioning verifies that the inclusion of a coarse correction results in a well-conditioned system. Note that the eigenvalues of the two-level preconditioned system have imaginary components of order \(10^{-16}\), which can be attributed to numerical round-off.

\begin{figure}
	\centering
	\begin{subfigure}[t]{0.25\textwidth}
		\includegraphics[width=\textwidth]{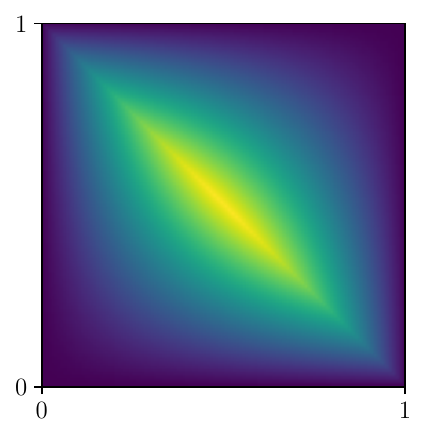}
		\caption{Exact}
	\end{subfigure}\hspace{1em}
	\begin{subfigure}[t]{0.25\textwidth}
		\includegraphics[width=\textwidth]{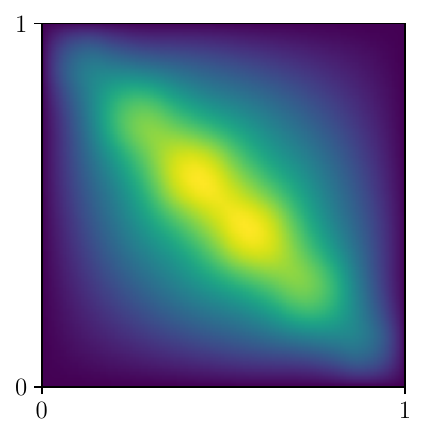}
		\caption{RINO / NGO}
	\end{subfigure}\hspace{1em}
	\begin{subfigure}[t]{0.25\textwidth}
		\includegraphics[width=\textwidth]{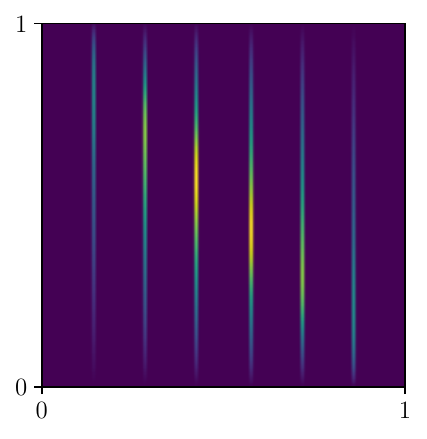}
		\caption{DeepONet / VarMiON}
	\end{subfigure}
	\caption{The Green's function of a 1D diffusion problem (left) and the structural form of the four learned preconditioners. The integration-based architectures (middle) are symmetric in their row/column structure; the sampling-based architectures (right) have a row space concentrated at the sensor nodes, breaking symmetry irrespective of training quality.}
	\label{fig:sampling-vs-integration}
\end{figure}

The spectra are exactly what Section~\ref{subsec:four-preconditioners} predicts. The DeepONet and VarMiON, i.e.~the two sampling-based architectures, produce non-normal matrices of which the eigenvalues have significant imaginary parts, even though the underlying matrix \(\mat{A}\) is symmetric. The reason is structural: the closed form \(\mat{Z}\,\mat{C}(\boldsymbol{\theta})\,\mat{S}\) of the VarMiON preconditioner (eq.~\eqref{eq:varmion-preconditioner}) maps from the span of (interpolated) Dirac deltas at the sensor nodes into the span of the basis functions, so its row and column spaces differ by construction (Figure~\ref{fig:sampling-vs-integration}). Such an operator can never be made symmetric. The DeepONet, sharing the same row-space structure (with a non-linear network in place of a fixed matrix), inherits the same property. By contrast, the integration-based architectures RINO and NGO, of which the closed-form structure is \(\mat{Z}\,\mat{C}(\boldsymbol{\theta})\,\mat{Z}\transpose\,\mat{W}\) (eq.~\eqref{eq:ngo-preconditioner}), have row and column space coinciding with \(\mathrm{span}(\mat{Z})\), and produce essentially real spectra. The remaining imaginary components of these two spectra can be attributed to training error: the trained models are generally not perfectly symmetric, resulting in eigenvalues with small but non-zero imaginary components. The eigenvalues of the RINO, VarMiON, and NGO preconditioned systems all have real parts that lie within the interval \([0.8, 2.2]\), so all three are reasonably well-conditioned.

\begin{remark}
	For general linear systems, the spectrum of a matrix does not directly determine the convergence behaviour of GMRES. In fact, a known result is that for \textit{any} spectrum and non-increasing sequence of residual norms, there exists a linear system with the given spectrum for which GMRES produces the given sequence of residual norms (see Greenbaum \etal\cite{greenbaum1996any}). As such, spectral plots such as Figure~\ref{fig:diff-prec-spectra} do not provide any guarantees regarding the convergence rate of GMRES. Nonetheless, the differences between spectra are consistent with the behaviour of GMRES for different preconditioners, meaning that these figures do offer some explanations for the observed convergence behaviour. For an overview of other matrix properties and their relation to GMRES, see Embree~\cite{embree2022how}.
\end{remark}

\begin{remark}
	The DeepONet and RINO are non-linear functions of the right-hand side, so investigating their spectra as linear operators requires a linear approximation. We linearise by applying the model to the Euclidean basis vectors \(\vec{e}_1, \vec{e}_2, \dots\) and assembling the outputs into a matrix; the non-linearity of the model is therefore not reflected in the spectral analysis. This is a necessary simplification to enable a like-for-like comparison of all four architectures.
\end{remark}

\begin{figure}
	\centering
	\begin{subfigure}[t]{0.45\textwidth}
		\includegraphics[width=\textwidth]{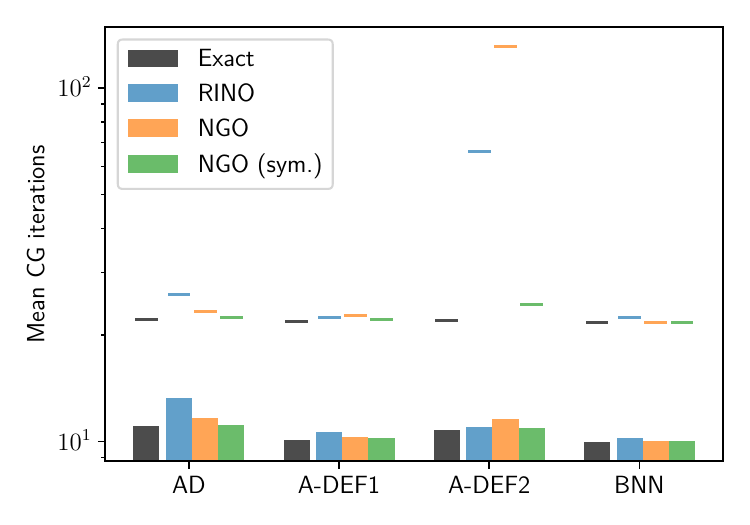}
	\end{subfigure}
	\begin{subfigure}[t]{0.45\textwidth}
		\includegraphics[width=\textwidth]{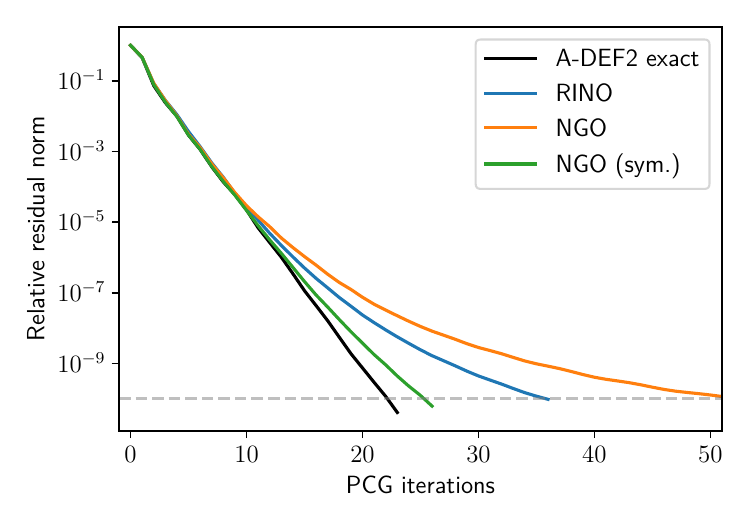}
	\end{subfigure}
	\caption{Left: mean PCG iterations for \(80 \times 80\) diffusion problems, by coarse-correction architecture and combination form. The DeepONet- and VarMiON-based corrections are omitted because PCG consistently failed on them. Right: convergence history for various A-DEF2 preconditioners on a representative problem, showing that the RINO and un-symmetrised NGO preconditioners slow down the convergence of PCG.}
	\label{fig:diff-prec-bars-cg}
\end{figure}

\paragraph*{Preconditioning with conjugate gradients}
The structural prediction of Section~\ref{subsec:four-preconditioners} has a practical consequence: only the integration-based architectures can be expected to produce preconditioners suitable for PCG. We previously derived a symmetrised NGO preconditioner \(\mat{P}_{\mathrm{NGO,sym}}\) (eq.~\eqref{eq:ngo-sym-preconditioner}) that is SPD by construction whenever the symmetric part of \(\mat{C}(\boldsymbol{\theta})\) is positive definite. Figure~\ref{fig:diff-prec-bars-cg} confirms both predictions: PCG converges with the RINO, the un-symmetrised NGO, and \(\mat{P}_{\mathrm{NGO,sym}}\), but fails to converge with the DeepONet and VarMiON in every test case (these are therefore omitted from the plot). Among the three working choices, iteration counts differ by at most a factor of two from the exact two-level preconditioner across all four combination forms (AD~\eqref{eq:prec-additive}, A-DEF1~\eqref{eq:prec-adef1}, A-DEF2~\eqref{eq:prec-adef2}, BNN~\eqref{eq:prec-bnn}). The most notable exception is A-DEF2 at low tolerance, where RINO and un-symmetrised NGO can converge several times more slowly than the symmetrised NGO and exact coarse solves. Of the four combination forms, A-DEF1 emerges as the best practical choice: it matches BNN in iteration count while requiring only one application of the coarse correction per step instead of two.

\paragraph*{Run-to-run variance}
The results above were obtained by training a single copy of each unique model. However, even for a fixed model architecture and dataset, the weights of the trained model are not deterministic due to the randomised initialisation of the model weights and the randomised order of the training data at each epoch. To test the strength of these randomised effects, we re-train four additional copies of each model, with identical architecture and training data but different initial weights.

In Table~\ref{tab:run-to-run-variance}, we show the average F-GMRES iterations needed for coarse corrections with five copies of each neural operator, all trained on the same data but with different initial weights. Note that since the averages are taken over 1000 test problems, and each problem produces an integer number of iterations, the averages are multiples of \(10^{-3}\) and are therefore listed \emph{exactly}. As can be seen in Table~\ref{tab:run-to-run-variance}, the RINO-, VarMiON-, and NGO-based coarse corrections all show very consistent performance. On the other hand, the performance of the DeepONet varies significantly per run, with one of the five models performing significantly worse than the other four. Note that the failure of the DeepONet and VarMiON on the finest mesh is consistent across all five runs.

Of all four neural operators, the NGOs appear to be the most consistent preconditioners. The largest difference between two NGOs on the same mesh size is that run 1 required \(0.06\) more iterations on average than run 3 on the \(40 \times 40\) mesh, which is less than 1\% of the average iteration count. Further analysis reveals the following:
\begin{itemize}
	\item For each mesh size, the minimum, median, and maximum iteration counts (over all 1000 test problems) are the same across all five NGOs.
	\item Over all 5000 test problems solved in this test (1000 per mesh size), 4599 require the exact same number of F-GMRES iterations independently of which of the five NGOs is used as a coarse correction.
	\item In the remaining 501 cases, all NGOs produce iteration counts that differ by only one. There is no test problem in which the weight initialisation of the NGO increased or decreased the iteration count by more than one.
\end{itemize}
From these results, it can be seen that NGOs produce very consistent preconditioners: the effect of the randomised weight initialisation is minimal. Of course, the effectiveness of the NGOs can still be affected by other parameters such as the optimiser and loss function (and indeed the dataset, as shown earlier). Nonetheless, the results in this section indicate that the training procedure used in this work produces highly consistent results.

\begin{table}
	\centering
	\caption{
		Mean F-GMRES iterations required to solve discretisations of 1000 diffusion problems on various mesh sizes. Each row corresponds to a model with differently initialised weights.
	}
	\label{tab:run-to-run-variance}
	\begin{tabular}{rr | rrrrr}
		\toprule
		& \multirow{2}*{Run} & \multicolumn{5}{c}{Mean F-GMRES iterations} \\
		&& \(n=20\) & \(n=40\) & \(n=80\) & \(n=160\) & \(n=320\) \\
		\midrule
		\multirow{5}*{\rotatebox{90}{DeepONet}}
		 & 1 &  50.675 &  19.717 & 23.429 & 41.015 &      - \\
		 & 2 &  43.091 &  18.650 & 20.767 & 35.935 &      - \\
		 & 3 &  56.896 &  26.914 & 25.828 & 48.312 &      - \\
		 & 4 & 264.660 & 238.572 & 76.098 & 53.139 &      - \\
		 & 5 &  46.638 &  17.720 & 22.447 & 42.996 &      - \\
		\midrule
		\multirow{5}*{\rotatebox{90}{RINO}}
		 & 1 &   6.326 &   7.146 &  9.933 & 15.988 & 26.379 \\
		 & 2 &   5.991 &   7.092 &  9.950 & 16.178 & 26.618 \\
		 & 3 &   6.039 &   7.115 &  9.937 & 16.120 & 26.500 \\
		 & 4 &   5.850 &   7.040 &  9.935 & 16.122 & 26.498 \\
		 & 5 &   5.807 &   7.048 &  9.933 & 16.129 & 26.546 \\
		\midrule
		\multirow{5}*{\rotatebox{90}{VarMiON}}
		 & 1 &   8.233 &   8.032 & 13.529 & 34.952 &      - \\
		 & 2 &   8.240 &   8.030 & 13.464 & 34.441 &      - \\
		 & 3 &   8.383 &   8.075 & 13.631 & 35.103 &      - \\
		 & 4 &   8.565 &   8.161 & 13.690 & 35.642 &      - \\
		 & 5 &   8.249 &   8.043 & 13.540 & 35.421 &      - \\
		\midrule
		\multirow{5}*{\rotatebox{90}{NGO}}
		 & 1 &   5.032 &   6.753 &  9.888 & 14.901 & 23.691 \\
		 & 2 &   5.029 &   6.740 &  9.881 & 14.893 & 23.657 \\
		 & 3 &   5.027 &   6.693 &  9.883 & 14.886 & 23.641 \\
		 & 4 &   5.030 &   6.745 &  9.891 & 14.896 & 23.651 \\
		 & 5 &   5.032 &   6.734 &  9.877 & 14.887 & 23.666 \\
		\bottomrule
	\end{tabular}
\end{table}

\paragraph*{Effect of smoother}
The smoother used throughout the diffusion preconditioning tests is based on the Green's function for the diffusion equation and described in Appendix~\ref{app:smoothers}. In real-world problems, many other smoothers may be used instead. To test whether the obtained results generalise to different smoothers, we compare the machine learning-based coarse corrections to exact coarse solves when used with a number of other smoothers: incomplete Cholesky decompositions, Jacobi smoothing, and block-Jacobi. The Jacobi and block-Jacobi smoothers are damped with an underrelaxation parameter of \(\omega = \frac{4}{5}\), as is the typical choice for discretisations of 2D diffusion problems (see Example~4.7.8 of Trottenberg \etal\cite{trottenberg2001multigrid}).

The results are shown in Table~\ref{tab:diffusion-effect-of-smoother}. For brevity, the results are only given for the \(160 \times 160\) mesh size. These results indicate that while different smoothers result in different convergence rates (which is expected), the NGO-based and exact coarse corrections show similar performance independently of the smoother used: for all smoothers tested, the mean iteration count when using the RINO or NGO coarse corrections is equal to or slightly higher than in the exact coarse correction case. We note that there is some variation: for example, consider the incomplete Cholesky and block-Jacobi smoothers, which perform similarly when paired with an exact coarse correction. When paired with an NGO, the average iteration count increases by only \(0.4\) for block-Jacobi, but by \(1.9\) for incomplete Cholesky. The level of ``deterioration'' in performance when the coarse solve is done inexactly differs significantly between these two smoothers. Nonetheless, the difference between exact and NGO coarse solves is small for each of the tested smoothers.

This contrasts sharply with the DeepONet- and VarMiON-based coarse corrections, which show some surprising behaviour. Most notably, the DeepONet-based coarse correction fails when used with a Jacobi or block-Jacobi smoother, and performs many times worse than other coarse corrections when combined with an incomplete Cholesky smoother. Furthermore, a combination of a block-Jacobi smoother and VarMiON coarse correction results in a higher average iteration count than block-Jacobi in isolation. Combining the VarMiON coarse correction with a Jacobi smoother also results in slightly faster convergence than a combination with block-Jacobi, despite block-Jacobi generally being a better smoother than Jacobi. Therefore, when using a VarMiON-based coarse correction, existing experience about the properties of various smoothers in combination with exact coarse corrections cannot be used to choose an appropriate smoother. The NGO, by contrast, again performs similarly to the exact coarse corrections, and uniformly outperforms the other models.

\begin{table}
	\centering
	\caption{
		Average F-GMRES iterations needed to solve \(160 \times 160\) discretisations of diffusion equations to a relative tolerance of \(4 \cdot 10^{-6}\) when using various smoothers and coarse corrections.
	}
	\label{tab:diffusion-effect-of-smoother}
	\begin{tabular}{l|rrrrrr}
		\toprule
		\multirow{2}*{Smoother} & \multicolumn{6}{c}{Coarse correction} \\
		                        & None   & Exact & DeepONet & RINO & VarMiON & NGO  \\
		\midrule
		Green's function        &  111.7 &  14.9 &     41.0 & 16.0 &    35.0 & 14.9 \\
		Incomplete Cholesky     &  325.4 &  25.6 &   1379.2 & 28.9 &    48.9 & 27.5 \\
		Jacobi                  & 2032.2 &  91.6 &        - & 97.1 &   256.0 & 93.0 \\
		Block-Jacobi            &  243.6 &  20.5 &        - & 21.1 &   263.0 & 20.9 \\
		\bottomrule
	\end{tabular}
\end{table}

\paragraph*{Preconditioner generalisation}
It has already been shown that if the training data contains a sufficiently rich distribution of source terms, then the resulting learned coarse corrections generalise to arbitrary right-hand-side vectors. However, there are a number of other ways in which a PDE problem may differ from what is contained in the training data: the \emph{domain geometry}, the \emph{boundary conditions}, and the \emph{PDE coefficients}. The better a single neural operator generalises over these features, the more broadly it can be applied and therefore the greater the benefit of training the neural operator.

The generalisation to different domain geometries is not considered in this work. This is for two reasons. First, the finite-volume and finite-difference discretisations used in this work are designed for rectangular domains and are not well-suited for discretising problems on domains with other shapes. Second, out of the four neural operators tested here the DeepONet is the only architecture for which variants exist that can be applied to varying geometries (see for example \cite{lu2022comprehensive}). While it is likely that the other model architectures can also be adapted to handle varying geometries, performing such architectural enhancements to the RINO, VarMiON, and NGO is beyond the scope of this work. However, earlier work by Kahana \etal\cite{kahana2023geometry} indicates that HINTS-style preconditioning with DeepONets \emph{can} generalise over geometries, albeit with a degradation in preconditioner effectiveness that can be mitigated by fine-tuning the DeepONet on the new geometry. Similarly, the generalisation over different boundary conditions first requires extensions to the architectures compared in this work, which is a relevant direction for future research, but is not in scope for this work.

The generalisation to different PDE parameters is tested in two different ways. First, we test the generalisation of the NGO-based coarse correction to diffusion equations with a higher \emph{contrast}, defined here as
\begin{align*}
	\text{contrast} = \frac{\sup_{\vec{x} \in \Omega} \theta(\vec{x})}{\inf_{\vec{x} \in \Omega} \theta(\vec{x})},
\end{align*}
i.e.~the ratio between the largest and smallest values that the diffusivity \(\theta\) takes on the domain. To sample diffusivity fields with a higher contrast, we first draw samples \(\theta\) from the same distribution as used in the previous data sets, but then raise them to some power, i.e.~we solve
\begin{align}
	-\nabla\cdot(\hat{\theta} \,\nabla u) & = f + \text{boundary conditions},
\end{align}
where \(\hat{\theta}(\vec{x}) = \theta(\vec{x})^\alpha\) for some \(\alpha \geq 1\).

Like most neural operators, NGOs are not expected to generalise well to PDEs with coefficient fields that lie far outside the training data distribution. As such, constructing an NGO-based preconditioner for the new field \(\hat{\theta}\) is unlikely to result in an effective preconditioner. However, we do expect that a preconditioner for one PDE remains reasonably effective for a `nearby' PDE, i.e.~a PDE with a similar diffusivity field. For this reason, we compare two strategies to precondition out-of-distribution problems:
\begin{enumerate}[label=(\alph*)]
	\item Apply the NGO to the out-of-distribution coefficient \(\hat\theta\). This tests the generalisability of the NGO itself.
	\item Apply the NGO to the in-distribution coefficient \(\theta\), and use the resulting preconditioner for the diffusion problem with \(\hat\theta\). Now, the NGO is only evaluated on in-distribution coefficients, and instead we test the generalisability of the \emph{preconditioner} to different PDEs.
\end{enumerate}
Table~\ref{tab:diffusion-powtheta} shows how both approaches perform. These results show that as the exponent \(\alpha\) is increased, the average contrast ratio of the diffusion field indeed grows. The mean iterations in the ``Smoother only'' and ``Exact'' columns indicate that two-level preconditioning remains very effective for high contrast ratios. The NGO-based coarse correction with approach (a), however, quickly deteriorates and fails entirely when the diffusivity is too far out of distribution. Instead, good results are obtained with approach (b): by applying the NGO to an in-distribution diffusivity field, we guarantee that the NGO produces an appropriate preconditioner. Although this preconditioner is now constructed for the ``wrong'' PDE, the last column of Table~\ref{tab:diffusion-powtheta} shows that the effectiveness of the preconditioner only drops very gradually as the difference between \(\theta\) and \(\hat\theta\) increases. For example, in the case of \(\alpha=5.0\), the mean contrast ratio is over \(100 \times\) greater than in the NGO's training data, but the NGO-based coarse correction is nonetheless very effective, requiring only about \(30\%\) more iterations on average than an exact coarse correction.

\begin{table}
	\centering
	\caption{Mean F-GMRES iterations needed to solve 1000 \(80\times80\) discretisations of diffusion problems, for various exponents \(\alpha\) and various preconditioners. When shown, numbers in parentheses indicate that a preconditioner resulted in convergence failure for the indicated number of equations.}
	\label{tab:diffusion-powtheta}
	\begin{tabular}{r r | r r r r}
		\toprule
		\multirow{2}*{\(\alpha\)} & \multirow{2}*{Mean contrast} & \multicolumn{4}{c}{Mean iterations}                                                    \\
		                          &                              & Smoother only                       & Exact & NGO (a)             & NGO (b)            \\
		\midrule
		1.0                       & 3.5                          & 70.9                                & 10.8  & 10.9\phantom{ (84)} & 10.9\phantom{ (1)} \\
		2.0                       & 13                           & 73.8                                & 11.0  & 11.7\phantom{ (84)} & 11.2\phantom{ (1)} \\
		3.0                       & 46                           & 78.1                                & 11.5  & 452.5 (84)          & 12.0\phantom{ (1)} \\
		4.0                       & \(1.7 \cdot 10^2\)           & 85.5                                & 12.0  & -\phantom{ (84)}    & 13.5\phantom{ (1)} \\
		5.0                       & \(6.6 \cdot 10^2\)           & 98.0                                & 12.7  & -\phantom{ (84)}    & 16.7\phantom{ (1)} \\
		6.0                       & \(2.5 \cdot 10^3\)           & 116.7                               & 13.6  & -\phantom{ (84)}    & 28.4 (1)           \\
		\bottomrule
	\end{tabular}
\end{table}

Second, we test the generalisation to diffusivity fields that oscillate more rapidly than in the training data. To do this, we generate 1000 diffusion problems where the diffusivity is sampled on a finer basis, resulting in fields that lie outside of the training distribution not due to their contrast, but due to their higher-frequency components.
\begin{figure}
	\centering
	\begin{subfigure}[t]{0.45\textwidth}
		\includegraphics[width=0.75\textwidth]{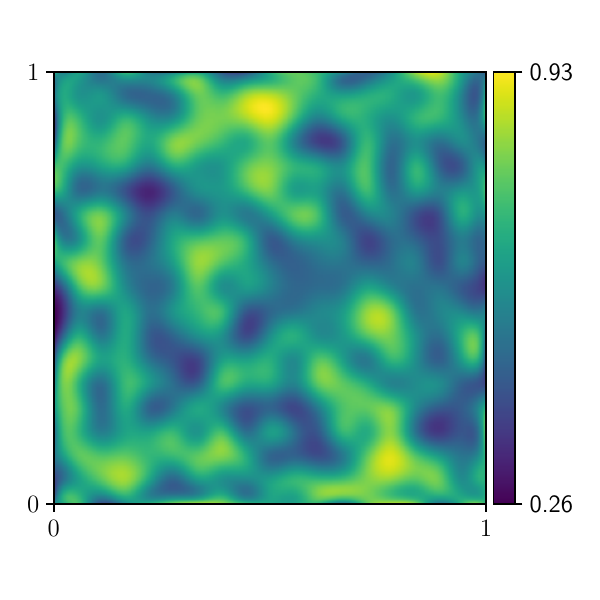}
	\end{subfigure}
	\begin{subfigure}[t]{0.45\textwidth}
		\includegraphics[width=\textwidth]{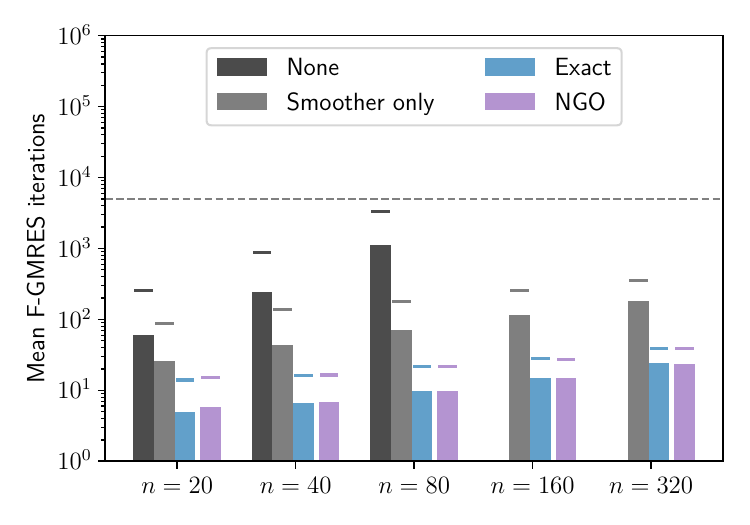}
	\end{subfigure}
	\caption{Left: an example of an out-of-distribution diffusivity field \(\theta\) with finer features. Right: the NGO preconditioner effectiveness on 1000 of such out-of-distribution problems.}
	\label{fig:diffusion-fine-theta}
\end{figure}
Figure~\ref{fig:diffusion-fine-theta} shows an example out-of-distribution diffusivity field, as well as the performance of the NGO-based coarse correction compared to the exact and smoother-only preconditioners. Comparing with Figure~\ref{fig:diff-prec-bars}, we see that the preconditioner performance is essentially identical to the in-distribution case.

\paragraph*{Effect of coarse-space dimension}
So far, all coarse corrections have used the same \(12 \times 12\) coarse space. In practical applications, however, the dimension of the coarse space typically grows with the size of the discretisation. This is done to improve computational efficiency: a larger coarse space preconditions more length scales of the problem, meaning that the other part of the preconditioner, the smoother, has to handle fewer length scales and therefore becomes more computationally efficient. Obtaining good computational performance for a two-level preconditioner therefore requires balancing the computational complexities of the coarse correction and the smoother.

In addition to the previously trained NGO, we also trained an NGO with a smaller \(6 \times 6\) cubic B-spline basis, as well as one with a larger \(24 \times 24\) basis. The results in this section come from three different NGOs: each one is used for a fixed value of the coarse-space size \(m\), and is then reused across all values of \(n\) as well as the PDE coefficients \(\theta\) and right-hand-side vectors \(\vec{b}\). Using either the smaller or larger NGO as a coarse correction, we then also double or halve the length scale \(H\) of the smoother, respectively (see Appendix~\ref{app:smoothers}). Table~\ref{tab:diffusion-coarse-space} shows the effectiveness of the resulting preconditioners. These results show that NGO-based coarse corrections remain about equally effective as exact coarse corrections. The only exception to this is the \(24 \times 24\) coarse correction applied to the \(20 \times 20\) discretisation: in this case, the dimension of the coarse space exceeds that of the discretisation being solved. As a result, the coarse-correction is actually the exact inverse: \(\mat{Q} = \mat{A}^{-1}\). Therefore, the Krylov method simply converges in one iteration. Of course, in practical applications the coarse space would not be chosen to be larger than the original system, so this data point is of little practical relevance, and is only included here for the sake of completeness.

Importantly, the results of Table~\ref{tab:diffusion-coarse-space} indicate that the preconditioners can be made effective for a wide range of coarse spaces, without requiring careful tuning of model hyper-parameters. This implies that the same methodology could be applied to even larger coarse spaces, producing effective preconditioners for very large linear systems. We note, however, that effectively preconditioning with larger coarse spaces also requires that the computational complexity of the NGO scales favourably in the dimension of the coarse space. We return to this point in Section~\ref{subsec:vs-amg}.

\begin{table}
	\centering
	\caption{Mean F-GMRES iterations over 1000 diffusion problems, for various discretisations, preconditioners, and coarse-space dimensions.}
	\label{tab:diffusion-coarse-space}
	\begin{tabular}{r r l | r r r r r}
		\toprule
		\multirow{2}*{\(m\)}    & \multirow{2}*{\(H\)}     & \multirow{2}*{Preconditioner} & \multicolumn{5}{c}{Mean F-GMRES iterations}                                                       \\
		                        &                          &                               & \(n = 20\)                                  & \(n = 40\) & \(n = 80\) & \(n = 160\) & \(n = 320\) \\
		\midrule
		\multirow{3}*{\(6^2\)}  & \multirow{3}*{\(0.1\)}   &
		Smoother only           & 17.0                     & 29.0                          & 45.8                                        & 74.2       & 123.3                                  \\
		                        &                          & Exact A-DEF1                  & 6.0                                         & 9.0        & 14.3       & 23.1        & 35.5        \\
		                        &                          & NGO A-DEF1                    & 6.0                                         & 9.0        & 14.2       & 22.9        & 35.3        \\
		\midrule
		\multirow{3}*{\(12^2\)} & \multirow{3}*{\(0.05\)}  &
		Smoother only           & 25.4                     & 43.7                          & 69.9                                        & 111.7      & 181.0                                  \\
		                        &                          & Exact A-DEF1                  & 4.1                                         & 6.4        & 9.8        & 14.9        & 23.8        \\
		                        &                          & NGO A-DEF1                    & 5.0                                         & 6.8        & 9.9        & 14.9        & 23.7        \\
		\midrule
		\multirow{3}*{\(24^2\)} & \multirow{3}*{\(0.025\)} &
		Smoother only           & 25.4                     & 65.2                          & 128.8                                       & 207.7      & 310.5                                  \\
		                        &                          & Exact A-DEF1                  & 1.0                                         & 4.1        & 6.1        & 9.5         & 14.5        \\
		                        &                          & NGO A-DEF1                    & 6.7                                         & 6.0        & 7.0        & 9.8         & 14.4        \\
		\bottomrule
	\end{tabular}
\end{table}

\FloatBarrier
\subsection{Advection-diffusion equation}
\label{subsec:advection-diffusion}

The advection-diffusion equation is non-self-adjoint: its discrete operator \(\mat{A}\) is no longer symmetric, even for symmetric domains and constant coefficients. The integration-vs-sampling axis is therefore expected to matter less than for diffusion, since there is no row/column-space symmetry of \(\mat{A}\) for an integration-based architecture to preserve. This means that the second architectural property, the linear dependence of the solution on the source, should now have a more visible effect. We use the same B-spline source distribution as Dataset~3 of the diffusion experiments (the only one effective for training). To test out-of-distribution generalisation we additionally construct a high-Péclet test set with substantially larger advective velocities than those appearing in training. The PDE definitions, the four-direction Gauss-Seidel smoother, and the data-set construction are detailed in Appendix~\ref{app:pde-advdiff}.

\begin{figure}
	\centering
	\begin{subfigure}[t]{0.45\textwidth}
		\includegraphics[width=\textwidth]{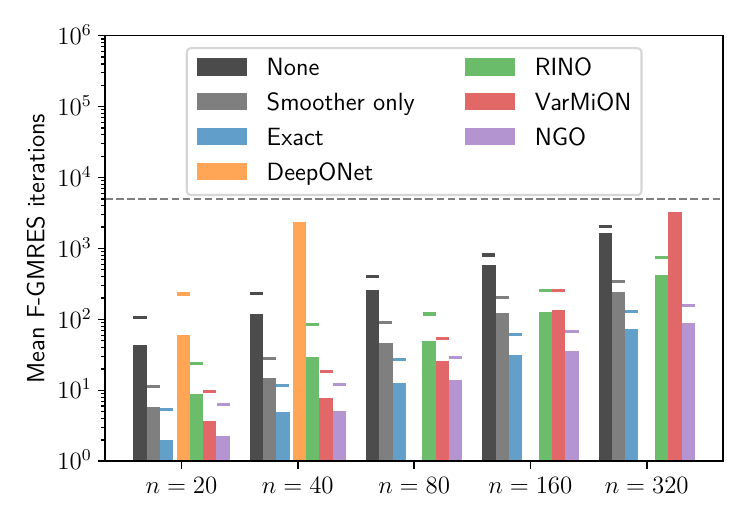}
		\caption{Low Péclet dataset (in-distribution)}
	\end{subfigure}
	\begin{subfigure}[t]{0.45\textwidth}
		\includegraphics[width=\textwidth]{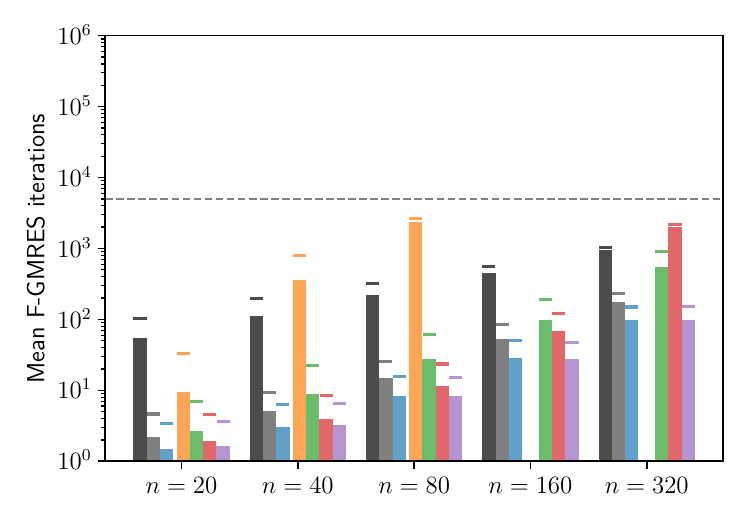}
		\caption{High Péclet dataset (out-of-distribution)}
	\end{subfigure}
	\caption{Mean F-GMRES iterations for advection-diffusion. Bars: discretisation-dependent tolerance, in-distribution right-hand side; lines: tolerance \(10^{-10}\), random right-hand side. The NGO closely tracks the exact two-level preconditioner on every mesh and on both data sets.}
	\label{fig:advdiff-prec-bars}
\end{figure}

The results are summarised in Figure~\ref{fig:advdiff-prec-bars}. Three behaviours are visible. First, the four-direction Gauss-Seidel smoother is already a strong preconditioner on its own, especially for high Péclet numbers. Second, the NGO-based correction closely tracks the exact coarse solve across all five mesh resolutions, both in- and out-of-distribution. Third, only the NGO produces a consistent speed-up over the smoother. The VarMiON helps on coarse meshes but is harmful elsewhere; the RINO is consistently counterproductive; the DeepONet typically fails to converge altogether. The advection-diffusion problem thus separates the four architectures along both design axes simultaneously. The NGO is the only consistently effective choice, with the second-best option being either the RINO or VarMiON depending on the mesh size, and the DeepONet consistently failing.

\begin{figure}
	\centering
	\includegraphics[width=0.5\textwidth]{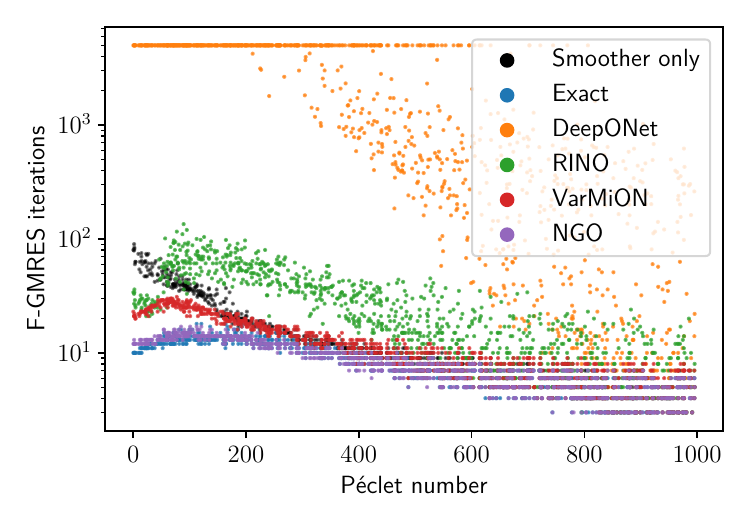}
	\caption{F-GMRES iterations for advection-diffusion as a function of the average Péclet number, on the \(80 \times 80\) discretisation. The NGO is the only architecture that consistently accelerates the smoother across the full Péclet range.}
	\label{fig:advdiff-prec-vs-pe}
\end{figure}

Figure~\ref{fig:advdiff-prec-vs-pe} shows the same data resolved by Péclet number on the \(80 \times 80\) discretisation. The four-direction Gauss-Seidel smoother becomes increasingly effective on its own as the Péclet number grows, so the gain from a coarse correction is largest in the diffusion-dominated regime. There, all three models provide a speed-up, with the NGO the most effective; the RINO, however, becomes counterproductive at high Péclet, and the VarMiON shows the same behaviour to a lesser extent. This leaves the NGO as the only choice that is consistently effective across the tested range of Péclet numbers.

\begin{figure}
	\centering
	\begin{subfigure}[t]{0.45\textwidth}
		\includegraphics[width=\textwidth]{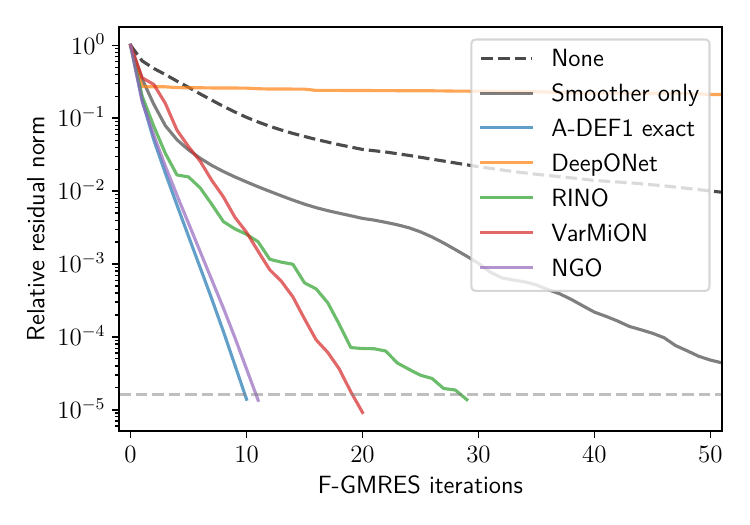}
		\caption{Diffusion-dominated, \(\mathrm{Pe} = 10\)}
	\end{subfigure}
	\begin{subfigure}[t]{0.45\textwidth}
		\includegraphics[width=\textwidth]{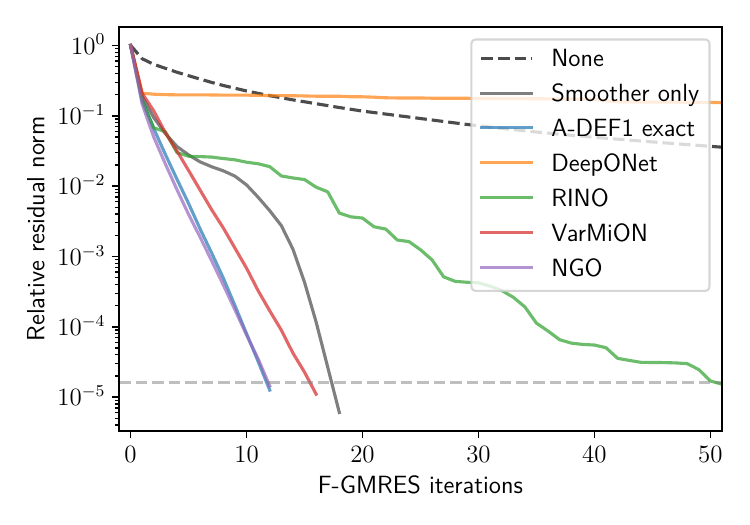}
		\caption{Advection-dominated, \(\mathrm{Pe} = 200\)}
	\end{subfigure}
	\caption{Convergence histories of F-GMRES on representative advection-diffusion problems at low and high Péclet number.}
	\label{fig:advdiff-conv-history}
\end{figure}

Figure~\ref{fig:advdiff-conv-history} confirms the per-iteration behaviour expected from the bar chart: in the diffusion-dominated regime (\(\mathrm{Pe} = 10\)), the NGO and RINO drive the residual down at a rate close to the exact two-level preconditioner, while the DeepONet stagnates after an initial drop. In the advection-dominated regime (\(\mathrm{Pe} = 200\)), the gap between the smoother and the smoother-plus-coarse-correction narrows substantially, but the NGO still matches the exact coarse solve. The RINO trace shows a plateau near iteration 15 followed by re-convergence, though at a slower rate.

Table~\ref{tab:advdiff-failures} shows the number of failures for each preconditioner on the \(80 \times 80\) discretisations, as well as the average number of F-GMRES iterations over the cases where convergence was achieved. The DeepONet-based preconditioner, which is the only one to lead to failures for these problems, also performs the worst on average when convergence failures are omitted from the average. As such, the presence of convergence failures does not significantly distort the comparison shown in Figure~\ref{fig:advdiff-prec-bars}.

\begin{table}
	\centering
	\caption{Average number of F-GMRES iterations and number of failures for different preconditioners on \(80 \times 80\) discretisations of the advection-diffusion equation. The average iteration counts now \textit{exclude} the cases where convergence was not achieved.}
	\label{tab:advdiff-failures}
	\begin{tabular}{l r r r r}
		\toprule
		\multirow{2}*{Preconditioner} & \multicolumn{2}{c}{Low Péclet number} & \multicolumn{2}{c}{High Péclet number}                              \\
		                              & Mean iterations                       & Failures                               & Mean iterations & Failures \\
		\midrule
		None                          & 258.2                                 & 0                                      & 218.5           & 0        \\
		Smoother only                 & 47.0                                  & 0                                      & 14.7            & 0        \\
		Exact                         & 12.8                                  & 0                                      & 8.2             & 0        \\
		DeepONet                      & -                                     & 1000                                   & 465.8           & 416      \\
		RINO                          & 50.3                                  & 0                                      & 27.4            & 0        \\
		VarMiON                       & 25.8                                  & 0                                      & 11.7            & 0        \\
		NGO                           & 13.9                                  & 0                                      & 8.3             & 0        \\
		\bottomrule
	\end{tabular}
\end{table}

\begin{figure}
	\centering
	\includegraphics[width=0.9\textwidth]{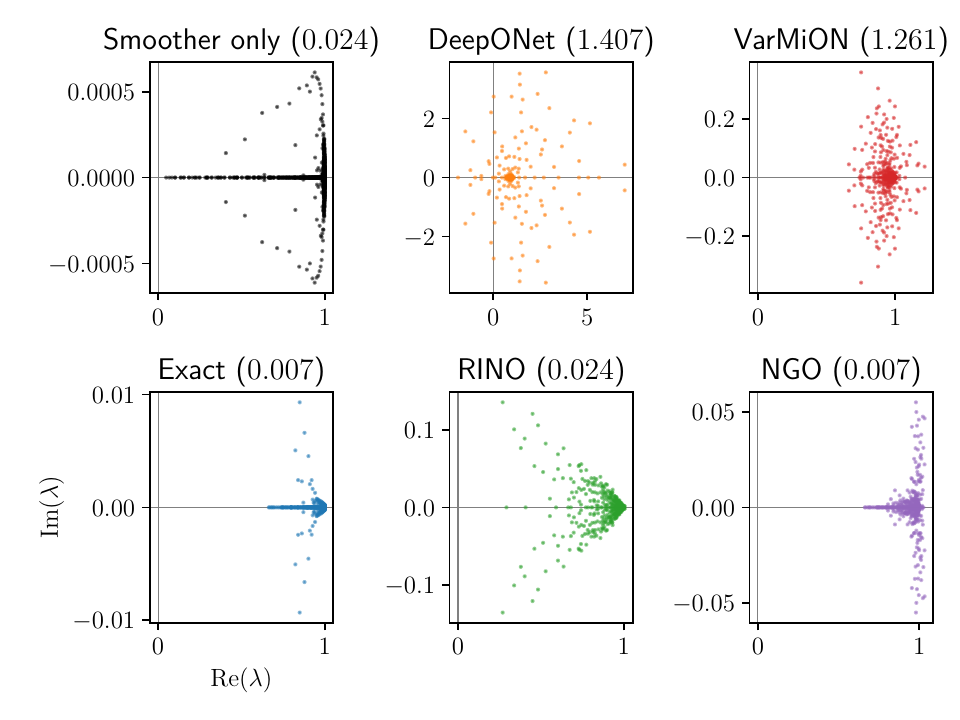}
	\caption{Spectra of preconditioned systems for advection-diffusion on a \(40 \times 40\) discretisation with \(\theta \equiv 0.1\) and \(\vec{c} = (0.6, 0.8)\) (\(\mathrm{Pe} = 10\)). The numbers in parentheses are the commutator errors as defined in \eqref{eq:commutator-error}. The NGO and RINO both produce well-clustered spectra; the NGO's is tighter.}
	\label{fig:advdiff-spectra}
\end{figure}

The spectral plot in Figure~\ref{fig:advdiff-spectra} explains the convergence behaviour. Complex eigenvalues are now expected, as \(\mat{A}\) itself is non-symmetric, so the question is no longer whether the spectrum is real but whether it is well-clustered around \(1\). The DeepONet preconditioner produces eigenvalues spread over a region containing both negative and large-imaginary values, which hinders convergence. The VarMiON spectrum is closer to the real axis but has a wide spread of moduli. The RINO and NGO both produce spectra concentrated near \(1\); the NGO's is tighter, consistent with its iteration-count advantage.

\begin{remark}
	When constructing neural operator-based coarse corrections for advection-diffusion equations with high Péclet number, the inputs to the neural operator are also regularised: the diffusion field \(\theta\) is scaled to ensure \(\text{Pe} \leq 60\). Importantly, this scaling is not applied to the linear system being solved, or in the construction of the exact coarse correction. It is purely done to avoid applying the neural operators to problems that lie far outside the training data distribution, and Figure~\ref{fig:advdiff-prec-vs-pe} shows that this does not negatively affect the preconditioner effectiveness for high Péclet numbers.
\end{remark}

\paragraph*{Sensitivity to training data set size}
The models evaluated in this work were trained on 8000 example solutions of the relevant PDE. From Figure~\ref{fig:advdiff-prec-bars}, it can be seen that for the advection-diffusion problem, this is sufficient for the NGO-based coarse corrections nearly to match their exact counterparts. Since generating training data sets may be a computationally expensive process, it is useful to know how large the training data set needs to be in order to be able to train accurate preconditioners. It turns out that for this family of advection-diffusion problems, far fewer than 8000 example solutions are already sufficient. To test this, we train NGOs on subsets of the 8000 solutions of the training data, and evaluate them on the same 1000 test equations. Both training and evaluation are done on the low-Péclet number data set, as low-Péclet number problems were found to be most sensitive to the coarse correction (Figure~\ref{fig:advdiff-prec-bars}).

\begin{figure}
	\centering
	\includegraphics[width=0.7\textwidth]{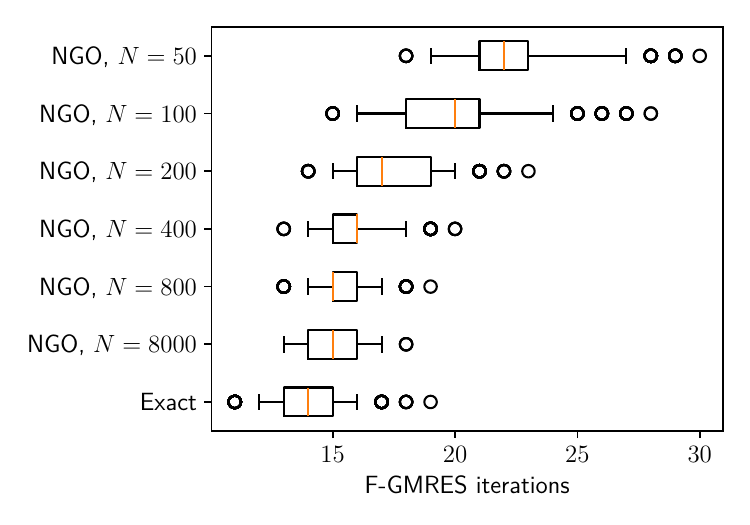}
	\caption{The distributions of F-GMRES iterations required to solve \(80 \times 80\) advection-diffusion problems, depending on the number \(N\) of solutions used to train the NGO.}
	\label{fig:advdiff-datasetsize}
\end{figure}

Figure~\ref{fig:advdiff-datasetsize} shows how the performance of NGO-based preconditioners depends on the size \(N\) of the training data set. From this figure, it appears that training on as little as 400 example solutions is sufficient to obtain an effective coarse correction. Even the NGOs trained on just \(50\) and \(100\) solutions produce coarse corrections that speed up GMRES compared to the smoother-only baseline. This is surprising, as the NGO uses a coarse space of \(144\) (or \(120\) with the Dirichlet basis functions excluded), meaning that the first two models of Figure~\ref{fig:advdiff-datasetsize} were trained on data sets that obviously cannot span the full coarse space. As was also found from the results in Figure~\ref{fig:diff-prec-bars}, the training data set does not need to span the full coarse space for NGOs to learn effective preconditioners. Nonetheless, increasing the dataset size to \(200\) or \(400\) greatly improves the learned preconditioner.

Naturally, training models with larger coarse spaces on more complex PDEs will likely require more than \(400\) samples. As such, we obviously cannot make the universal recommendation that \(400\) samples is sufficient. Nevertheless, these results show that the training data set size may only need to be a few times larger than the coarse-space dimension in order to obtain near-optimal preconditioners.

\begin{remark}
	For the NGOs trained in this section, we maintained an approximate 8:1 split of training and validation data. For example, the NGO trained on 400 samples additionally used 50 samples as validation data to prevent overfitting.
\end{remark}

\FloatBarrier
\subsection{Helmholtz equation}
\label{subsec:helmholtz}

The Helmholtz equation is qualitatively harder than the previous two cases. It is indefinite: \(-\Delta - k^2\) has both positive and negative eigenvalues, with the negative ones corresponding to oscillatory modes of which the number scales with \(k^d\) in dimension \(d\). The discrete operator \(\mat{A}\) inherits this indefiniteness, and the resulting linear systems are notoriously hard to solve unless \(k\) is small.

\paragraph*{Classical challenges}
The classical literature on Helmholtz coarse spaces is unanimous on one point: an effective coarse space for indefinite Helmholtz \emph{cannot be wave-number-independent}. Standard analyses of two-level domain-decomposition methods for Helmholtz \cite{bonazzoli2018twolevel, bootland2023overlapping, conen2014coarse, dolean2026deltakgeneo, graham2017domain} all show, by analysis or by numerical experiment, that the dimension of an effective coarse space must grow with $k$. The empirical comparisons of Bootland \etal\cite{bootland2021comparison} make the consequences concrete: only spectral coarse spaces such as DtN \cite{conen2014coarse} or GenEO, where the coarse-space dimension is \emph{automatically} chosen to track \(k\), yield wave-number-robust convergence. SPD-based coarse spaces can be made to work for indefinite Helmholtz \cite{dolean2026deltakgeneo} but only for low to moderate frequencies, with a delicate threshold choice and, again, with a coarse-space dimension that grows with \(k\).

This sets a hard constraint on what we can hope to achieve. Our coarse basis is a fixed \(12 \times 12\) tensor-product B-spline basis, meaning its dimension does not grow with \(k\). Classical theory therefore predicts that no neural operator built on top of this basis can yield a wave-number-robust preconditioner, regardless of architecture. The Helmholtz section of our experiments tests \emph{not} whether learned coarse corrections can solve high-\(k\) Helmholtz (they cannot, and we should not expect them to) but whether machine learning-based coarse corrections can achieve similar convergence rates as their exact counterparts, and whether the architectural ranking established in the previous two sections persists, in the wave-number range where the basis is still rich enough for the exact coarse solve to work.

\begin{figure}
	\centering
	\begin{subfigure}[t]{0.45\textwidth}
		\includegraphics[width=\textwidth]{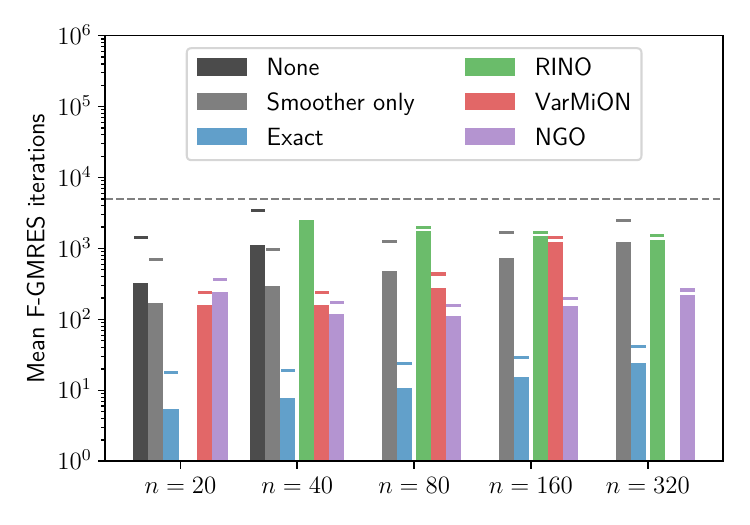}
		\caption{F-GMRES iterations by mesh size}
	\end{subfigure}
	\begin{subfigure}[t]{0.45\textwidth}
		\includegraphics[width=\textwidth]{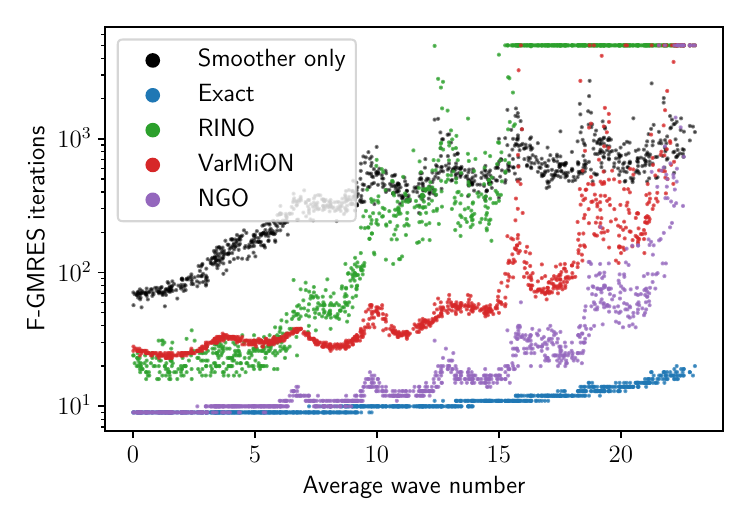}
		\caption{F-GMRES iterations by wave number}
	\end{subfigure}
	\caption{Helmholtz preconditioning. Left: average F-GMRES iterations for various discretisations and preconditioners. Right: F-GMRES iterations on the \(80 \times 80\) mesh as a function of the average wave number. The DeepONet is omitted because it failed to converge in almost every test.}
	\label{fig:helmholtz-prec-bars}
\end{figure}

\paragraph*{Results}
The results in Figure~\ref{fig:helmholtz-prec-bars} show that the fixed-size coarse space, although a theoretical hard limit, is not the limiting factor in the effectiveness of the learned preconditioners. The DeepONet-based preconditioner fails to converge in almost every test case and is omitted from the figures (for a complete table of results, please see Table~\ref{tab:all-results}). Of the three remaining architectures, the NGO is consistently the most effective. As the wave number grows, however, all three learned preconditioners deteriorate: the RINO becomes counterproductive for \(\langle k \rangle > 15\) and ceases to converge at all; the VarMiON converges across the wave-number range but is uniformly worse than the RINO at low \(k\); the NGO outperforms both but still deteriorates compared to the exact coarse solve as \(k\) grows. The exact two-level preconditioner with the same fixed-size basis remains effective on this entire range, which is why the range is interesting: within it, the failure of the learned methods is not due to an insufficient coarse space, but due to the poor quality of the learned coarse solution operator. As such, the fixed size of the coarse space is not the limiting factor in using machine learning-based coarse corrections for this equation: the learned coarse corrections already perform significantly worse than exact coarse solves even in the wave number regime where the coarse space is sufficiently rich.

Table~\ref{tab:helmholtz-failures} shows the failure counts and mean number of iterations excluding failures for the \(80 \times 80\) discretisations. As was the case for earlier test problems, preconditioners that use fewer iterations on average also fail less often. The only exception to this rule is the RINO, which uses fewer iterations on average than the smoother without coarse correction, even though the latter always converges. The reason for this is visible in Figure~\ref{fig:helmholtz-prec-bars}: the RINO-based preconditioner consistently fails for high wave numbers, which means these difficult cases are excluded from the RINO's average in Table~\ref{tab:helmholtz-failures}.

\begin{table}
	\centering
	\caption{Average number of F-GMRES iterations and number of failures for different preconditioners on \(80 \times 80\) discretisations of the Helmholtz equation. The average iteration counts now \textit{exclude} the cases where convergence was not achieved.}
	\label{tab:helmholtz-failures}
	\begin{tabular}{l r r}
		\toprule
		Preconditioner & Mean iterations & Failures \\
		\midrule
		Smoother only  & 481.8           & 0        \\
		Exact          & 10.9            & 0        \\
		RINO           & 216.4           & 323      \\
		VarMiON        & 127.9           & 31       \\
		NGO            & 37.5            & 15       \\
		\bottomrule
	\end{tabular}
\end{table}

\begin{figure}
	\centering
	\begin{subfigure}{0.9\textwidth}
		\includegraphics[width=0.9\textwidth]{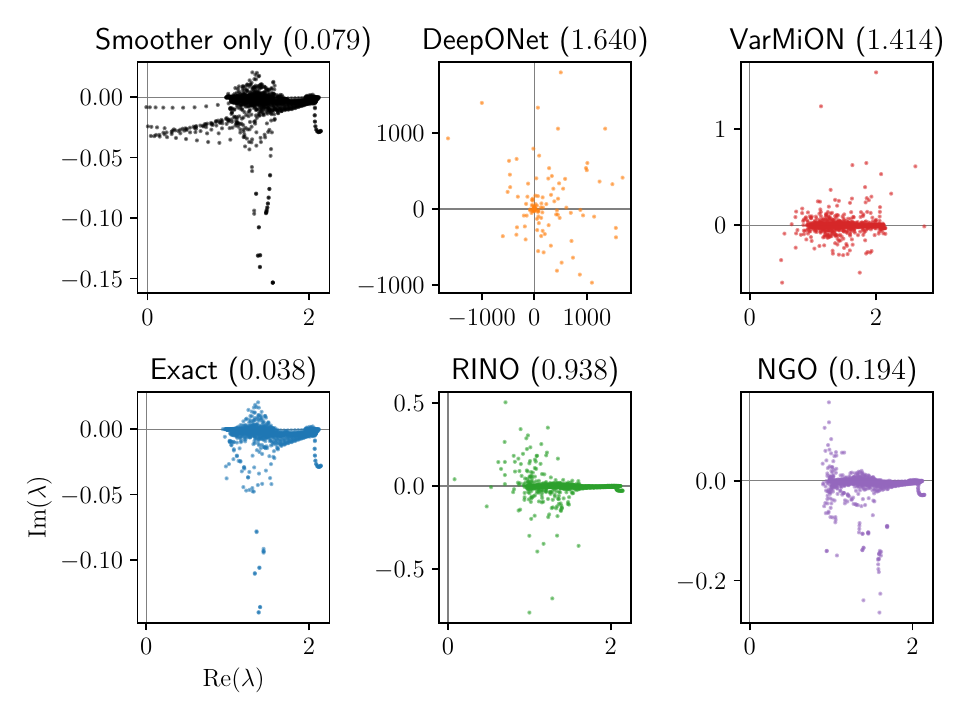}
		\caption{\(k = 5\)}
	\end{subfigure}
	\begin{subfigure}{0.9\textwidth}
		\includegraphics[width=0.9\textwidth]{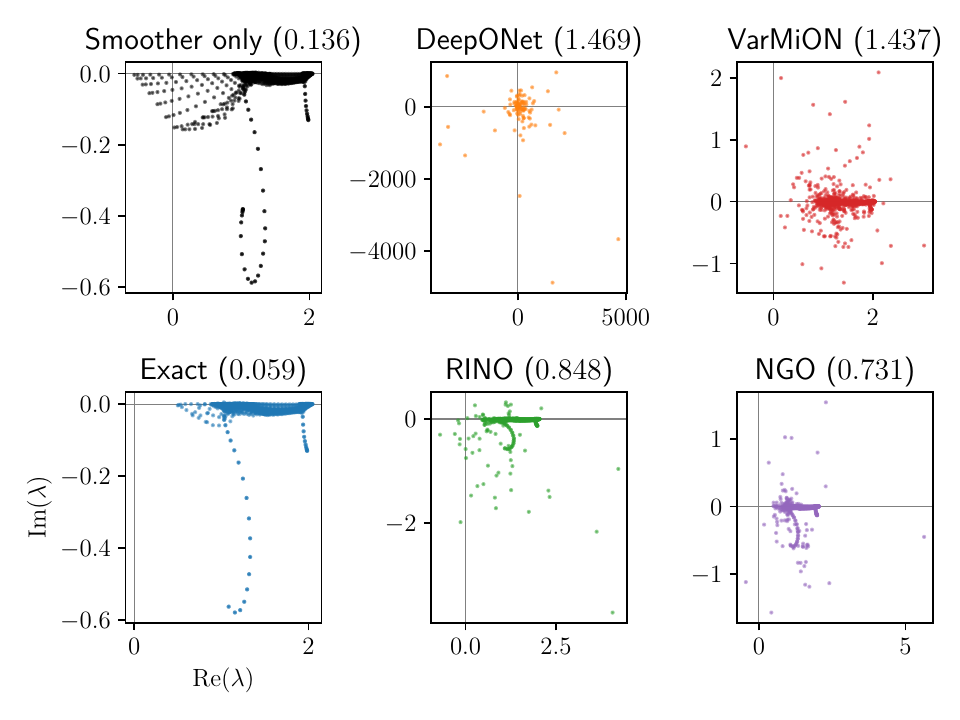}
		\caption{\(k = 20\)}
	\end{subfigure}
	\caption{Helmholtz preconditioned spectra on a \(40 \times 40\) mesh. The numbers in parentheses are the commutator errors as defined in \eqref{eq:commutator-error}. At \(k = 5\) the problem is close to positive definite and all three viable architectures cluster near \([1, 2]\); at \(k = 20\) the exact preconditioner remains well-conditioned but all three learned ones acquire near-imaginary eigenvalues.}
	\label{fig:helmholtz-spectra}
\end{figure}

The spectral analysis (Figure~\ref{fig:helmholtz-spectra}) provides more insight into the convergence results. At low wave number (\(k = 5\)), the problem is still close to positive definite and all three viable architectures produce spectra clustered near the real interval \([1, 2]\), consistent with the fast convergence of F-GMRES for these problems. At high wave number (\(k = 20\)), the exact A-DEF1-preconditioned spectrum still avoids the origin, meaning the basis is still rich enough for the exact coarse solve. However, all three learned-preconditioned spectra now have eigenvalues near or across the imaginary axis, leading to indefiniteness in the preconditioned system. The DeepONet produces a wildly diffuse spectrum extending into all four quadrants of the complex plane, which explains its outright failure at all \(k\).

\begin{figure}
	\centering
	\begin{subfigure}[t]{0.45\textwidth}
		\includegraphics[width=\textwidth]{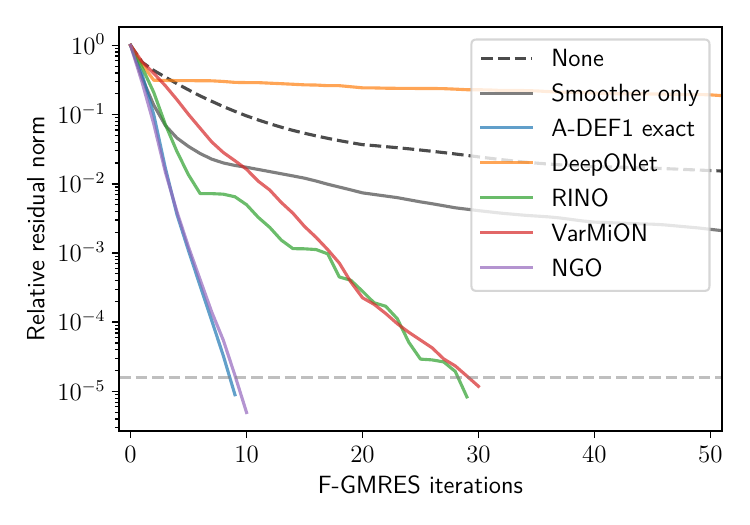}
		\caption{Low wave number, \(k = 5\)}
	\end{subfigure}
	\begin{subfigure}[t]{0.45\textwidth}
		\includegraphics[width=\textwidth]{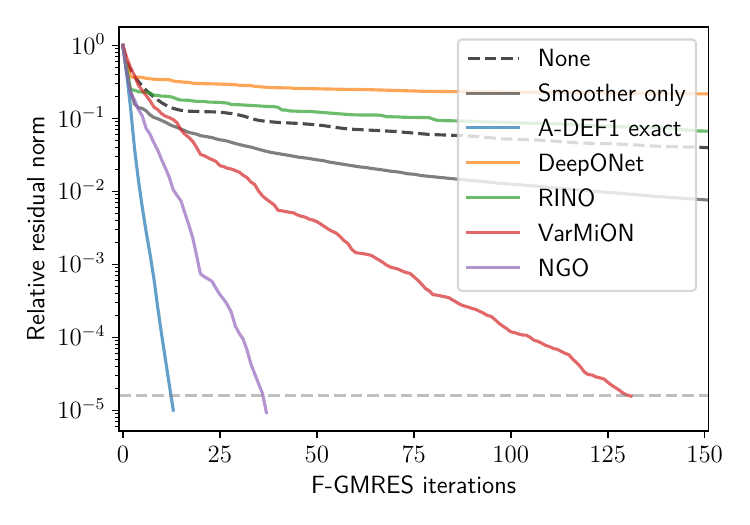}
		\caption{High wave number, \(k = 20\)}
	\end{subfigure}
	\caption{Convergence histories of F-GMRES for representative Helmholtz problems on the \(80 \times 80\) mesh. The exact coarse solve maintains monotone convergence at both wave numbers; the learned methods deteriorate at high \(k\).}
	\label{fig:helmholtz-conv-history}
\end{figure}

Figure~\ref{fig:helmholtz-conv-history} shows the corresponding F-GMRES convergence histories. At \(k = 5\), the NGO and exact coarse-solve curves are nearly indistinguishable, with the VarMiON and RINO converging steadily, albeit slower than the NGO, while the DeepONet quickly stagnates. At \(k = 20\), the learned methods slow down dramatically, while the exact coarse solve continues to drive the residual down quickly. The NGO is now roughly a factor 2 slower than exact A-DEF1, while the VarMiON converges much more slowly and the RINO fails entirely.

\paragraph*{Failure modes and prospects}
The failure of learned coarse corrections at high \(k\) therefore has two components: first, within the regime where the basis is rich enough (i.e.\ where the exact coarse solve still converges) learned coarse corrections significantly deteriorate compared to the exact coarse solve. However, the architectural ranking of Table~\ref{tab:no-design-axes} continues to hold, with the NGO consistently the best of the four learned options. Second, there is the theoretical limit that the fixed basis is incapable of resolving the propagating modes that the Helmholtz operator supports at high \(k\), and no choice of network can recover what is missing from the basis itself.

Thus, producing learned coarse corrections that are scalable to high \(k\) requires solving two problems: first, the deterioration compared to exact coarse solves must be resolved, by improving network architecture, training procedure, or other factors. Second, scaling to higher wave numbers requires a model of which the coarse space can grow with increasing wave number. We do not pursue any of these here; the contribution of the present section is to identify the failure modes precisely, and connect these to the classical Helmholtz coarse-space literature so that the obvious next steps become visible.

\begin{remark}
	The Helmholtz wave-number tests do not involve out-of-distribution evaluation: the models are trained and tested on the same range of \(k\). The deterioration with \(k\) is therefore not a generalisation problem in the machine-learning sense.
\end{remark}

\FloatBarrier
\subsection{Computational cost compared to AMG}
\label{subsec:vs-amg}
While the comparisons so far have focused on the number of Krylov iterations, the ultimate goal is of course to reduce the wall-clock time taken to obtain solutions. In this section, we present some analysis and experimental results indicating that two-level preconditioners with neural operators as coarse solves may indeed be competitive with other methods when measured in terms of wall-clock time.

First, we briefly analyse the computational complexity of various coarse corrections. For the VarMiON and NGO, one application of the preconditioner can be written explicitly as in~\eqref{eq:varmion-preconditioner} and~\eqref{eq:ngo-preconditioner}. For both models, the matrix-vector product with the \(m \times m\) matrix \(\mat{C}(\boldsymbol{\theta})\) has a computational complexity of \(\O{m^2}\). The computational cost of constructing the matrix \(\mat{C}(\boldsymbol\theta\)), although more difficult to analyse, is at least of order \(m^2\) as the output is of size \(m^2\). The computational complexity of the DeepONet and RINO is less obvious, but can still be considered \(\O{m^2}\) under the assumptions that (1) the depth of the neural networks is bounded, and (2) the size of the hidden layers in the neural network is proportional to \(m\). This quadratic complexity is significantly worse than the \(\O{m}\) complexity of a multigrid cycle, meaning that neural coarse corrections of this form are not expected to be competitive with multigrid for very large systems. However, the performance of the neural operators is highly architecture-dependent, and due to the similarities between U-Nets and multigrid methods (as observed in earlier works \cite{he2024mgno, hsieh2019learning, özbay2021poisson}) it is possible to construct neural operator architectures for which a single application only has a linear cost. For example, convolutional neural networks (CNNs) have a time complexity that is linear in the input size (for fixed convolutional kernel size and channel count), as pointed out by \"Ozbay \etal\cite{özbay2021poisson}. U-Nets, which only combine convolutions and pointwise operations, therefore also have linear complexity. Neural operator-based coarse corrections that are computationally efficient will require combining the architectural features discussed in this work with such a convolution-based architecture to ensure that the resulting operator is both effective as coarse correction and efficient to evaluate.

\begin{figure}
	\centering
	\includegraphics[width=0.45\textwidth]{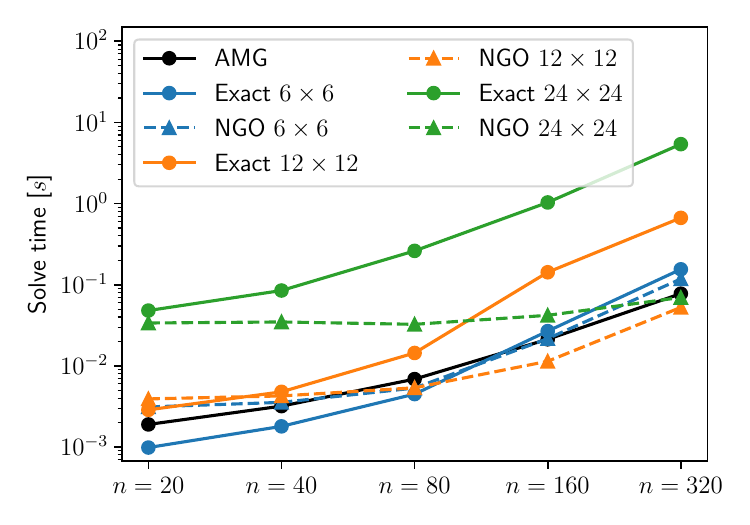}
	\includegraphics[width=0.45\textwidth]{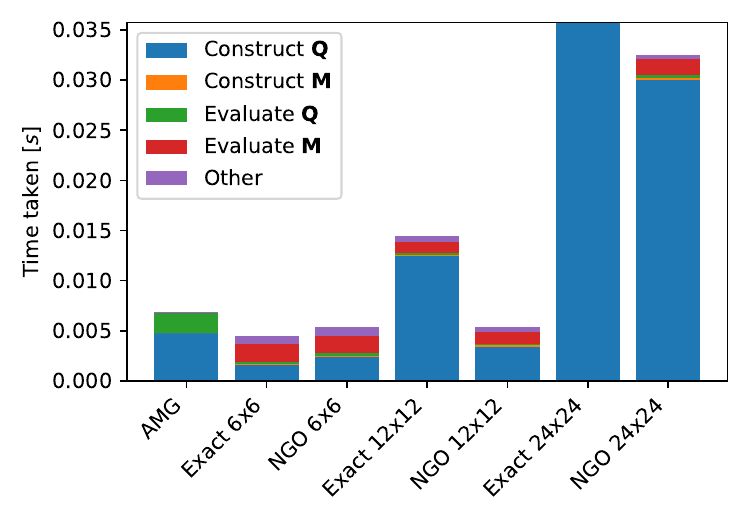}
	\caption{Left: average time taken to solve 1000 diffusion problems using conjugate gradients with various preconditioners. Right: breakdown of the computation time required to solve \(80 \times 80\) discretisations of diffusion problems. The exact \(24 \times 24\) two-level preconditioner is not shown fully.}
	\label{fig:vs-amg}
\end{figure}

Figure~\ref{fig:vs-amg} shows how the exact and learned two-level preconditioners compare to algebraic multigrid (AMG) when solving diffusion equations with preconditioned conjugate gradient. One figure shows the average time taken to solve 1000 linear systems for various mesh sizes and preconditioners, while the other breaks down the time taken for the \(80 \times 80\) discretisations into the construction of the preconditioners and their application inside the CG algorithm. For clarity, Table~\ref{tab:vs-amg} shows the same data as well as the average number of CG iterations required. From this data, we can see that NGO-based coarse corrections provide a significant speed-up over exact coarse corrections: this is largely due to the fact that the coarse matrix \(\mat{C}(\boldsymbol{\theta})\) is cheaper to compute than \(\mat{E}^{-1}\), while being equally effective as a preconditioner. The difference grows as the mesh size increases, due to the fact that the computation of \(\mat{E} = \mat{Z}\transpose\,\mat{A}\,\mat{Z}\) becomes slower as the mesh size increases, while the time taken to apply an NGO remains approximately constant. However, it is clear that the presence of dense \(m \times m\) matrices in the current NGO architecture hinders the scalability of this approach: for example, the NGO with \(24 \times 24\) basis must construct a dense \(576 \times 576\) matrix, which is computationally expensive and results in uncompetitive performance except on the finest grid. This is an architectural limitation that we do not aim to address in this work. Nonetheless, we see that the \(12 \times 12\) NGO, which produces a smaller \(144 \times 144\) matrix, actually results in the fastest overall solves for the \(160 \times 160\) and \(320 \times 320\) discretisations. This indicates that preconditioners with machine learning-based coarse corrections can indeed result in performance that is competitive with classical methods, but that the construction of the matrix \(\mat{C}(\boldsymbol\theta)\) and subsequent matrix-vector products must be replaced by a more computationally efficient representation of a linear function in order to remain performant on problems where larger coarse spaces are required.

We must note that the coarse space and smoother were not tuned for optimal performance, so better performance may be obtained with more careful tuning of parameters such as the depth of the neural networks and the sizes of the hidden layers. Furthermore, to ensure a fair comparison of the computational work, all the reported timings were obtained by running all code on one \emph{CPU} thread. In practice, the parallelisable nature of neural networks means that the performance comparison may be more in favour of ML-based approaches when utilising multi-threading or GPU-acceleration. Nevertheless, these results show that neural coarse corrections have the potential to produce preconditioners that are competitive with classical methods. In particular, compared to classical two-level preconditioners, the use of neural operators can reduce the setup time by eliminating the expensive construction of the coarse system by a computationally efficient neural network. Although the resulting preconditioned systems may need more iterations to solve, the lower computational cost per iteration almost compensates for the slower convergence.

\begin{table}
	\centering
	\caption{Comparison of iteration counts and wall-clock time for NGO-based two-level preconditioning and multigrid with PyAMG. Reported times are in milliseconds. Note that the ``Total'' column does not equal the sum of the other columns due to additional overhead of the Conjugate Gradient algorithm.}
	\label{tab:vs-amg}
	\begin{tabular}{l r | r r r r | r}
		\toprule
		\multirow{2}*{Preconditioner} & \multirow{2}*{CG iterations} & \multicolumn{2}{c}{Setup time} & \multicolumn{2}{c|}{Solve time} & \multirow{2}*{Total}                       \\
		                              &                              & \(\mat{Q}\)                    & \(\mat{M}\)                     & \(\mat{Q}\)          & \(\mat{M}\) &       \\
		\midrule
		AMG                           & 4.0                          & 4.7                            & 0.0                             & 2.0                  & 0.0         & 6.9   \\
		Exact \(6\times6\)            & 16.4                         & 1.5                            & 0.1                             & 0.2                  & 1.8         & 4.5   \\
		NGO \(6\times6\)              & 16.4                         & 2.4                            & 0.1                             & 0.2                  & 1.8         & 5.3   \\
		Exact \(12\times12\)          & 10.1                         & 12.5                           & 0.1                             & 0.2                  & 1.1         & 14.4  \\
		NGO \(12\times12\)            & 10.2                         & 3.4                            & 0.1                             & 0.2                  & 1.1         & 5.4   \\
		Exact \(24\times24\)          & 7.3                          & 258.9                          & 0.1                             & 0.3                  & 1.4         & 261.2 \\
		NGO \(24\times24\)            & 8.0                          & 30.0                           & 0.1                             & 0.4                  & 1.5         & 32.5  \\
		\bottomrule
	\end{tabular}
\end{table}

\FloatBarrier
\subsection{A 3D problem}
\label{subsec:3d-diffusion}
The performance study performed in the previous section also has implications for 3D PDEs: in discretisations of problems with a three-dimensional domain, the number of unknowns is typically much higher. As such, the unfavourable computational complexity of the current models limits us to very small 3D problems. Nevertheless, to show that the proposed approach is applicable to three-dimensional problems, we train an NGO on a 3D variant of the diffusion problem \eqref{eq:diffusion-pde} over the domain \(\Omega = (0, 1)^3\) (see Appendix~\ref{app:3d-diffusion} for details). Our NGO uses an \(8 \times 8 \times 8\) cubic B-spline basis, resulting in a coarse space of dimension \(m = 512\), similar in size to the \(576\)-dimensional coarse space of the \(24 \times 24\) NGO of the previous section. However, after removing the Dirichlet-supported basis functions, the actual size of the coarse space is only \(6 \times 8 \times 8 = 384\).
\begin{table}
	\centering
	\caption{Mean F-GMRES iteration counts for discretisations of 3D diffusion with various preconditioners. Note that \(n\) is the number of discretisation points in one dimension, meaning that the total number of unknowns in the discretisation is approximately \(n^3\).}
	\label{tab:prec-diffusion3d}
	\begin{tabular}{l | rrr}
		\toprule
		Preconditioner      & \(n=20\) & \(n=40\) & \(n=80\) \\
		\midrule
		None                & 65.4     & 244.6    & 1026.0   \\
		Smoother only (IC)  & 18.1     & 45.0     & 124.1    \\
		A-DEF1 (IC + exact) & 5.5      & 11.0     & 23.0     \\
		A-DEF1 (IC + NGO)   & 8.0      & 16.6     & 37.3     \\
		\bottomrule
	\end{tabular}
\end{table}
Table~\ref{tab:prec-diffusion3d} shows the results. Note that the incomplete Cholesky preconditioner, although not robust, is used here as a representative smoother, as the choice of smoother was found to have little to no effect on the comparison between NGO-based and exact coarse corrections. The 3D results differ from those of the 2D case above; substituting the exact coarse solve by an NGO now comes at a significant increase in the required number of F-GMRES iterations.

The reason for the deterioration is not clear; one explanation may be that the NGO should have a larger neural network with more parameters in order to approximate the exact coarse solves more accurately. An insufficient training data set may also be a cause (the data set consists of 10000 problems as in the 2D case).

Nonetheless, these results indicate that the same approach that is successful for 2D diffusion may also be used on 3D PDEs. As mentioned earlier, a significant challenge to overcome first is the development of more computationally efficient neural operator models, to allow scaling to larger coarse spaces.

\FloatBarrier
\subsection{Limitations and strengthening directions}
\label{subsec:limitations}

The numerical experiments above are deliberately scoped to isolate the architectural axes of Table~\ref{tab:no-design-axes}: the coarse basis, smoother, training procedure, and evaluation metrics are all kept the same. This isolation comes at the cost of two limitations that are worth stating clearly, together with the directions in which an extended study could go.

\paragraph*{2D-focused}
Almost all test problems are posed on the unit square, with only a brief experiment on a 3D diffusion problem which indicates that the NGO-based coarse corrections can indeed be used in 3D. A full 3D study would test whether the architectural ranking established in 2D persists in higher dimension, whether the integration-based architectures retain their resolution-invariance benefit on tensorised 3D quadrature rules, and whether the high-frequency basis-coverage condition identified in Section~\ref{subsec:diffusion} (Dataset~3) carries the same weight when the coarse-space dimension grows from \(m = 144\) to \(m = \mathcal{O}(10^3)\). We have not pursued this here because the goal of the present paper is to establish the architectural principles cleanly; an extended 3D study is the natural sequel.

\paragraph*{Exact coarse solves as reference}
We benchmark all four learned methods against the exact two-level preconditioner (with the same basis and a direct solve of \(\mat{E}\,\vec{w} = \vec{v}\)) and against the smoother alone. This is the right benchmark for the question we ask, but it is not the same as benchmarking against the strongest classical preconditioner available for each PDE. Two concrete extensions would tighten this story.
\begin{itemize}
	\item For all three PDEs: combining the learned coarse corrections with more advanced smoothers (such as those based on domain-decomposition) to obtain mesh-independent iteration counts.
	\item For Helmholtz: comparison against shifted-Laplacian multigrid  or a classical wave-number-dependent coarse space such as DtN or \(H_k\)-GenEO. As argued in Section~\ref{subsec:helmholtz}, no fixed-basis learned method can be expected to match these at high \(k\); the comparison would quantify the gap.
\end{itemize}
None of these extensions is technically difficult; we omit them because they would not affect the architectural conclusions of the paper.

\FloatBarrier\section{Conclusion}
\label{sec:conclusions}

We have studied a clean architectural question: when a neural operator is used as a coarse-space correction inside a two-level preconditioner with a fixed coarse basis, which of its design choices matter? The four-corner design grid of Table~\ref{tab:no-design-axes} isolated two axes: whether the input fields are discretised by sampling or by integration against a basis, and whether the network preserves the linear dependence of the solution on the source term. The numerical experiments confirmed the prediction made in Section~\ref{subsec:four-preconditioners}: integration is the dominant axis, controlling whether the resulting preconditioner can be made symmetric and therefore whether it is even compatible with PCG; linearity is the secondary but consequential axis, strongly influencing the preconditioned spectrum for non-self-adjoint problems. Among the four architectures tested, only the Neural Green's Operator (NGO) sits at the favourable corner of both axes, and only the NGO matches an exact coarse solve in iteration count on 2D diffusion and advection-diffusion.

Our additional evaluation of NGO-based coarse corrections shows that training such a model only needs little training data, and that the resulting model can even be applied to some problems outside of the training data distribution: generalisation to different source terms works without any modifications to the model, and generalisation over the PDE coefficients can be achieved by applying the model to a ``nearby'' in-distribution coefficient field. Generalisation to different domain geometries and boundary conditions requires architectural extensions and is therefore not considered here.

\paragraph*{Practical recommendations}
Two design choices are recommended for any neural operator-based coarse-space correction. First, the source-term inputs should be \emph{integrated} against the coarse basis rather than sampled at fixed sensor nodes; this aligns the row and column spaces of the learned operator and is what allows it to be symmetrised when the underlying problem is self-adjoint. Second, when the PDE is linear, the architecture should preserve the linear dependence of the solution on the source. Two methodological conditions are also necessary: the training data should populate the entire coarse space (not just a low-dimensional subset of it), and basis functions supported on Dirichlet boundaries must be excluded from the coarse correction. The latter is a small price: it can be enforced trivially by zeroing out the corresponding coefficients in the network output, requires no re-training, and is not associated with any loss of accuracy on the interior of the domain.

\paragraph*{Deployment of neural coarse corrections}
While coarse corrections based on neural networks, especially NGOs, are found to be very effective on average, there are still cases where they may result in slower convergence than expected, or even fail to converge entirely. While these issues mostly occur when applying the NGO to an out-of-distribution equation, such inputs may not be avoidable in practice. Therefore, a software library that uses neural operator-based preconditioners should also include a mechanism by which stagnation can be avoided. A simple solution could be to detect when the residual decreases too slowly (for example, when the residual drop \(\Vert \vec{r}^{(k)}\Vert / \Vert \vec{r}^{(k-l)} \Vert\) exceeds a given threshold), and switch back to exact two-level preconditioning in those cases. A well-designed fail-safe mechanism would ensure that the worst-case behaviour of the learned preconditioners is avoided, while still enabling a speed-up in most cases. However, in applications where \emph{consistently} fast linear solves are important, such as in real-time control problems, learned coarse corrections may still be unsuitable due to the inherent uncertainty in their performance.

\paragraph*{Helmholtz}
The Helmholtz experiments make precise where the presented approaches break down: at high wave numbers, all four architectures deteriorate and lose to the exact coarse solve. Although the fixed size of the coarse space poses a theoretical limit on the preconditioner performance for increasing wave number, a different problem arises first: neural operators are unable to match exact coarse corrections for moderate wave numbers. Resolving this issue is likely a matter of choosing the neural network architecture, training data, and training procedure, as there is no fundamental limitation that prevents neural operators from being effective in this regime. If this issue is resolved, generalising to even higher wave numbers requires re-thinking the architecture: a robust neural operator requires a basis of which the dimension increases with the wave number.

\paragraph*{Other directions}
Beyond the extensions discussed in Section~\ref{subsec:limitations}, four further directions stand out. First, application to PDEs whose coarse-space behaviour is qualitatively different, such as saddle-point systems arising from Stokes and Oseen problems. Second, alternatives to data-driven training: manufactured solutions \cite{hasani2025generating} or physics-informed losses \cite{li2024physics} would remove the need to solve many PDE realisations to generate training data. Third, training objectives aimed specifically at producing good preconditioners; the present results show that an \(L^2\) loss on the solution already produces useful preconditioners, but a loss that targets the spectral properties relevant to Krylov convergence could plausibly close the residual gap to the exact coarse solve. For example, positive results with such loss functions have been obtained by Dimola \etal\cite{dimola2025neural}. Finally, as discussed in Section~\ref{subsec:vs-amg}, the use of more computationally efficient architectures is necessary for learned coarse corrections to be competitive with classical preconditioners in terms of wall-clock time.

Code reproducing all experiments (including generation of the datasets, training of the models, and evaluation of the preconditioners) is available on GitHub at \url{https://github.com/HugoMelchers/neural-coarse-corrections}.

\paragraph*{Acknowledgments}
This project has received funding from the Eindhoven Artificial Intelligence Systems Institute (EAISI) under the Exploratory Multidisciplinary AI Research program for the NGO-PDE project. The authors thank Harald van Brummelen, Barry Koren, and Joost Prins for their valuable insights and helpful comments.

\printbibliography

\FloatBarrier\appendix
\section{PDE details}
\label{app:pde-details}
\subsection{Diffusion}
\label{app:pde-diffusion}
The exact equation is as follows:
\begin{subequations}
	\label{eq:diffusion-pde}
	\begin{align}
		-\nabla\cdot\left( \theta\nabla u \right) & = f \text{ on } \Omega = (0, 1)^2,                                      \\
		\theta \nabla u \cdot \mathbf{n}          & = \eta \text{ on } \Gamma_{\mathrm{N}} = \{\text{top}, \text{bottom}\}, \\
		\theta u                                  & = g \text{ on } \Gamma_{\mathrm{D}} = \{\text{left}, \text{right}\}.
	\end{align}
\end{subequations}

The training data is generated by solving the equations with the stabilised weak form:
\begin{align}
	\left\langle \theta\nabla v, \nabla u\right\rangle_\Omega  - \left\langle\theta \nabla v\cdot\mathbf{n}, u\right\rangle_{\Gamma_\mathrm{D}} - \left\langle v, \theta\nabla u\cdot\mathbf{n}\right\rangle_{\Gamma_\mathrm{D}} + \gamma\left\langle v, \theta u\right\rangle_{\Gamma_\mathrm{D}} \\
	= \left\langle v, f\right\rangle_\Omega - \left\langle \nabla v\cdot\mathbf{n}, g\right\rangle_{\Gamma_\mathrm{D}} + \left\langle v, \eta\right\rangle_{\Gamma_\mathrm{N}} + \gamma\left\langle v, g\right\rangle_{\Gamma_\mathrm{D}}.
\end{align}

\begin{figure}
	\centering
	\begin{subfigure}{\textwidth}
		\includegraphics[width=0.9\textwidth]{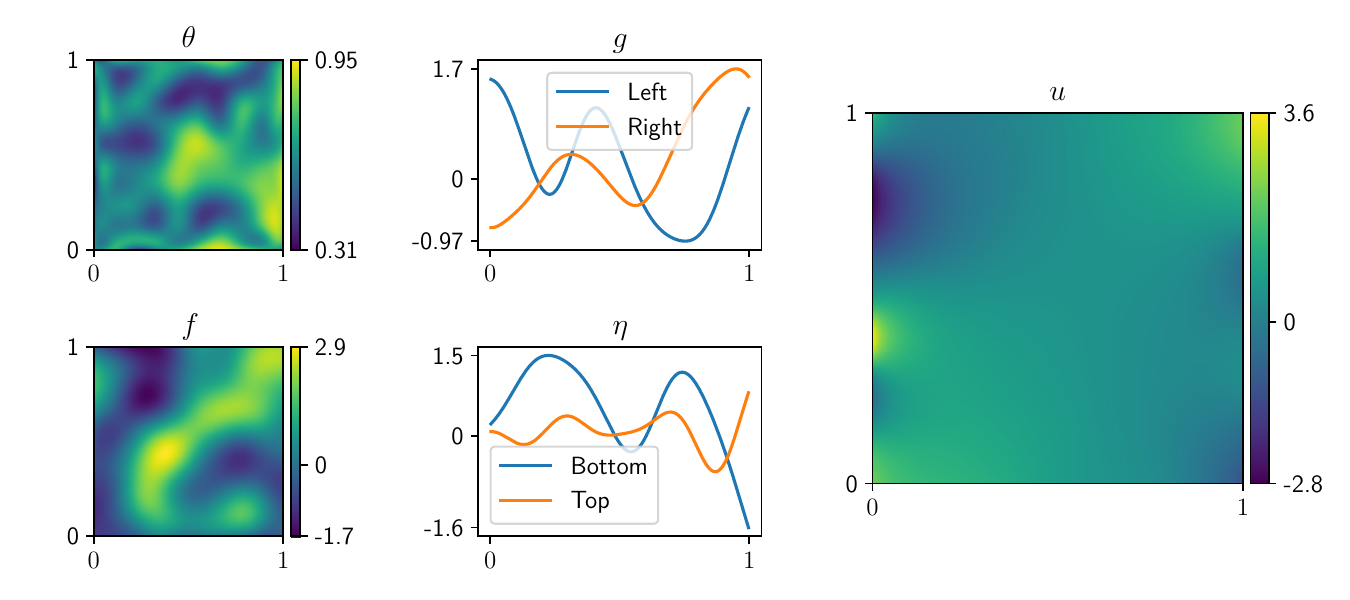}
		\caption{Dataset 1}
	\end{subfigure}
	\begin{subfigure}{\textwidth}
		\includegraphics[width=0.9\textwidth]{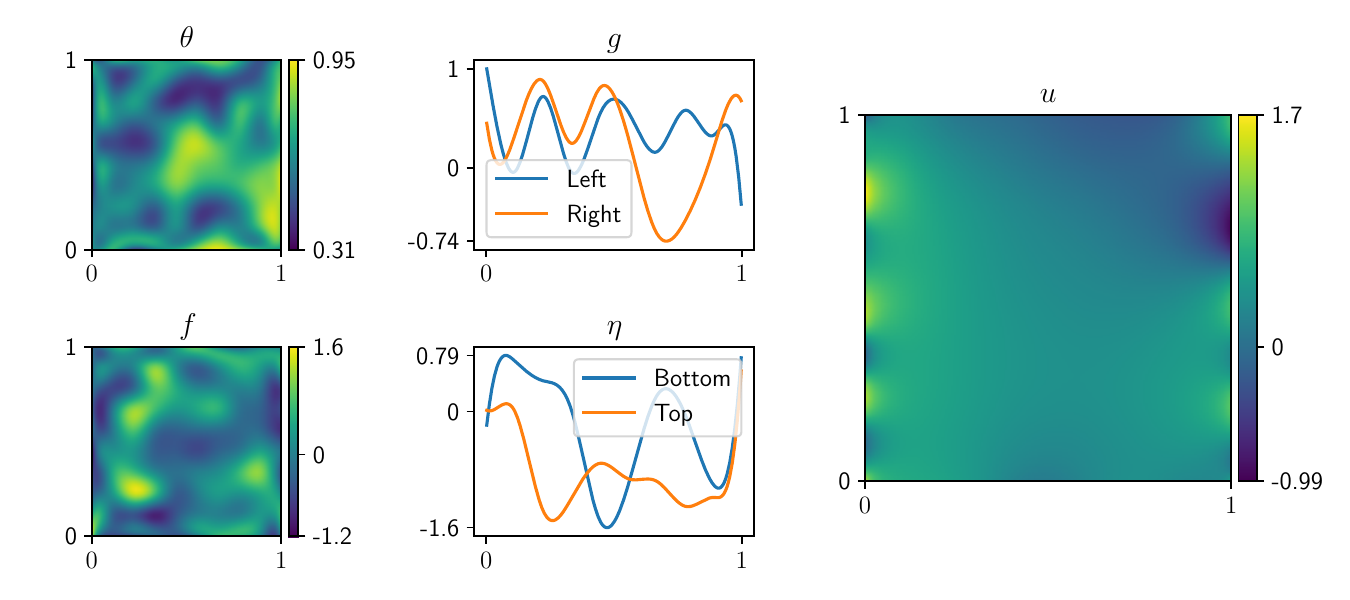}
		\caption{Dataset 2}
	\end{subfigure}
	\begin{subfigure}{\textwidth}
		\includegraphics[width=0.9\textwidth]{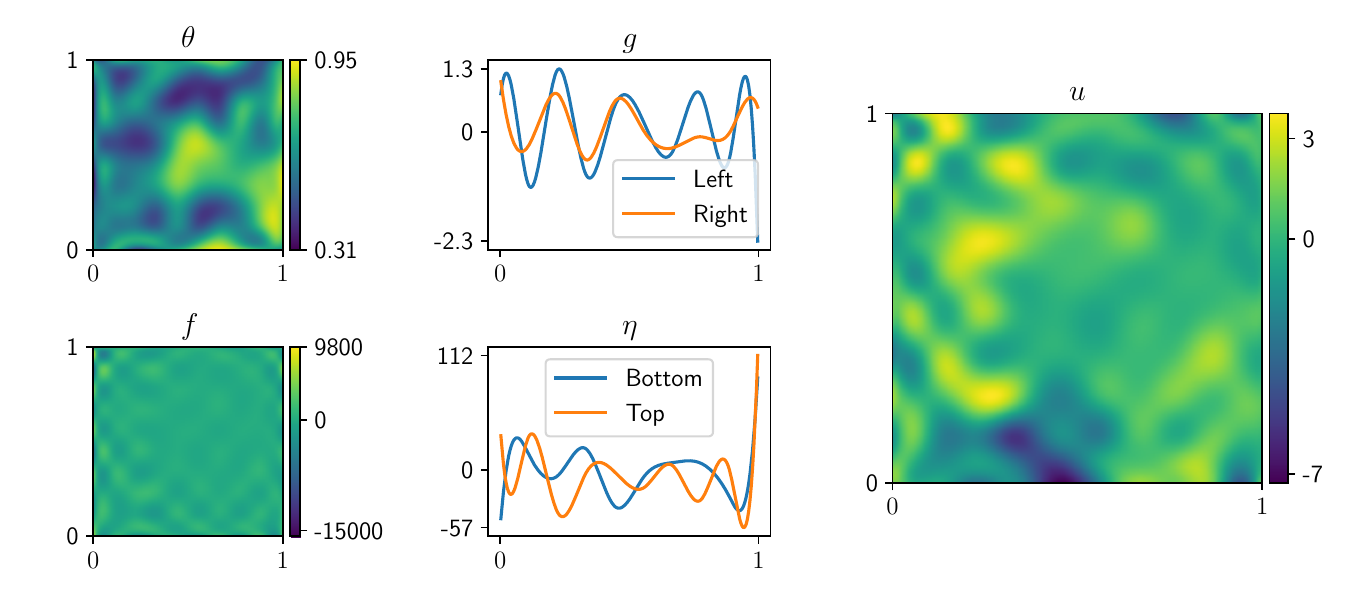}
		\caption{Dataset 3}
	\end{subfigure}
	\caption{Example of a diffusion problem from the three training data sets.}
	\label{fig:diffusion-data-sample}
\end{figure}

The above weak form is discretised on a \(50 \times 50\) cubic B-spline basis, and the resulting linear systems are solved using a direct sparse solver. Example problems from all three data sets are shown in Figure~\ref{fig:diffusion-data-sample}. The preconditioners are applied to finite-volume discretisations using the second-order accurate central difference scheme:
\begin{subequations}
	\label{eq:preconditioning-diffusion-discretisation}
	\begin{align*}
		-\nabla\cdot\left(\theta\nabla u\right)             & = -\partial_x\left(\theta\,\partial_xu\right) - \partial_y\left(\theta\,\partial_yu\right),                                                                                                       \\
		\left[\partial_x (\theta\,\partial_xu)\right]_{i,j} & \approx \frac{1}{\Delta x}\left(\frac{\theta_{i,j} + \theta_{i+1,j}}{2}\frac{u_{i+1,j} - u_{i,j}}{\Delta x} - \frac{\theta_{i-1,j} + \theta_{i,j}}{2}\frac{u_{i,j} - u_{i-1,j}}{\Delta x}\right), \\
		\left[\partial_y (\theta\,\partial_yu)\right]_{i,j} & \approx \frac{1}{\Delta y}\left(\frac{\theta_{i,j} + \theta_{i,j+1}}{2}\frac{u_{i,j+1} - u_{i,j}}{\Delta y} - \frac{\theta_{i,j-1} + \theta_{i,j}}{2}\frac{u_{i,j} - u_{i,j-1}}{\Delta y}\right),
	\end{align*}
\end{subequations}

\subsubsection{3D diffusion}
\label{app:3d-diffusion}
The 3D diffusion test case is similar, with a unit cube as the domain and Dirichlet boundary conditions on two of the six sides:
\begin{subequations}
	\label{eq:diffusion-pde-3d}
	\begin{align}
		-\nabla\cdot\left( \theta\nabla u \right) & = f \text{ on } \Omega = (0, 1)^3,                                      \\
		\theta \nabla u \cdot \mathbf{n}          & = \eta \text{ on } \Gamma_{\mathrm{N}} = \{\text{front}, \text{back}, \text{top}, \text{bottom}\}, \\
		\theta u                                  & = g \text{ on } \Gamma_{\mathrm{D}} = \{\text{left}, \text{right}\}.
	\end{align}
\end{subequations}
The PDE coefficients and source terms are generated using the three-dimensional analogue of the tensorised B-spline construction used for Dataset~3. The same stabilised weak form is used, but is now discretised using a \(15 \times 15 \times 15\) cubic B-spline basis. The finite-volume discretisations are analogous to those in 2D.

\subsection{Advection-diffusion}
\label{app:pde-advdiff}
\begin{figure}
	\centering
	\begin{subfigure}{\textwidth}
		\includegraphics[width=\textwidth]{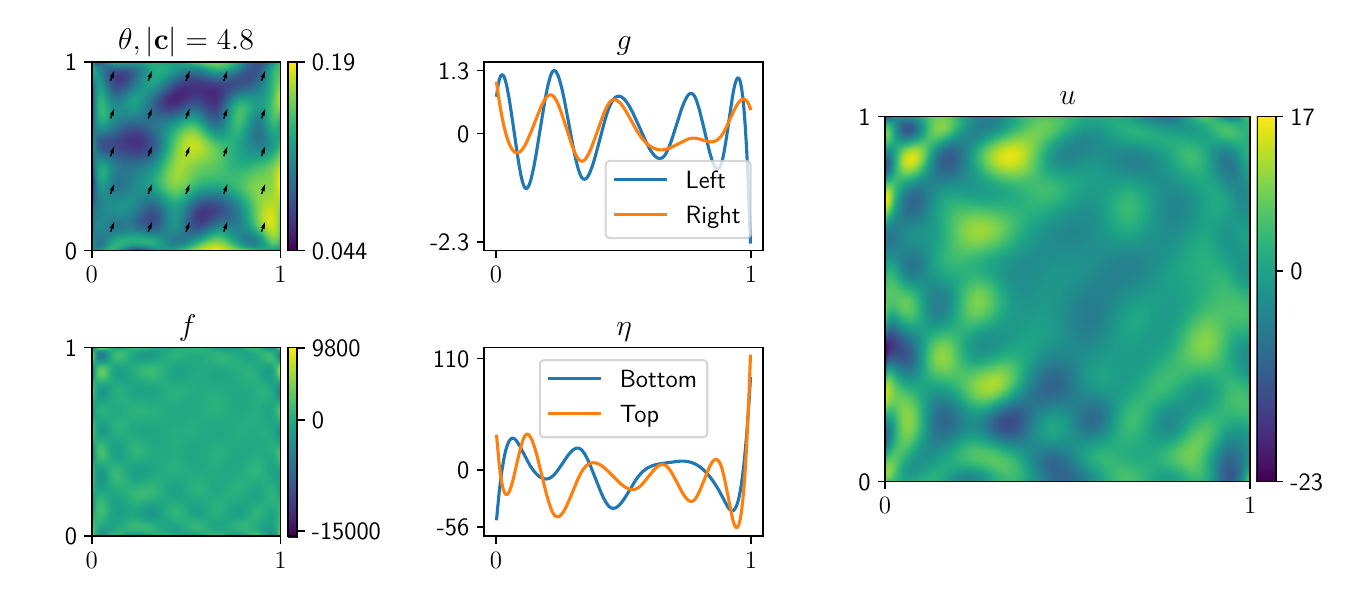}
		\caption{Low Péclet data set}
	\end{subfigure}
	\begin{subfigure}{\textwidth}
		\includegraphics[width=\textwidth]{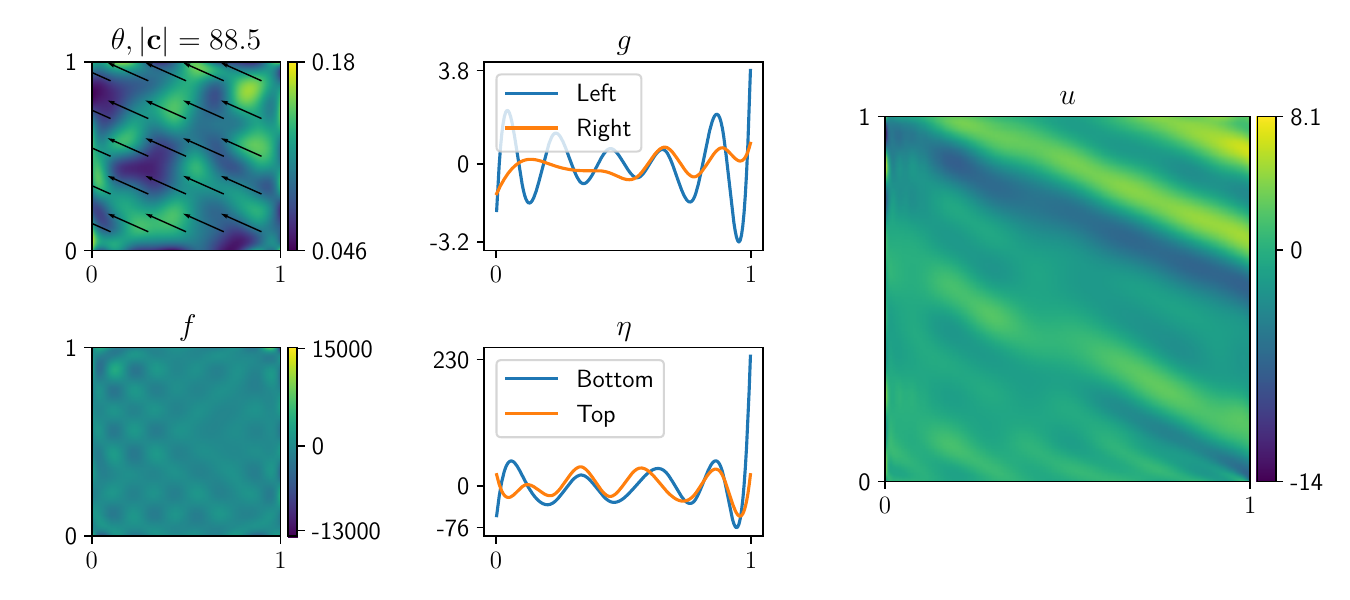}
		\caption{High Péclet data set}
	\end{subfigure}
	\caption{Example of an advection-diffusion problem from the training data set (top) and the out-of-distribution testing data set (bottom).}
	\label{fig:advdiff-data-sample}
\end{figure}

The exact equation is as follows:
\begin{subequations}
	\begin{align}
		-\nabla\cdot(\theta\nabla u) + \nabla\cdot(\mathbf{c}u)          & = f \text{ on } \Omega = (0, 1)^2,                                    \\
		\mathbf{n} \cdot \theta\nabla u - (\mathbf{c}\cdot\mathbf{n})_-u & = \eta \text{ on } \Gamma_\mathrm{N} = \{\text{bottom}, \text{top}\}, \\
		\theta u                                                         & = g \text{ on } \Gamma_\mathrm{D} = \{\text{left}, \text{right}\}.
	\end{align}
\end{subequations}
Note that \(\mathbf{c}\) is taken spatially constant in our experiments; this restriction is for simplicity only and is not inherent to either the discretisation or the smoother (see Appendix~\ref{app:smoothers}). The chosen setting is sufficiently complicated to investigate the effect of the Péclet number on the preconditioning. In the low-Péclet data set, the velocity field is generated with \(\vec{c}_x, \vec{c}_y\) sampled from a Gaussian distribution with mean \(\mu = 0\) and variance \(\sigma^2 = 25\). In the high-Péclet data set, we sample \(\vec{c} = (R\cos\beta, R\sin\beta)\) where \(R\) is uniformly sampled in \([0, 100]\) and \(\beta\) is uniformly sampled in \([0, 2\pi]\). Example problems from both data sets are shown in Figure~\ref{fig:advdiff-data-sample}.

Training data is generated using the stabilised weak form adapted from Bazilevs and Hughes~\cite{bazilevs2005weak}. We write \(\mathcal{L}\) for the advection-diffusion operator,
\begin{equation}
	\label{eq:advdiff-operator}
	\mathcal{L}\,u = -\nabla\cdot(\theta\,\nabla u) + \nabla\cdot(\mathbf{c}\,u),
\end{equation}
\(\tau\) for the SUPG (streamline-upwind) stabilisation parameter, and split the Dirichlet boundary into inflow and outflow parts according to the sign of \(\mathbf{c}\cdot\mathbf{n}\):
\begin{equation*}
	\Gamma_{\mathrm{D},\text{in}}  = \{\vec{x} \in \Gamma_\mathrm{D} : \mathbf{c}\cdot\mathbf{n} <    0\}, \quad
	\Gamma_{\mathrm{D},\text{out}} = \{\vec{x} \in \Gamma_\mathrm{D} : \mathbf{c}\cdot\mathbf{n} \geq 0\}.
\end{equation*}
The constants \(\gamma\) and \(C/h\) are Nitsche penalty parameters that weakly enforce the Dirichlet condition \cite{bazilevs2005weak}, with \(h\) the local mesh size. With these conventions the stabilised weak form reads
\begin{align*}
	 & (-\nabla v, \mathbf{c}u - \theta\nabla u)_\Omega
	+ (\mathcal{L}v\,\tau, \mathcal{L}u)_\Omega
	+ (v, -\theta\nabla u\cdot\mathbf{n} + \mathbf{c}\cdot\mathbf{n}u)_{\Gamma_\mathrm{D}}
	+ (-\gamma\theta\nabla v\cdot\mathbf{n} - \mathbf{c}\cdot\mathbf{n}v, u)_{\Gamma_{\mathrm{D},\text{in}}}                          \\
	 & + (-\gamma\theta\nabla v\cdot\mathbf{n}, u)_{\Gamma_{\mathrm{D},\text{out}}}
	+ \frac{C}{h}(v\theta, u)_{\Gamma_\mathrm{D}}
	+ \left(v, (\mathbf{c}\cdot\mathbf{n})_+u\right)_{\Gamma_\mathrm{N}} = (v, f)_\Omega
	+ (v, \eta)_{\Gamma_\mathrm{N}}
	+ (\mathcal{L}v\,\tau, f)_{\Omega}                                                                                                \\
	 & + \left(-\gamma\nabla v\cdot\mathbf{n} - \frac{1}{\theta}\mathbf{c}\cdot\mathbf{n}v, g\right)_{\Gamma_{\mathrm{D}, \text{in}}}
	+ (-\gamma\nabla v\cdot\mathbf{n}, g)_{\Gamma_{\mathrm{D},\text{out}}} + \frac{C}{h}(v, g)_{\Gamma_\mathrm{D}}.
\end{align*}

As is the case for diffusion, the discretisation is done on a \(50\times50\) cubic B-spline basis with direct sparse solves. The preconditioners are applied to finite-volume discretisations obtained using the complete flux scheme of ten Thije-Boonkkamp and Anthonissen~\cite{tenthijeboonkkamp2011finite}.

\subsection{Helmholtz}
\label{app:pde-helmholtz}

The equation is as follows:
\begin{subequations}
	\label{eq:helmholtz}
	\begin{align}
		-\Delta u - k^2 u                   & = f    \text{ on } \Omega = (0, 1)^2, \label{eq:helmholtz-pde}                                            \\
		\nabla u \cdot \mathbf{n} - i\,k\,u & = \eta \text{ on } \Gamma_{\mathrm{N}} = \{\text{top}, \text{bottom}\}, \label{eq:helmholtz-bc-neumann}   \\
		u                                   & = g    \text{ on } \Gamma_{\mathrm{D}} = \{\text{left}, \text{right}\}, \label{eq:helmholtz-bc-dirichlet}
	\end{align}
\end{subequations}
where \(k\) is the spatially varying wave number and \(i = \sqrt{-1}\). The Robin condition \eqref{eq:helmholtz-bc-neumann} is the first-order Sommerfeld absorbing boundary condition, modelling out-going waves through \(\Gamma_{\mathrm{N}}\). Note that, unlike for diffusion and advection-diffusion, the equation \eqref{eq:helmholtz-pde} is unweighted: there is no diffusion coefficient \(\theta\) and the Dirichlet datum \(g\) is imposed directly on \(u\).

\begin{figure}
	\centering
	\includegraphics[width=\textwidth]{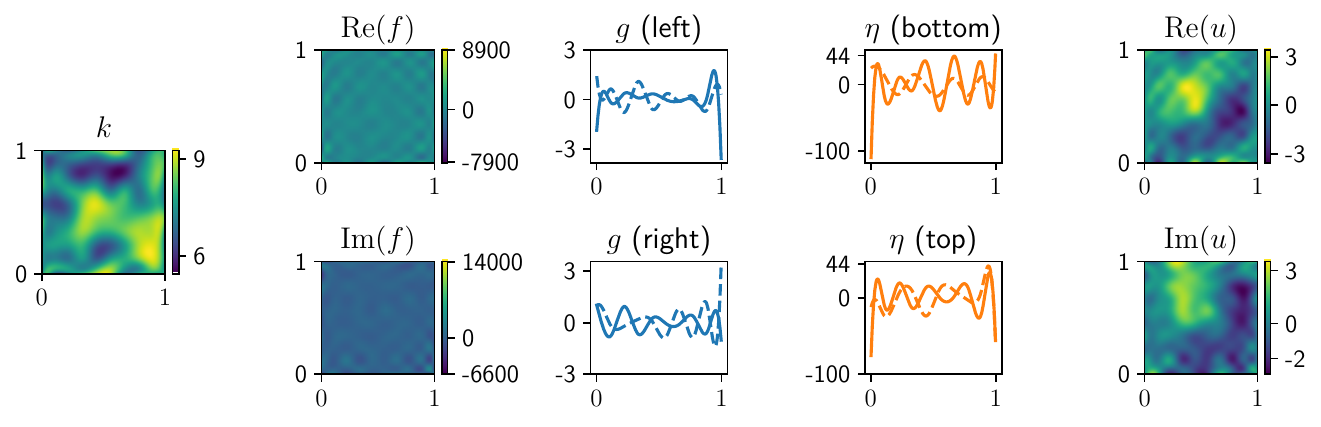}
	\caption{Example of a Helmholtz problem from the training data set. For \(g\) and \(\eta\), the solid and dashed lines indicate the real and imaginary parts, respectively.}
	\label{fig:helmholtz-data-sample}
\end{figure}

The data is generated with the finite element method using the following weak form:
\begin{align*}
	\int_{\Omega} \overline{\nabla v} \cdot \nabla u - k^2\,\overline{v}\,u\,d\vec{x} - i\int_{\Gamma_N} k\,\overline{v}\,u\,ds = \int_\Omega \overline{v}\,f\,d\vec{x} + \int_{\Gamma_N}\overline{v}\,\eta\,ds,
\end{align*}
where the Dirichlet boundary conditions are imposed strongly. The discretisation is done on a \(50 \times 50\) cubic B-spline basis and the resulting linear systems are solved using a direct sparse solver. A sample from the data set is shown in Figure~\ref{fig:helmholtz-data-sample}.

The preconditioner is applied to finite difference discretisations using the second-order accurate central difference scheme:
\begin{align*}
	-\Delta u - k^2u & \approx \frac{1}{h^2}\left( 4u_{i,j} - u_{i-1,j} - u_{i+1,j} - u_{i,j-1} - u_{i,j+1} \right) - k^2 u_{i,j}.
\end{align*}

\FloatBarrier
\section{Smoothers}
\label{app:smoothers}
In the three different PDEs considered, the NGO plays the role of the coarse-scale solver, while the smoother is chosen a priori and does not involve any machine learning.

\begin{figure}
	\centering
	\includegraphics[width=0.7\textwidth]{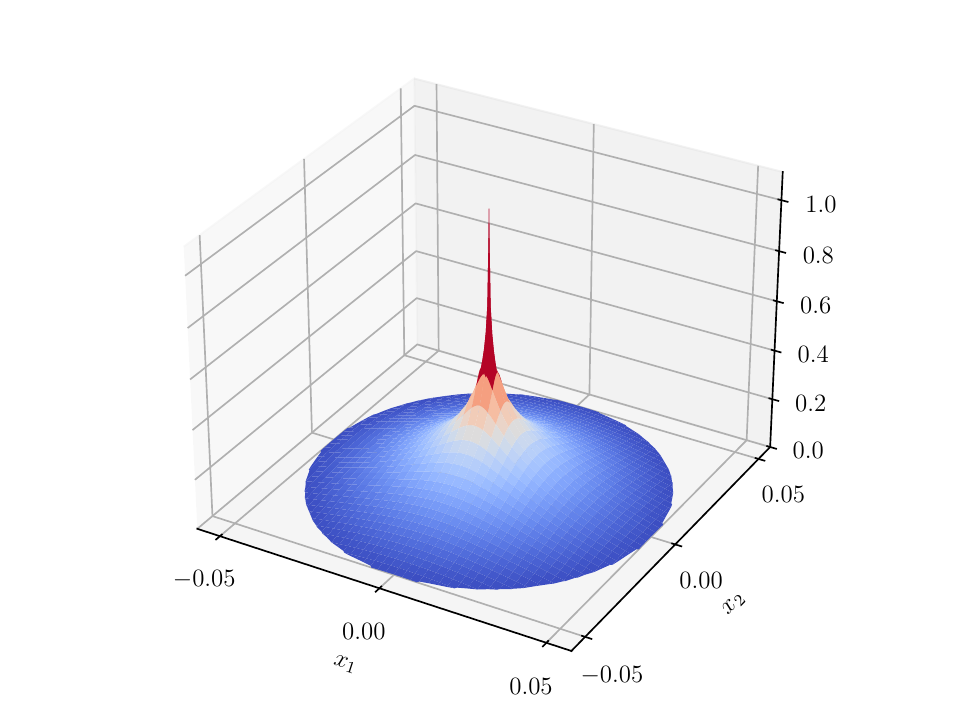}
	\caption{The convolutional kernel used as a smoother for the diffusion and Helmholtz equations.}
	\label{fig:diffusion-smoother}
\end{figure}

\paragraph*{Diffusion}
For the diffusion equation, the smoother is chosen as a convolution with a kernel based on the free-space Green's function for the 2D diffusion equation:
\begin{align*}
	G(\vec{x}, \vec{y}) & = \frac{1}{2\pi}\log\frac{1}{\abs{\vec{x} - \vec{y}}}.
\end{align*}
Since the coarse solver is expected to accurately capture the global effects of the source terms on the solution, the smoother must handle the localised effects. This is done by basing the smoother on a truncation of the Green's function to a radius \(H\):
\begin{align*}
	G_{\text{trunc}}(\vec{x}, \vec{y}) & = \frac{1}{2\pi}\max\left(\log\frac{H}{\abs{\vec{x} - \vec{y}}}, 0\right),
\end{align*}
which is locally supported, as \(G_{\text{trunc}}(\vec{x}, \vec{y}) = 0\) when \(\abs{\vec{x} - \vec{y}} > H\) (see Figure~\ref{fig:diffusion-smoother}). The smoother therefore computes the ``local solution component''. Since the Green's function above applies to the constant-coefficient case \(\theta \equiv 1\), the varying diffusion coefficient is used here by pre- and post-multiplication by \(\theta^{-1/2}\). This scaling has the effect of taking the magnitude of the diffusion coefficient into account, while still resulting in an SPD preconditioner:
\begin{align*}
	u_{\text{loc}}(\vec{x}) & = \frac{1}{2\pi\sqrt{\theta(\vec{x})}}\int_{B_{H}(\vec{x})} \frac{f(\vec{y})}{\sqrt{\theta(\vec{y})}} \max\left(\log\frac{H}{\abs{\vec{x} - \vec{y}}}, 0\right) \, d\vec{y}.
\end{align*}
In our implementation, this continuous integral is approximated by a fixed-stencil convolution on the discretisation grid: the kernel shown in Figure~\ref{fig:diffusion-smoother} is precomputed once from the truncated Green's function evaluated at the grid offsets \(\Delta\vec{x}\) with \(\abs{\Delta\vec{x}} \leq H\); pre- and post-multiplication by \(\theta^{-1/2}\) is implemented as a diagonal scaling. The resulting smoother is therefore a banded matrix.

This smoother is similar to the Green's function-based preconditioner constructed by Ichimura \etal\cite{ichimura2020fast} for equations of diffusion-reaction type. For such equations, the Green's function is more localised around \(\vec{x} \approx \vec{y}\) (see Loghin~\cite{loghin1999greens}), meaning that such a smoother is an effective preconditioner on its own. However, for the diffusion equation the addition of a coarse correction is found to greatly improve the convergence speed. Here, \(H = 0.05\) is used.

The other smoothers are as follows:
\begin{itemize}
	\item Incomplete Cholesky: based on an incomplete Cholesky factorisation with zero fill-in, i.e.~\(\mat{M}^{-1} = \mat{L}^{-\mathsf{T}}\mat{L}^{-1}\) where \(\mat{L}\mat{L}\transpose \approx \mat{A}\) but \(\mat{L}\) is restricted to have the same sparsity pattern as the lower-triangular part of \(\mat{A}\).
	\item Jacobi: standard point Jacobi with underrelaxation parameter \(\omega = \frac{4}{5}\), i.e.~\(\mat{M}^{-1} = \frac{4}{5}\text{diag}(\mat{A})^{-1}\).
	\item Block-Jacobi: using 81 non-overlapping blocks based on a subdivision of the domain \((0, 1)^2\) into \(9 \times 9\) non-overlapping subdomains containing approximately \(\frac{n}{9} \times \frac{n}{9}\) unknowns each (where \(n\) is the number of discretisation points on one axis). Therefore, the block-Jacobi preconditioner inverts matrix blocks that are approximately of size \(\frac{N}{81} \times \frac{N}{81}\) where \(N\) is the total number of unknowns of the system. Again, an underrelaxation parameter \(\omega = \frac{4}{5}\) is used.
\end{itemize}

\begin{figure}
	\centering
	\includegraphics[width=0.7\textwidth]{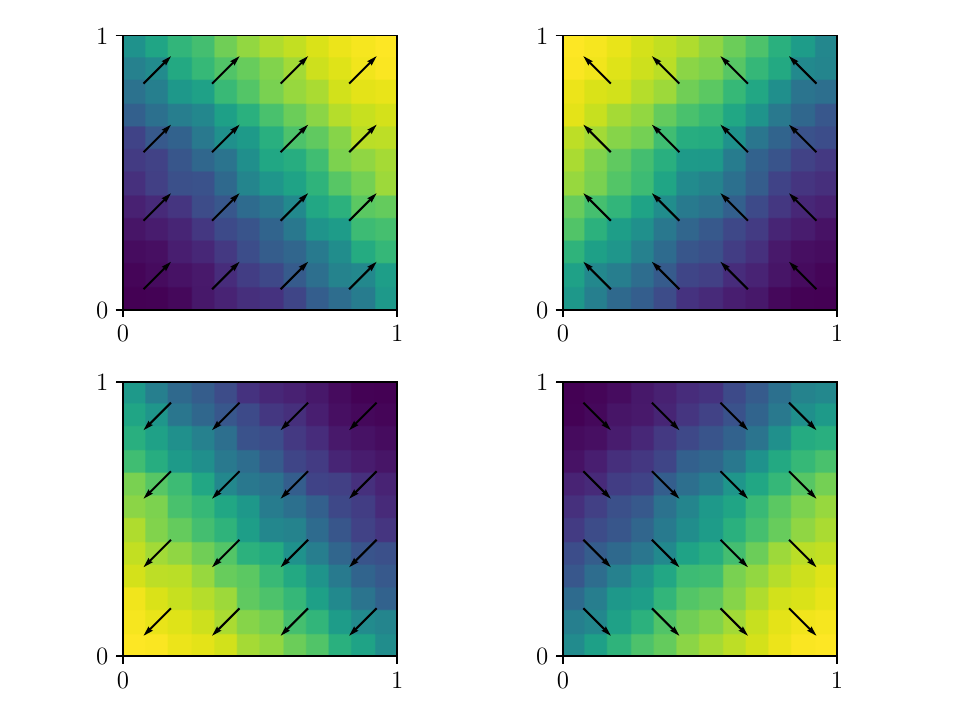}
	\caption{The four sweep directions used by the four-direction Gauss-Seidel smoother.}
	\label{fig:advdiff-smoother}
\end{figure}

\paragraph*{Advection-diffusion}
In the advection-diffusion equation, a four-direction Gauss-Seidel smoother is used (see Section~4.3 of \cite{wesseling2004introduction}). This preconditioner is a multiplicative combination of four separate point-Gauss-Seidel sweeps, each using a different permutation of the unknowns, as shown in Figure~\ref{fig:advdiff-smoother}. This smoother therefore depends algebraically on the matrix \(\mat{A}\) but otherwise does not depend on the velocity field, and is also applicable to advection-diffusion equations with a spatially varying velocity field.

\paragraph*{Helmholtz}
Finally, the smoother for the Helmholtz equation is chosen to be the same as that of the diffusion equation, i.e.~based on the free-space Green's function for the diffusion equation. While it may seem more reasonable to use the Green's function for the Helmholtz equation, this was found to work poorly as a smoother.

\FloatBarrier
\section{Model architecture and training procedure}
\label{app:model-architecture-training}
\begin{figure}
	\centering
	\begin{subfigure}{0.7\textwidth}
		\includegraphics[width=\textwidth]{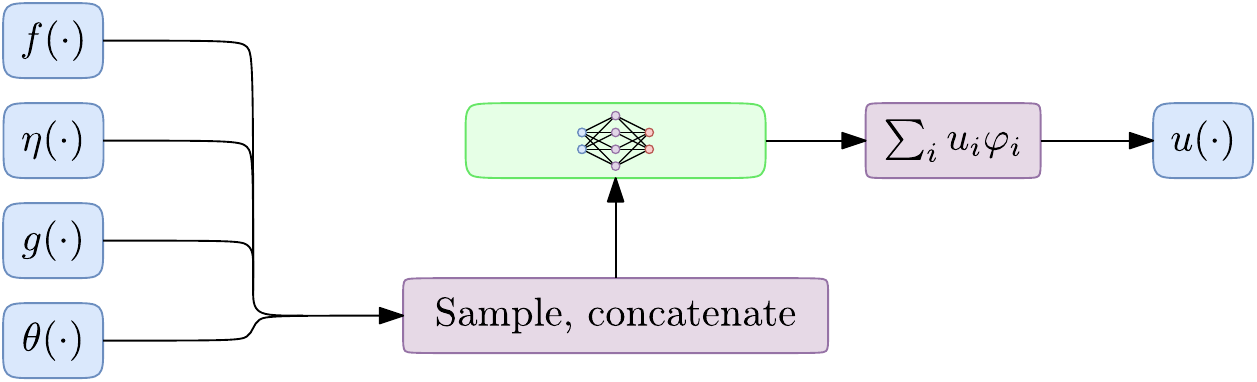}
		\caption{DeepONet}
	\end{subfigure}
	\begin{subfigure}{0.7\textwidth}
		\includegraphics[width=\textwidth]{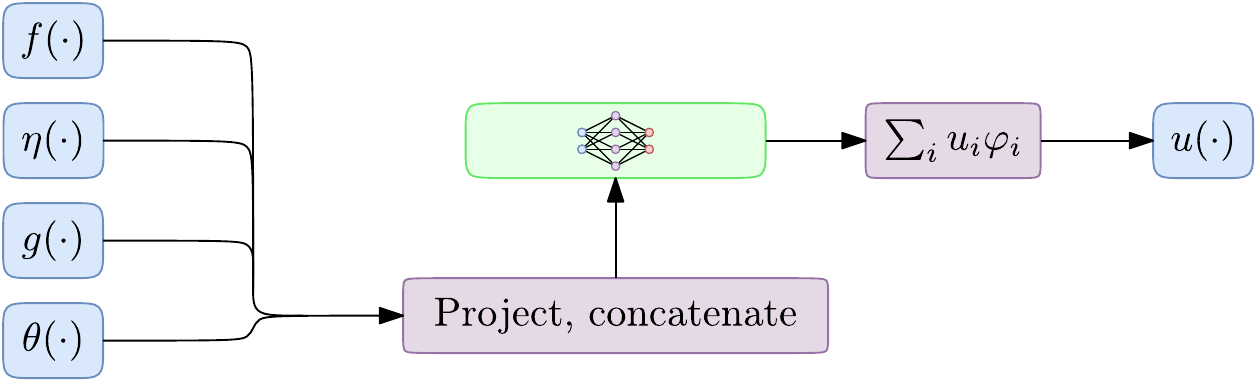}
		\caption{RINO}
	\end{subfigure}
	\begin{subfigure}{0.7\textwidth}
		\includegraphics[width=\textwidth]{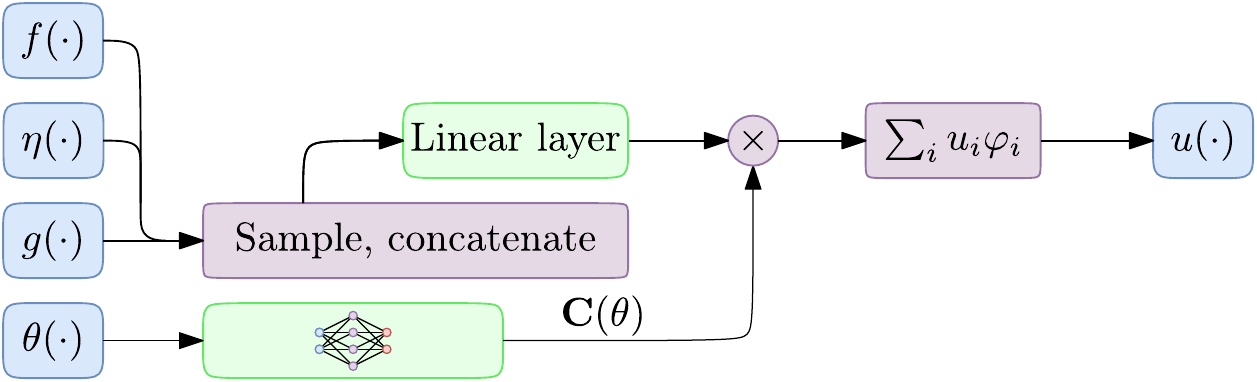}
		\caption{VarMiON}
	\end{subfigure}
	\begin{subfigure}{0.7\textwidth}
		\includegraphics[width=\textwidth]{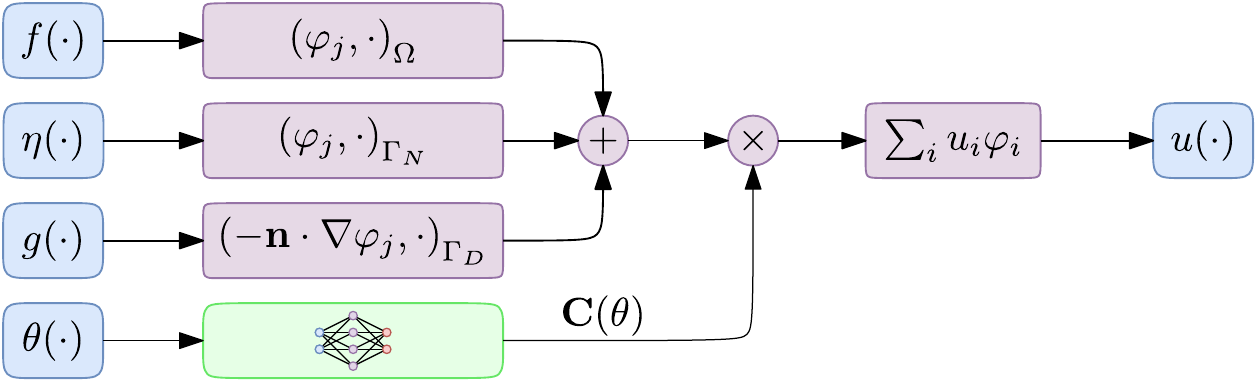}
		\caption{NGO}
	\end{subfigure}
	\caption{Architectures of the four neural operators considered in this work. The order matches the design grid of Table~\ref{tab:no-design-axes}: DeepONet (sample-and-non-linear), RINO (integrate-and-non-linear), VarMiON (sample-and-linear), NGO (integrate-and-linear).}
	\label{fig:model-architectures}
\end{figure}

Figure~\ref{fig:model-architectures} shows the architecture of the four neural operators tested. All symbols appearing in the figure are as defined in Section~\ref{subsec:four-preconditioners}. Note that the NGO architecture used here is referred to as ``Data-NGO'' in \cite{melchers2026neural}. An important detail is that the DeepONet and RINO use the exact same architecture, the only difference being whether the input fields are sampled (DeepONet) or projected onto the basis (RINO). This means that the consistently poor performance of the DeepONet-based preconditioners is not due to insufficiencies in the parameter count, number of hidden layers, or other architectural properties: the RINO, which uses the exact same neural network architecture, is typically much better as a preconditioner.

To help the DeepONet and RINO, i.e.~the non-linear models, to generalise to right-hand-side vectors outside of the training data distribution, these models scale their inputs so as to ensure positive homogeneity in the source terms. Concretely, the source terms \(f\), \(g\), and \(\eta\) are uniformly scaled such that their maximal \(L^2\)-norm equals unity. This scaling is then undone on the approximate solution obtained by the model:
\begin{align*}
	\alpha  & = \max\left( \norm{f}_{L^2(\Omega)}, \norm{g}_{\Gamma_D}, \norm{\eta}_{\Gamma_N} \right),        \\
	\hat{f} & = \frac{1}{\alpha}f, \quad\hat{g} = \frac{1}{\alpha}{g}, \quad\hat{\eta} = \frac{1}{\alpha}\eta, \\
	\hat{u} & = \text{model}(\theta[, \vec{c}], \hat{f}, \hat{g}, \hat{\eta}),                                 \\
	u       & = \alpha \hat{u}.
\end{align*}
This scaling is especially important as the source terms corresponding to vectors produced by Krylov methods may be scaled differently than those in the training data set. For the diffusion and advection-diffusion equations, a similar scaling is performed in the PDE coefficients \(\theta\) (and \(\vec{c}\) if applicable). Furthermore, all test problems are symmetric with respect to mirroring in the \(x\)- and \(y\)-directions: mirroring the input fields of the PDE horizontally or vertically (including negating the \(x\)- or \(y\)-component of \(\vec{c}\)) results in a new PDE the solution of which is a mirrored version of the original solution. This symmetry is embedded into each model by passing all four mirrored versions of each PDE problem to the model, `un-mirroring' the resulting predictions, and averaging.

In all cases, the neural operators are trained to minimise the relative squared \(L^2\) error of the PDE solutions:
\begin{align*}
	\text{Loss} = \frac{1}{N}\sum_{i=1}^{N} \frac{\norm{\text{NN}(\theta_i[, \vec{c}_i], f_i, g_i, \eta_i) - u_i^h}_{L^2}^2}{\norm{u_i^h}_{L^2}^2},
\end{align*}
where \(u_i^h\) is the high-fidelity PDE solution obtained by the finite element discretisation of the \(i\)-th problem.

Table~\ref{tab:parameter-counts} shows the parameter counts of the compared models for each of the test equations. Since the DeepONet and RINO use the exact same architecture and therefore have the same parameter counts by construction, their parameter counts are shown together.

\begin{table}
	\centering
	\caption{Parameter counts of the compared models for each of the test equations. Unless otherwise specified, these parameter counts refer to the 2D test problems with the \(12 \times 12\) basis.}
	\label{tab:parameter-counts}
	\begin{tabular}{l|rrr}
		\toprule
		Equation                              & DeepONet/RINO & VarMiON & NGO   \\
		\midrule
		Diffusion                             & 32664         & 35365   & 30409 \\
		\quad \(m = 6 \times 6\)              &               &         & 26881 \\
		\quad \(m = 24 \times 24\)            &               &         & 35581 \\
		\quad 3D, \(m = 8 \times 8 \times 8\) &               &         & 70832 \\
		Advection-diffusion                   & 32784         & 35389   & 30457 \\
		Helmholtz                             & 33336         & 35414   & 30506 \\
		\bottomrule
	\end{tabular}
\end{table}

Table~\ref{tab:training-config} gives other details about the training procedure. The validation data is used to prevent overfitting: after the 20000 epochs have passed, the trainable parameters of the model are restored to the state that achieved the lowest error on the \emph{validation data}, as opposed to the \emph{training data} over which the error is minimised with ADAM. Therefore, if a model overfits, its error on the validation data set will begin to increase, meaning that after training the model is restored to a state before it began overfitting.
\begin{table}
	\centering
	\caption{An overview of the other training parameters used in this work}
	\label{tab:training-config}
	\begin{tabular}{r l}
		\toprule
		Loss function   & Relative squared \(L^2\) error \\
		Optimiser       & ADAM                           \\
		Learning rate   & Default (\(10^{-3}\))          \\
		Batch size      & 100                            \\
		Training time   & 20000 epochs                   \\
		Training data   & 8000 samples                   \\
		Validation data & 1000 samples                   \\
		Testing data    & 1000 samples                   \\
		\bottomrule
	\end{tabular}
\end{table}

\FloatBarrier
\section{Results in table form}
\label{app:results-table}
Table~\ref{tab:all-results} shows, for each of the three 2D PDEs, the main results also shown in Figures~\ref{fig:diff-prec-bars},~\ref{fig:advdiff-prec-bars}, and \ref{fig:helmholtz-prec-bars}. In this table, the average iteration counts include only the problems in which F-GMRES converged within 5000 iterations. The numbers in parentheses count how many problems were not solved within 5000 iterations.
\begin{table}
	\caption{
		Mean F-GMRES iterations required for various models to solve diffusion (Dataset 3), advection-diffusion, and Helmholtz equations. The numbers in parentheses list the number of problems (out of 1000) on which F-GMRES did not converge in 5000 iterations.
	}
	\label{tab:all-results}
	\resizebox{\columnwidth}{!}{
		\begin{tabular}{l|l | rr|rr|rr|rr|rr}
		    \toprule
		    & \multirow{2}*{Preconditioner} & \multicolumn{10}{c}{F-GMRES iters (failures)} \\
		    && \multicolumn{2}{c|}{\(n = 20\)} & \multicolumn{2}{c|}{\(n = 40\)} & \multicolumn{2}{c|}{\(n = 80\)} & \multicolumn{2}{c|}{\(n = 160\)} & \multicolumn{2}{c}{\(n = 320\)} \\
		    \midrule
			\multirow{7}*{\rotatebox{90}{Diffusion}}
		    & None          &   61.0 &    (0) &  248.5 &    (0) & 1105.4 &    (0) & 4247.3 &  (490) &      - & (1000) \\
		    & Smoother only &   25.4 &    (0) &   43.7 &    (0) &   69.9 &    (0) &  111.7 &    (0) &  181.0 &    (0) \\
		    & Exact         &    4.1 &    (0) &    6.4 &    (0) &    9.8 &    (0) &   14.9 &    (0) &   23.8 &    (0) \\
		    & DeepONet      &   50.7 &    (0) &   19.7 &    (0) &   23.4 &    (0) &   41.0 &    (0) &      - & (1000) \\
		    & RINO          &    6.3 &    (0) &    7.1 &    (0) &    9.9 &    (0) &   16.0 &    (0) &   26.4 &    (0) \\
		    & VarMiON       &    8.2 &    (0) &    8.0 &    (0) &   13.5 &    (0) &   35.0 &    (0) &      - & (1000) \\
		    & NGO           &    5.0 &    (0) &    6.8 &    (0) &    9.9 &    (0) &   14.9 &    (0) &   23.7 &    (0) \\
		    \midrule
			\multirow{7}*{\rotatebox{90}{Advection-diffusion}}
		    & None          &   44.0 &    (0) &  117.4 &    (0) &  258.2 &    (0) &  581.0 &    (0) & 1496.8 &   (37) \\
		    & Smoother only &    5.7 &    (0) &   14.7 &    (0) &   47.0 &    (0) &  124.7 &    (0) &  246.6 &    (0) \\
		    & Exact         &    2.0 &    (0) &    5.0 &    (0) &   12.8 &    (0) &   32.0 &    (0) &   74.3 &    (0) \\
		    & DeepONet      &   59.6 &    (0) & 1515.8 &  (238) &      - & (1000) &      - & (1000) &      - & (1000) \\
		    & RINO          &    9.0 &    (0) &   29.6 &    (0) &   50.3 &    (0) &  127.0 &    (0) &  416.2 &    (0) \\
		    & VarMiON       &    3.7 &    (0) &    7.7 &    (0) &   25.8 &    (0) &  135.9 &    (0) &  689.7 &  (500) \\
		    & NGO           &    2.3 &    (0) &    5.1 &    (0) &   13.9 &    (0) &   35.9 &    (0) &   88.9 &    (0) \\
		    \midrule
			\multirow{7}*{\rotatebox{90}{Helmholtz}}
		    & None          &  330.7 &    (0) & 1119.4 &    (1) &      - & (1000) &      - & (1000) &      - & (1000) \\
		    & Smoother only &  168.7 &    (0) &  297.9 &    (0) &  481.8 &    (0) &  732.2 &    (0) & 1219.1 &    (2) \\
		    & Exact         &    5.4 &    (0) &    7.7 &    (0) &   10.9 &    (0) &   15.7 &    (0) &   24.6 &    (0) \\
		    & DeepONet      &      - & (1000) &      - & (1000) &      - & (1000) & 3580.0 &  (998) &      - & (1000) \\
		    & RINO          &  355.3 &  (737) &  586.0 &  (437) &  216.4 &  (323) &  447.8 &  (230) &  710.2 &  (144) \\
		    & VarMiON       &   34.3 &   (25) &   57.7 &   (21) &  127.9 &   (31) &  255.0 &  (203) &      - & (1000) \\
		    & NGO           &   51.8 &   (39) &   35.6 &   (17) &   37.5 &   (15) &   45.5 &   (22) &   60.2 &   (32) \\
		    \bottomrule
		\end{tabular}
	}
\end{table}

\end{document}